\documentclass{article}%
\usepackage{amssymb}
\usepackage{amsfonts}
\usepackage{amsmath}
\usepackage{graphicx}%
\providecommand{\U}[1]{\protect\rule{.1in}{.1in}}
\newtheorem{theorem}{Theorem}

\newtheorem{corollary}[theorem]{Corollary}

\newtheorem{definition}[theorem]{Definition}

\newtheorem{lemma}[theorem]{Lemma}
\newtheorem{notation}[theorem]{Notation}

\newtheorem{proposition}[theorem]{Proposition}
\newtheorem{remark}[theorem]{Remark}

\begin{document}

\title{Quantum Background Independence and Witten Geometric Quantization of the
Moduli of CY\ Threefolds.}
\author{Andrey Todorov\\University of California, \\Department of Mathematics \\Santa Cruz, CA 95064\\Bulgarian Academy of Sciences\\Institute of Mathematics\\Sofia, Bulgaria\\\ \ \ \ \ \ \ \ \ \ \ \ \ \ \ \ \ \ \ \ \ \ \ \ \ \ \ \ \ \ \ \ \ \emph{Dedicated
to Betty (1949-2002)} }
\maketitle

\begin{abstract}
In this paper we study two different topics. The first topic is the
applications of the geometric quantization scheme of Witten introduced in
\cite{ADW} and \cite{EMSS} to the problem of the quantum background
independence in string theory. The second topic is the introduction of a
$\mathbb{Z}$ structure on the tangent space of the moduli space of polarized
CY threefolds $\mathcal{M(}$M). Based on the existence of a $\mathbb{Z}$
structure on the tangent space of the moduli space of polarized CY threefolds
we associate an algebraic integrable structure on the tangent bundle of
$\mathcal{M(}$M). In both cases it is crucial to construct a flat
$Sp(2h^{2,1},\mathbb{R)}$ connection on the tangent bundle of the moduli space
$\mathcal{M}$(M) of polarized CY threefolds. In this paper we define a Higgs
field on the tangent bundle of the moduli space of CY threefolds. Combining
this Higgs field with the Levi-Cevita connection of the Weil-Petersson metrics
on the moduli space of three dimensional CY manifolds, we construct a
new\ \ $Sp(2h^{2,1},\mathbb{R)}$ connection, following the ideas of Cecotti
and Vafa. Using this flat connection, we apply the scheme of geometric
quantization introduced by Axelrod, Della Pietra and Witten to the tangent
bundle of the moduli space of three dimensional CY manifolds to realize Witten
program in \cite{W93} of solving the problem of background quantum
independence for topological string field theories. By modifying the
calculations of E. Witten done on the flat bundle $R^{3}\pi_{\ast}\mathbb{C}$
to the tangent bundle of the moduli space of CY threefolds, we derive the
holomorphic anomaly equations of Bershadsky, Cecotti, Ooguri and Vafa as flat
projective connection.

\end{abstract}
\tableofcontents

\section{Introduction}

By definition a Calabi-Yau (CY) manifold is a compact complex $n-$dimensional
K\"{a}hler manifold M with a holomorphic $n-$form $\Omega_{\text{M}}$ which
has no zeroes and $H^{0}\left(  \text{M},\Omega_{\text{M}}^{k}\right)  =0$ for
$0<k<n$. Calabi-Yau manifolds are playing important role in string theory. The
powerful ideas from string theory played a very important role in the recent
developments in some branches of mathematics and especially in the study of
moduli of CY manifolds. In this paper we will study moduli space of CY
threefolds based on the ideas introduced in \cite{CV}, \cite{BCOV} and
\cite{W93}.

In \cite{to89} and \cite{tian} it was proved that there are no obstructions to
the deformations of the complex structures on CY manifolds. This means that
the local moduli space of CY manifolds is smooth of dimension%
\[
h^{n-1,1}=\dim_{\mathbb{C}}H^{1}\left(  \text{M},\Omega_{\text{M}}%
^{n-1}\right)  .
\]
From the theory of moduli of polarized algebraic manifolds developed by
Viehweg in \cite{vw} it follows that the moduli space of polarized CY
manifolds is a quasi-projective variety.

The moduli space $\mathcal{M}$ of three dimensional CY manifolds has a very
rich structure. According to the theory of variations of Hodge structures
there exists is a well defined map from the moduli space of polarized CY
manifolds to $\mathbb{P}\left(  H^{n}\left(  \text{M},\mathbb{Z}\right)
\otimes\mathbb{C}\right)  $ which is called the period map$.$ It assigns to
each point $\tau$ of the moduli space the line in $H^{n}\left(  \text{M}%
,\mathbb{Z}\right)  \otimes\mathbb{C)}$ spanned by the cohomology class
represented by the non-zero holomorphic n-form. According to local Torelli
Theorem the period map is a local isomorphism. Local Torelli Theorem implies
that locally the moduli space of CY manifolds can be embedded in
$\mathbb{P}\left(  H^{n}\left(  \text{M},\mathbb{Z}\right)  \otimes
\mathbb{C}\right)  .$ When the dimension $n$ of the CY manifold is odd,
Griffiths and Bryant\ noticed that the intersection form on $H^{n}%
($M,$\mathbb{Z})$ defines on the standard charts $U_{i}=\mathbb{C}^{n-1}$ of
$\mathbb{P}\left(  H^{n}\left(  \text{M},\mathbb{Z}\right)  \otimes
\mathbb{C}\right)  $ holomorphic one forms $\alpha_{i}$ such that $d\alpha
_{i}$ is a skew symmetric form of maximal rank on $\mathbb{C}^{n-1}.$ This
means that on $U_{i}=\mathbb{C}^{n-1}$\ a natural contact structure is
defined. In the case of three dimensional CY manifolds Griffiths and Bryant
proved that the restrictions of $d\alpha_{i}$ on the tangent space of the
image of the local moduli space of CY manifolds is zero. Thus the image of the
local moduli space is a Legandre submanifold. See \cite{BG}. Arnold described
the local structure of the Legandre submanifolds in a contact manifold in
\cite{Ar}. This description implies the existence of a generating holomorphic
function for the local moduli space.

Based on the work \cite{BG} A. Strominger noticed that the potential of the
Weil-Petersson metric on the local moduli space can be expressed through the
generating holomorphic function. See \cite{str}. By using this observation,
Strominger introduced the notion of special K\"{a}hler geometry. V. Cortes
showed that on the tangent bundle of the special K\"{a}hler manifold one can
introduce a Hyper-K\"{a}hler structure. See \cite{cor} and \cite{f}. From here
it follows that on the tangent bundle of $\mathcal{M}\left(  \text{M}\right)
$ one can introduce a Hyper-K\"{a}hler structure. Earlier R. Donagi and E.
Markman constructed in \cite{DM} an analytically completely integrable
Hamiltonian system which is canonically associated with the family of CY
manifolds over the relative dualizing line bundle over the moduli space
$\mathcal{M}\left(  \text{M}\right)  $. They showed that the space of the
Griffiths intermediate Jacobians, associated with the family of three
dimensional CY manifolds on $\mathcal{M}$ carries a Hyper-K\"{a}hler
structure. B. Dubrovin introduced the notion of \ Frobenius manifolds in
\cite{D}. The relations of the structure of Frobenius manifolds and
Gromov-Witten invariants were studied by Yu. I. Manin and M. Kontsevich in
\cite{M}.

The importance of all these structures is justified by the work of Candelas
and coauthors in their seminal paper \cite{CO}. In this paper Candelas and his
coauthors gave an explicit formula for the number of rational curves on the
quintic hypersurface in the four dimensional projective space. M. Kontsevich
defined the correct compactification of the stable maps and realized that one
can use the localization formula for the computations of the rational curves.
See \cite{KO}. Recently B. Lian, K. Liu and Yau gave a rigorous mathematical
proof of the Candelas formula in \cite{LLY}. See also the important paper of
Givental \cite{Giv}.

In this paper we study two different topics. The first topic is the
applications of the geometric quantization scheme of Witten introduced in
\cite{ADW} to the problem of the quantum background independence in string
theory. The second topic is the introduction of a $\mathbb{Z}$ structure on
the tangent space of the moduli space of polarized CY threefolds
$\mathcal{M}\left(  \text{M}\right)  $ and thus we associate an algebraic
integrable structure on the tangent bundle of $\mathcal{M}\left(
\text{M}\right)  $. For both topics it is crucial to construct a flat
$\mathbb{S}p(2h^{2,1},\mathbb{R)}$ connection on the tangent bundle of the
moduli space $\mathcal{M}\left(  \text{M}\right)  $ of polarized CY threefolds.

The problem of the quantum background independence was addressed in
\cite{W93}. In the paper \cite{W93} Witten wrote: \ 

"Finding the right framework for intrinsic, background independent formulation
of string theory is one of the main problems in the subject, and so far has
remained out of reach..."

In fact in \cite{W93} a program was outlined how one can solve the problem of
the background independence in the topological field theory:

"Though the interpretation of the holomorphic anomaly as an obstruction to
background independence eliminates some thorny puzzles, it is not satisfactory
to simply leave matters as this. Is there some sophisticated sense in which
background independence does hold? In thinking about this question, it is
natural to examine the all orders generalization of the holomorphic anomaly
equation, which in the final equation of their paper \cite{BCOV}) Besrshadsky
et. al. write the following form. Let $F_{g}$ be the genus $g$ free energy.
Then
\begin{equation}
\overline{\partial_{i^{^{\prime}}}}F_{g}=\overline{C}_{i^{^{\prime}%
}j^{^{\prime}}k^{^{\prime}}}e^{2K}G^{jj^{^{\prime}}}G^{kk^{^{\prime}}}\left(
D_{j}D_{k}F_{g-1}+\frac{1}{2}%
{\displaystyle\sum\limits_{r}}
D_{j}F_{r}\cdot D_{k}F_{g-r}\right)  . \label{Z0}%
\end{equation}
This equation can be written as a linear equation for
\begin{equation}
Z=\exp\left(  \frac{1}{2}%
{\displaystyle\sum\limits_{g=0}^{\infty}}
\lambda^{2g-2}F_{g}\right)  , \label{Z}%
\end{equation}
namely%
\begin{equation}
\left(  \overline{\partial_{i^{^{\prime}}}}-\lambda^{2}\overline
{C}_{i^{^{\prime}}j^{^{\prime}}k^{^{\prime}}}e^{2K}G^{jj^{^{\prime}}%
}G^{kk^{^{\prime}}}D_{j}D_{k}\right)  Z=0 \label{Mase}%
\end{equation}
This linear equation is called a master equation by Bershadsky et. al.; it is
similar to the structure of the heat equations obeyed by theta functions...

It would be nice to interpret $\left(  \ref{Mase}\right)  $ as a statement of
some sophisticated version of background independence. In thinking about this
equation, a natural analogy arises with Chern-Simon gauge theory in $2+1$
dimensions. In this theory, an initial value surface is a Riemann surface
$\Sigma.$ In the Hamiltonian formulation of the theory, one construct a
Hilbert space $\mathcal{H}$ upon quantization on $\Sigma.\mathcal{H}$ should
be obtained from some physical space $\mathbf{W}$ (a moduli space of flat
connections on $\Sigma$). Because the underlying Chern-Simon Lagrangian does
not depend on the choice of the metric, one would like to construct
$\mathcal{H}$ in a natural, background independent way. In practice, however,
quantization of $\mathbf{W}\,$\ requires a choice of polarization, and there
is no natural way or background independent choice of polarization.

The best that one can do is to pick a complex structure $J$ on $\Sigma,$
whereupon $\mathbf{W}$ gets a complex structure. Then a Hilbert space
$\mathcal{H}_{J}$ is constructed as a suitable space of holomorphic functions
(really sections of a line bundle) over $\mathbf{W}$. We denote such function
as $\psi(a^{i},t^{\prime a})$ where $a^{i}$ are complex coordinates on
$\mathbf{W}$ and $t^{\prime a}$ are coordinates parametrizing the choice of
$J.$ Now background independence does not hold in a naive sense; $\psi$ can
not be independent of $t^{^{\prime i}}$ (given that it is to be holomorphic on
$\mathbf{W}$ in a complex structure dependent on $t^{\prime a}).$ But there is
a more sophisticated sense in which background independence can be formulated.
See \cite{ADW} and \cite{EMSS}. The $\mathcal{H}_{J}$ can be identified with
each other (projectively) using a (projectively) flat connection over the
space of $J^{\prime}$s. This connection $\nabla$ is such that a covariant
constant wave function should have the following property: as $J$ changes,
$\psi$ should change by Bogoliubov transformation, representing the effect of
a change in the representation used for the canonical commutation relations.
Using parallel transport by $\nabla$ to identify the various $\mathcal{H}%
_{J}^{\prime}$ s are realizations determined by a $J$-dependent choice of the
representation of the canonical commutators. Background independence of
$\psi(a^{i},t^{\prime a})$ should be interpreted to mean that the quantum
state represented by $\psi$ is independent of $t^{\prime a},$ or equivalently
that $\psi$ is invariant under parallel transport by $\nabla.$ Concretely,
this can be written as an equation:
\begin{equation}
\left(  \frac{\partial}{\partial t^{\prime a}}-\frac{1}{4}\left(
\frac{\partial J}{\partial t^{\prime a}}\omega^{-1}\right)  ^{ij}\frac
{D}{Da^{i}}\frac{D}{Da^{j}}\right)  \psi=0. \label{Mase1}%
\end{equation}
that is analogous to the heat equation for theta functions..."

In \cite{W93} the above program is realized on the space $\mathbf{W}%
=H^{3}\left(  \text{M},\mathbb{R}\right)  $. Bershadsky, Cecotti, Ooguri and
Vafa work on $H^{1}($M,$T_{\text{M}}^{1,0}),$ i.e. the tangent space to the
moduli of CY$.$ The space $\mathbf{W}=H^{3}\left(  \text{M},\mathbb{R}\right)
$ has a natural symplectic form structure given by the intersection pairing%
\[
\omega(\alpha,\beta):=%
{\displaystyle\int\limits_{\text{M}}}
\alpha\wedge\beta.
\]
The complex structure on M defines a complex structure on $H^{3}\left(
\text{M},\mathbb{R}\right)  .$ On the vector bundle $R^{3}\pi_{\ast}%
\mathbb{R}$ over the moduli space with a fibre $\mathbf{W}=H^{3}\left(
\text{M},\mathbb{R}\right)  $ we have a natural flat $\mathbb{S}%
p(2h^{2,1}+2,\mathbb{R)}$ connection. The tangent bundle to $\omega
_{\mathcal{X}\text{/}\mathcal{M}\text{ }\left(  \text{M}\right)  \text{ }}$ is
naturally isomorphic to $\pi^{\ast}\left(  R^{3}\pi_{\ast}\mathbb{R}\right)
$. Thus on it we have a natural flat $\mathbb{S}p(2h^{2,1}+2,\mathbb{R)}$ connection.

In the present paper the program of Witten is realized for the tangent bundle
of the moduli space of polarized CY threefolds. One of the most important
ingredient in the realization of the Witten program is the construction of a
flat $\mathbb{S}p(2h^{2,1},\mathbb{R)}$ connection on the tangent bundle of
the moduli space $\mathcal{M}\left(  \text{M}\right)  $ of polarized CY
threefolds. In the present article we constructed such flat $\mathbb{S}%
p(2h^{2,1},\mathbb{R)}$ connection.

The idea of the construction of the\ flat $\mathbb{S}p(2h^{2,1},\mathbb{R)}$
connection on the tangent bundle of the moduli space $\mathcal{M}\left(
\text{M}\right)  $ is to modify the unitary connection of the Weil-Petersson
metric on $\mathcal{M}\left(  \text{M}\right)  $ with a Higgs field to a
$\mathbb{S}p(2h^{2,1},\mathbb{R)}$ connection and then prove that the
$\mathbb{S}p(2h^{2,1},\mathbb{R)}$ connection is flat. The construction of the
Higgs field defined on $R^{1}\pi_{\ast}\Omega_{\mathcal{Y}\left(
\text{M}\right)  \left/  \mathcal{M}\left(  \text{M}\right)  \right.  }^{2}$
is done by using the cup product $\phi_{1}\wedge\phi_{2}\in$ $H^{2}($%
M,$\wedge^{2}T^{1,0})$ for $\phi_{i}\in H^{1}($M,$T^{1,0}),$ the
identifications of $H^{1}\left(  \text{M},\Omega_{\text{M}}^{2}\right)  $ with
$H^{1}($M,$T^{1,0})$, $H^{2}\left(  \text{M},\Omega_{\text{M}}^{1}\right)  $
with $H^{2}($M,$\wedge^{2}T^{1,0})$ and the identification of $H^{1}\left(
\text{M},\Omega_{\text{M}}^{2}\right)  $ with $H^{2}\left(  \text{M}%
,\Omega_{\text{M}}^{1}\right)  $ by the Poincare duality.

The construction of a flat $\mathbb{S}p(2h^{2,1},\mathbb{R)}$ connection on
the tangent bundle of $\mathcal{M}\left(  \text{M}\right)  $ is related to the
important example of \ a special Hyper-K\"{a}hler manifold which occurs in
four dimensional gauge theories \ with $N=2$ supersymmetry: the scalars in the
vector multiplet lie in a special K\"{a}hler manifolds. The moduli space of
such theories was studied by Cecotti and Vafa in \cite{CV}. They introduced
the tt* equations. One of the observation in this paper is that the analogue
of the tt* equations\ in case of the moduli of polarized CY threefolds is the
same as the Yang-Mills equations coupled with Higgs fields that were studied
by Hitchin in case of Riemann surfaces in \cite{H} and by C. Simpson in
general in \cite{sim2}.

The flat $\mathbb{S}p(2h^{2,1},\mathbb{R)}$ connection is crucial to apply the
geometric quantization method of Witten to the tangent space of the moduli
space $\mathcal{M}\left(  \text{M}\right)  $ of polarized CY threefolds to
solve the problem of the ground quantum independence in the topological field
theory. On the basis of the geometric quantization of the tangent bundle of
$\mathcal{M}\left(  \text{M}\right)  $ we are able to modify the beautiful
computations of E. Witten in \cite{W93} to obtain a projective connection on
some infinite dimensional Hilbert space bundle. We prove that the holomorphic
anomaly equations $\left(  \ref{Z0}\right)  $ of Bershadsky, Cecotti, Ooguri
and Vafa imply that the free energy obtained from the \textquotedblright
counting functions\textquotedblright\ F$_{g}$ of curves of genus g on a CY
manifold M is a parallel section of a projective flat connection. Our
computations are based on the technique developed in \cite{to89}.

The projective connection constructed in \cite{W93} is different from ours
since we work on different spaces. The difference appeared in the computation
of the formula for $\left(  dJ\omega^{-1}\right)  .$ On the space
$\mathbf{W}=H^{3}\left(  \text{M},\mathbb{R}\right)  $ Witten obtained that
\[
\left(  dJ\right)  _{a}^{\overline{b}}=2\sum_{c,d}\overline{C}_{acd}%
g^{d,\overline{b}}%
\]
where $\left(  g_{a,\overline{b}}\right)  $ defines the symplectic structure
on $\mathbf{W}$ coming from the cup product and $C_{acd}$ is the Yukawa
coupling. Our formula on $\mathbf{W}=H^{1}\left(  \text{M},\Omega_{\text{M}%
}^{2}\right)  $ is%
\[
\left(  dJ\right)  _{a}^{\overline{b}}=\sum_{c,d}\overline{C}_{acd}%
g^{d,\overline{b}},
\]
where $\left(  g_{a,\overline{b}}\right)  $ is the symplectic form obtained
from the restriction of the cup product on $H^{1}\left(  \text{M}%
,\Omega_{\text{M}}^{2}\right)  \subset H^{3}\left(  \text{M},\mathbb{R}%
\right)  .$ At the end we obtained exactly the formula $\left(  \ref{Mase1}%
\right)  $ suggested by E. Witten.

In \cite{BCOV} two equations are derived. One of them is $\left(
\ref{Z0}\right)  $. It gives a recurrent relation between $F_{g}$ 's$.$ The
other equation in \cite{BCOV} is $\left(  \ref{Zc}\right)  .$ These two
equations are marked as $(3.6)$ and $(3.8)$ in \cite{BCOV}. According to
\cite{BCOV} $\ $the free energy $Z$ satisfy the equation$\ \left(
\ref{Zc}\right)  $. One can notice that there is a difference between our
equation and the equation $\left(  \ref{Zc}\right)  $ for the free energy $Z$
in \cite{BCOV}. The holomorphic anomaly equation $\left(  \ref{Zc}\right)  $
in \cite{BCOV} involves the term $F_{1}$ while ours do not.

It was pointed out in \cite{W93} that the anomaly equations are the analogue
of the heat equations for the classical theta functions. Thus they are of
second order. From here one can deduce that if we know the functions F$_{0}$
and F$_{1}$ that count the rational and elliptic curves on M we will know the
functions F$_{g}$ that count all curves of given genus $g>1.$ It was Welters
who first noticed that the heat equation of theta functions can be interpreted
as a projective connection. See \cite{Wel83}. Later N. Hitchin used the
results of \cite{Wel83} to construct a projectively flat connection on a
vector bundle over the Teichm\"{u}ller space constructed from the symmetric
tensors of stable bundle over a Riemann surface. See \cite{H}. For other
useful applications of the geometric approach to quantization see \cite{EMSS}.

The second problem discussed in this paper is about the existence of
$\mathbb{Z}$ structure on the tangent bundle of the moduli space
$\mathcal{M}\left(  \text{M}\right)  $\ of polarized CY threefolds. This
problem is suggested by the mirror symmetry conjecture since it suggests that
\[
H^{3,0}\oplus H^{2,1}\oplus H^{1,2}\oplus H^{0,3}%
\]
can be "identified" on the mirror side with
\[
H^{0}\oplus H^{2}\oplus H^{4}\oplus H^{6}.
\]
Thus since by the mirror conjecture $H^{2,1}$ can be identified with $H^{2}$
one should expect some natural $\mathbb{Z}$ structure on $H^{2,1}$ invariant
under the flat $\mathbb{S}p(2h^{2,1},\mathbb{R)}$ connection. Thus we need to
define at a fixed point of the moduli space $\tau_{0}\in\mathcal{M}\left(
\text{M}\right)  $ a $\mathbb{Z}$ structure on the tangent space $T_{\tau
_{0},\mathcal{M}\left(  \text{M}\right)  }=H^{1}($M$_{\tau_{0}},\Omega
_{\text{M}_{\tau_{0}}}^{2})$. One way to obtain a natural $\mathbb{Z}$
structure on $R^{1}\pi_{\ast}\Omega_{\mathcal{Y}\left(  \text{M}\right)
\left/  \mathcal{M}\left(  \text{M}\right)  \right.  }^{2}$ is the following
one. Suppose that there exists a point $\tau_{0}\in\mathcal{M}\left(
\text{M}\right)  $ such that
\begin{equation}
H^{3,0}\left(  \text{M}_{\tau_{0}}\right)  \oplus H^{0,3}\left(
\text{M}_{\tau_{0}}\right)  =\Lambda_{0}\otimes\mathbb{C}, \label{ic}%
\end{equation}
where $\Lambda_{0}$ is a rank two sublattice in $H^{3}($M$_{\tau_{0}%
},\mathbb{Z}).$ Once we construct such $\mathbb{Z}$ structure on $T_{\tau
_{0},\mathcal{M}\left(  \text{M}\right)  },$ we can use the parallel transport
to define a $\mathbb{Z}$ structure on each tangent space of $\mathcal{M}%
\left(  \text{M}\right)  $. Unfortunately the existence of points that
$\left(  \ref{ic}\right)  $ is satisfied is a very rare phenomenon for CY
manifolds. There is a conjecture due to B. Mazur and Y. Andre which states
that if the moduli space of CY manifold is a Shimura variety then such points
are everywhere dense subset. The moduli space of polarized CY manifolds that
are not locally symmetric spaces probably will not contain everywhere dense
subset of points that correspond to CY manifolds for which
\[
H^{3,0}\left(  \text{M}_{\tau_{0}}\right)  \oplus\overline{H^{3,0}\left(
\text{M}_{\tau_{0}}\right)  }=\Lambda_{0}\otimes\mathbb{C},
\]
where $\Lambda_{0}$ is a rank two sublattice in $H^{3}\left(  \text{M}%
_{\tau_{0}},\mathbb{Z}\right)  .$

The idea of the introduction of the $\mathbb{Z}$ structure on the tangent
space of the moduli space $\mathcal{M}\left(  \text{M}\right)  $ is to
consider the deformation space of M$\times\overline{\text{M}}.$ We relate the
local deformation space on M$\times\overline{\text{M}}$ to the variation of
Hodge structure of weight two with $p_{g}=1$. Such Variations of Hodge
structures for the products M$\times\overline{\text{M}}$ are defined by the
two dimensional real subspace $H^{3,0}\left(  \text{M}\right)  \oplus
\overline{H^{3,0}\left(  \text{M}\right)  }$ in $H^{3}\left(  \text{M}%
,\mathbb{R}\right)  $ is generated by $\operatorname{Re}\Omega_{\tau}$ and
$\operatorname{Im}\Omega_{\tau}$ and they are parametrized by the symmetric
space
\[
\mathbb{SO}_{0}(2,2h^{2,1})\left/  \mathbb{SO}(2)\times\mathbb{SO}%
(2h^{2,1})\right.
\]
where the set of points for which $\left(  \ref{ic}\right)  $ holds is an
everywhere dense subset$.$ Thus we are in situation similar to the moduli of
algebraic polarized K3 surfaces. For the $(\tau,\overline{\nu})$ in the local
moduli space of M$\times\overline{\text{M}}$ that corresponds to M$_{\tau
}\times\overline{\text{M}_{\nu}}$ we construct a Hodge structure of weight
two
\[
H_{\tau,\nu}^{2,0}\oplus H_{\tau,\nu}^{1,1}\oplus\overline{H_{\tau,\nu}^{2,0}}%
\]
where $H_{\tau,\nu}^{2,0}\oplus\overline{H_{\tau,\nu}^{2,0}}$ is the two
dimensional subspace in $H^{3}\left(  \text{M},\mathbb{R}\right)  $ generated
by
\[
\operatorname{Re}\left(  \Omega_{\tau_{1}}+\overline{\Omega_{\tau_{2}}%
}\right)  \text{ and }\operatorname{Im}\left(  \Omega_{\tau_{2}}%
-\overline{\Omega_{\tau_{2}}}\right)  .
\]
It is not difficult to show that the points $(\tau,\nu)$ in the local moduli
space of M$\times\overline{\text{M}}$ such that%
\[
H_{\tau,\nu}^{2,0}\oplus\overline{H_{\tau,\nu}^{2,0}}=\Lambda_{1}%
\otimes\mathbb{R},
\]
where $\Lambda_{1}$ is a rank four sublattice in $H^{3}\left(  \text{M}%
,\mathbb{Z}\right)  $ is an everywhere dense subset$.$ Each point $(\tau,\nu)$
of this everywhere dense subset defines a natural $\mathbb{Z}$ structure on
$H_{\tau,\nu}^{1,1}.$ Then by using the flat $\mathbb{S}p(2h^{2,1}%
,\mathbb{R})$ connection on the tangent space of $\mathcal{M}\left(
\text{M}\right)  $ we define a $\mathbb{Z}$ structure on $H_{\tau,\tau}^{2,0}$
and thus on $R^{1}\pi_{\ast}\Omega_{\mathcal{Y}\left(  \text{M}\right)
\left/  \mathcal{M}\left(  \text{M}\right)  \right.  }^{2}.$ By using this
$\mathbb{Z}$ structure we introduce an algebraic integrable structure on the
tangent bundle of $\mathcal{M}\left(  \text{M}\right)  $. In \cite{DM} the
authors introduced algebraic integrable structure on the tangent bundle of the
relative dualizing line bundle of $\omega_{\mathcal{Y}\left(  \text{M}\right)
\left/  \mathcal{M}\left(  \text{M}\right)  \right.  }$ of the moduli space
$\mathcal{M}\left(  \text{M}\right)  $ of polarized CY threefolds.

All the results in the \textbf{Sections 3, 4, 5, 6 }and \textbf{7} are new.
Next we will describe the ideas and the content of each section.

In \textbf{Section 2} we review the results of \cite{to89} and \cite{tian}.
The Teichm\"{u}ller space of the CY\ manifolds is constructed too.

In \textbf{Section 3} we show that the analogue of the tt* equations on
$\mathcal{M}\left(  \text{M}\right)  $ are the same self dual equations that
were studied by N. Hitchin and C. Simpson's in \cite{H87} and \cite{sim1}.
Thus tt* equations define a flat $\mathbb{S}p(2h^{2,1},\mathbb{R)}$ connection
on the bundle $R^{1}\pi_{\ast}\Omega_{\mathcal{Y}\left(  \text{M}\right)
\left/  \mathcal{M}\left(  \text{M}\right)  \right.  \text{ }}^{2}.$ On the
other hand we know that the tangent bundle $\mathcal{T}_{\mathcal{Y}\left(
\text{M}\right)  \left/  \mathcal{M}\left(  \text{M}\right)  \right.  }$ is
isomorphic to $\mathfrak{L}^{\ast}\otimes$ $R^{1}\pi_{\ast}\Omega
_{\mathcal{Y}\left(  \text{M}\right)  \left/  \mathcal{M}\left(
\text{M}\right)  \right.  }^{2}$, where\ $\mathfrak{L}$ is isomorphic to
$\pi_{\ast}\Omega_{\mathcal{Y}\left(  \text{M}\right)  \left/  \mathcal{M}%
\left(  \text{M}\right)  \right.  }^{3}$. We constructed by using the theory
of determinant bundles a holomorphic non-vanishing section $\eta_{\tau}%
\in\Gamma\left(  \mathcal{M}(\text{M}),(\mathfrak{L})\right)  $ in
\cite{BT}$.$ Thus $\eta_{\tau}$ defines a flat structure on the tangent bundle
$\mathcal{T}_{\mathcal{Y}\left(  \text{M}\right)  \left/  \mathcal{M}\left(
\text{M}\right)  \right.  }$ of the moduli space of \ three dimensional CY
manifolds \ $\mathcal{M}\left(  \text{M}\right)  $. Using the flat structure
defined by tt* equations on $R^{1}\pi_{\ast}\Omega_{\mathcal{Y}\left(
\text{M}\right)  \left/  \mathcal{M}\left(  \text{M}\right)  \right.  \text{
}}^{2}$ and the flat structure defined by the section $\eta_{\tau}$ on
$\mathfrak{L}$, we define a flat $\mathbb{S}p(2h^{2,1},\mathbb{R)}$ connection
on the tangent bundle $\mathcal{T}_{\mathcal{Y}\left(  \text{M}\right)
\left/  \mathcal{M}\left(  \text{M}\right)  \right.  }$ of the moduli space
$\mathcal{M}\left(  \text{M}\right)  $ of \ three dimensional CY manifolds. We
will call this connection the Cecotti-Hitchin-Simpson-Vafa connection and will
refer to it as the CHSV connection.

A beautiful theorem of Simpson proved in \cite{sim1} shows when a
quasi-projective variety is covered by a symmetric domain. One can show that
the tt* equations can be interpreted in the same way. This will be done in
\cite{to1}.

In \textbf{Section 3} we interpreted the holomorphic connection which is
defined by the Frobenius Algebra structure on the bundle $R^{1}\pi_{\ast
}\Omega_{\mathcal{Y}\left(  \text{M}\right)  \left/  \mathcal{M}\left(
\text{M}\right)  \right.  }^{2}$ as a Higgs field. It seems that the paper by
Deligne is suggesting that the Higgs field that we constructed is related to
variation of mixed Hodge Structure of CY threefolds, when there exists a
maximal unipotent element in the mapping class group. See \cite{Del}.

In \textbf{Section 4} we review some basic constructions from \cite{ADW}.

In \textbf{Section 5 }we used the ideas from \cite{ADW} and some modifications
of the beautiful computations done by E. Witten in \cite{W93} to quantize the
tangent bundle of $\mathcal{M}\left(  \text{M}\right)  $. This can be done
since we can identify the tangent spaces\ at each point of the moduli space of
CY manifolds by using the parallel transport defined by the flat connection
$\mathbb{S}p(2h^{2,1},\mathbb{R)}$ defined by Cecotti-Hitchin-Simpson-Vafa and
the existence of the non-zero section $\eta_{\tau}$ of the relative dualizing
line bundle over $\mathcal{M}\left(  \text{M}\right)  $. We will show that the
symplectic structure defined by the imaginary part of the Weil-Petersson
metric is parallel with respect to the CHSV connection. In this section we
construct a projective flat connection on some Hilbert vector bundle
associated with the tangent bundle on the moduli space $\mathcal{M}\left(
\text{M}\right)  $. Based on these results, the method from \cite{W93} and the
technique developed in \cite{to89},\ we show that holomorphic anomaly
equations $\left(  \ref{Z0}\right)  $ of Bershadsky, Cecotti, Ooguri and Vafa
of the genus $g\geq2$ imply that the free energy $Z$ defined by $\left(
\ref{Zb}\right)  $ is a parallel with respect to a flat projective connection
constructed in \textbf{Section 6}.

In \textbf{Section 6} we will introduce a natural $\mathbb{Z}$ structure on
the tangent space of each point of the moduli space of CY threefolds by using
the flat $\mathbb{S}p(2h^{2,1},\mathbb{R})$ connection constructed in
\textbf{Section 3}. In order to do that we introduce the notion of the
extended period space of CY threefolds which is similar to the period domain
of marked algebraic polarized K3 surfaces. We know from the moduli theory
of\ algebraic polarized K3 surfaces that the points that define K3 surfaces
with CM structure form an everywhere dense subset. This follows from the fact
that the period domain is an open set on a quadric defined over $\mathbb{Q}$
in the projective space $\mathbb{P}\mathbf{(}\mathbb{Z}^{20}\otimes
\mathbb{C)}$. This fact together with the existence of a flat $\mathbb{S}%
p(4h^{2,1},\mathbb{R)}$ connection on the extended period domain will define
in a natural way a lattice of maximal rank in the tangent space at each point
of the moduli space of CY threefolds. Using the existence of the
$\mathbb{S}p(2h^{2,1},\mathbb{R)}$ connection defined by the tt* equation on
$R^{1}\pi_{\ast}\Omega_{\mathcal{Y}\left(  \text{M}\right)  \left/
\mathcal{M}\left(  \text{M}\right)  \right.  \text{ }}^{2}$, we define by the
parallel translations a $\mathbb{Z}$ structure on the fibres of the bundle
$R^{1}\pi_{\ast}\Omega_{\mathcal{Y}\left(  \text{M}\right)  \left/
\mathcal{M}\left(  \text{M}\right)  \right.  }^{2}$at each point of the moduli
space of M.

The mirror symmetry suggests that we can identify the second cohomology group
of the mirror CY\ M' with $H^{2,1}$ of the original CY manifold. Since the
second cohomology group of a CY manifold has a natural $\mathbb{Z}$ structure,
then $H^{2,1}$ of the original CY manifold should also carry a natural
$\mathbb{Z}$ structure. This construction suggests that the existence of the
natural $\mathbb{Z}$ structure on $H^{2,1}$ is equivalent to the tt* equations.

In \textbf{Section 7 }we obtain an algebraic integrable system in the sense of
R. Donagi and E. Markman using the flat Cecotti-Hitchin-Simpson-Vafa
connection. From that we obtain a map from the moduli space of CY manifold M
to the moduli space of principally polarized abelian varieties and the CHSV
connection is the pull back of the connection defined by R. Donagi and E.
Markman on the moduli space of principally polarized abelian varieties. See
\cite{DM}. We also construct a Hyper-K\"{a}hler structure on the tangent
bundle of the moduli space $\mathcal{M}\left(  \text{M}\right)  $ of polarized
CY threefolds. D. Freed constructed Hyper-K\"{a}hler structure on the tangent
bundle of the relative dualizing sheaf of the moduli space $\mathcal{M}\left(
\text{M}\right)  $ of polarized CY threefolds in \cite{f}.

\textbf{Acknowledgements } I want to express my special thanks to S.
Shatashvili, who pointed out the important paper of Witten, where the anomaly
equations of BCOV were interpreted as a projective connection. The author want
to thank Y. Eliashberg, J. Li and K. Liu for stimulating conversations and
their patience with me. Special thanks to N. Nekrasov whose remarks and
suggestions helped me enormously.

\section{Deformation Theory for CY manifolds.}

\subsection{Review of \cite{to89}}

\begin{definition}
\label{com str}Let M be an even dimensional C$^{\infty}$ manifold. We will say
that M has an almost complex structure if there exists a section $I\in
C^{\infty}($M$,Hom(T^{\ast},T^{\ast}))$ such that $I^{2}=-id.$ $T$ is the
tangent bundle \ and $T^{\ast}$ is the cotangent bundle on M.
\end{definition}

This definition is equivalent to the following one:

\begin{definition}
Let M be an even dimensional C$^{\infty}$ manifold. Suppose that there exists
a global splitting of the complexified cotangent bundle $T^{\ast}%
\otimes\mathbf{C}=\Omega^{1,0}\oplus\Omega^{0,1},$ where $\Omega
^{0,1}=\overline{\Omega^{1,0}}.$ Then we will say that M has an almost complex structure.
\end{definition}

We are going to define the almost integrable complex structure.

\begin{definition}
We will say that an almost complex structure is an integrable one if for each
point $x\in$M there exists an open set $U\subset$M such that we can find local
coordinates $z^{1},..,z^{n}$ such that $dz^{1},..,dz^{n}$ \ are linearly
independent in each point $m\in U$ and they generate $\Omega^{1,0}\left\vert
_{U}\right.  .$
\end{definition}

It is easy to see that any complex manifold has an almost integrable complex structure.

\begin{definition}
\label{BELT}Let M be a complex manifold. $\phi\in\Gamma($M$,Hom(\Omega
^{1,0},\Omega^{0,1}))$ is called a Beltrami differential.
\end{definition}

Since $\Gamma($M$,Hom(\Omega^{1,0},\Omega^{0,1}))\backsimeq\Gamma($%
M$,T^{1,0}\otimes\Omega^{0,1}),$ we deduce that locally $\phi$ can be written
as follows:%
\[
\phi\left\vert _{U}\right.  =\sum\phi_{\overline{\alpha}}^{\beta}\overline
{dz}^{\alpha}\otimes\frac{\partial}{\partial z^{\beta}}.
\]
From now on we will denote by%
\[
A_{\phi}=\left(
\begin{array}
[c]{cc}%
id & \phi(\tau)\\
\overline{\phi(\tau)} & id
\end{array}
\right)  :T^{\ast}\otimes\mathbb{C}\rightarrow T^{\ast}\otimes\mathbb{C}.
\]
We will consider only those Beltrami differentials $\phi$ such that
$\det(A_{\phi})\neq0.$

\begin{definition}
\label{ics}It is easy to see that the Beltrami differential $\phi$\ defines a
new almost complex structure operator $I_{\phi}=A_{\phi}^{-1}\circ I\circ
A_{\phi}.$
\end{definition}

With respect to this new almost complex structure the space $\Omega_{\phi
}^{1,0}$ is defined as follows; if \ $dz^{1},..,dz^{n}$ generate $\Omega
^{1,0}\left\vert _{U}\right.  $, then%
\[
dz^{1}+\phi(dz^{1}),..,dz^{n}+\phi(dz^{n})
\]
generate\textit{\ }$\Omega_{\phi}^{1,0}|_{U}$\textit{\ \ }and, moreover we
have: $\overline{\Omega_{\phi}^{1,0}}\cap\Omega_{\phi}^{1,0}=0.$ The Beltrami
differential \ $\phi$\ defines an integrable complex structure on M if and
only if the following equation holds:%
\[
\overline{\partial}\phi+\frac{1}{2}\left[  \phi,\phi\right]  =0.
\]
where%
\[
\left[  \phi,\phi\right]  \left\vert _{U}\right.  :=\sum_{\nu=1}^{n}%
\sum_{1\leqq\alpha,\beta\leqq n}\left(  \sum_{\mu=1}^{n}\left(  \phi
_{\overline{\alpha}}^{\mu}\left(  \partial_{\mu}\phi_{\overline{\beta}}^{\nu
}\right)  -\phi_{\overline{\beta}}^{\mu}\left(  \partial_{\mu}\phi
_{\overline{\alpha}}^{\nu}\right)  \right)  \right)  \overline{dz}^{\alpha
}\wedge\overline{dz}^{\beta}\otimes\frac{\partial}{dz^{\nu}}.
\]
(See \cite{NN}.)

The main results in \cite{to89} are the two theorems stated bellow:

\begin{theorem}
\label{tod1}Let M be a CY manifold and let $\left\{  \phi_{i}\right\}  $ be
harmonic (with respect to the CY metric $g$) representative of the basis in
$\mathbb{H}^{1}($M$,T^{1,0}),$ then the equation: $\overline{\partial}%
\phi+\frac{1}{2}\left[  \phi,\phi\right]  =0$ has a solution in the form:%
\[
\phi(\tau_{1},..,\tau_{N})=\sum_{i=1}^{N}\phi_{i}\tau^{i}+\sum_{|I_{N}|\geqq
2}\phi_{I_{N}}\tau^{I_{N}}%
\]
\textit{where } $I_{N}=(i_{1},..,i_{N})$\ \ \textit{is a multi-index}, \
\[
\phi_{I_{N}}\in C^{\infty}(M,\Omega^{0,1}\otimes T^{1,0}),
\]
$\tau^{I_{N}}=(\tau^{i})^{i_{1}}...(\tau^{N})^{i_{N}}$ \textit{and there
exists} $\varepsilon>0$ \textit{such that }
\[
\phi(\tau)\in C^{\infty}(M,\Omega^{0,1}\otimes T^{1,0})
\]
for $|\tau^{i}|<\varepsilon$ \ for $i=1,..,N.$ $See$ $\cite{to89}.$
\end{theorem}

\begin{theorem}
\label{tod2}Let $\Omega_{0}$ \ be a holomorphic n-form on the n dimensional CY
manifold M. Let \ $\left\{  U_{i}\right\}  $be a covering of M and let
$\left\{  z_{1}^{i},..,z_{n}^{i}\right\}  $ be local coordinates in $U_{i}$
such that $\Omega_{0}\left\vert _{U_{i}}\right.  =dz_{1}^{i}\wedge...\wedge
dz_{n}^{i}.$ Then for each $\tau=(\tau^{1},..,\tau^{N})$ such that $|\tau
_{i}|<\varepsilon$ the forms on M defined as:
\[
\Omega_{t}\left\vert _{U_{i}}\right.  :=(dz_{1}^{i}+\phi(\tau)(dz_{1}%
^{i}))\wedge..\wedge(dz_{n}^{i}+\phi(\tau)(dz_{n}^{i}))
\]
are globally defined complex n forms $\Omega_{\tau}$ on M and, moreover,
$\ \Omega_{\tau}$\ are closed holomorphic n forms with respect to the complex
structure on M defined by $\phi(\tau).$
\end{theorem}

\begin{corollary}
\label{tod3}We have the following Taylor expansion for%
\begin{equation}
\Omega_{\tau}\left\vert _{U}\right.  =\Omega_{0}+\sum_{k=1}^{n}(-1)^{\frac
{k(k-1)}{2}}\left(  \wedge^{k}\phi\right)  \lrcorner\Omega_{0}. \label{T1}%
\end{equation}
(See \cite{to89}.)
\end{corollary}

From here we deduce the following Taylor expansion for the cohomology class
\ $\left[  \Omega_{\tau}\right]  $\ $\in H^{n}($M,$\mathbf{C}):$

\begin{corollary}
\label{tod4}
\begin{equation}
\left[  \Omega_{\tau}\right]  =[\Omega_{0}]-\sum_{i=1}^{N}[\left(  \phi
_{i}\lrcorner\Omega_{0}\right)  ]\tau^{i}+\frac{1}{2}\sum_{i,j=1}^{N}[(\left(
\phi_{i}\wedge\phi_{j}\right)  \lrcorner\Omega_{0})]\tau^{i}\tau^{j}%
+O(\tau^{3}) \label{T2}%
\end{equation}
(See \cite{to89}.)
\end{corollary}

We are going to define the Kuranishi family for CY manifolds of any dimension.

\begin{definition}
\label{coor}Let $\mathcal{K\subset}\mathbb{C}^{N}$ be the polydisk defined by
$|\tau^{i}|<\varepsilon$ for every $i=1,..,N$, where $\varepsilon$ is chosen
such that for every $\tau\in\mathcal{K}$ , $\phi(\tau)\in C^{\infty}%
($M,$\Omega^{0,1}\otimes T^{1,0})$, where $\phi(\tau)$ is defined as in
\ Definition \ref{BELT}. On the trivial $C^{\infty}$ family M$\times
\mathcal{K}$ we will define for each $\tau\in\mathcal{K}$ an integrable
complex structure I$_{\phi(\tau)}$ on the fibre over $\tau$ of the family
M$\times\mathcal{K}$ $,$where I$_{\phi(\tau)}$ was defined in Definition
\ref{ics}. Thus we will obtain a complex analytic family $\pi
:\mathcal{X\rightarrow K}$ of CY manifolds. We will call this family the
Kuranishi family. Thus we introduce also a coordinate system in $\mathcal{K}$.
We call this coordinate system a flat coordinate system.
\end{definition}

\subsection{Construction of the Teichm\"{u}ller Space of CY Manifolds}

\begin{definition}
We will define the Teichm\"{u}ller space $\mathcal{T}\left(  \text{M}\right)
$ of M as follows:
\[
\mathcal{T}\left(  \text{M}\right)  :=\{\text{all integrable complex
structures on M}\}\left/  \mathbf{Diff}_{0}\left(  \text{M}\right)  \right.
,
\]
where $\mathbf{Diff}_{0}\left(  \text{M}\right)  $ is the group of
diffeomorphisms of M isotopic to identity.
\end{definition}

$\mathbf{Diff}_{0}\left(  \text{M}\right)  $ acts on the complex structures as
follows: let $\psi\in\mathbf{Diff}_{0}\left(  \text{M}\right)  $ and let
\[
I\in C^{\infty}(Hom(T^{\ast}\left(  \text{M}\right)  ,T^{\ast}\left(
\text{M}\right)  ),
\]
such that $I^{2}=-id,$ then clearly $\psi^{\ast}(I)$ is such that $(\psi
^{\ast}(I))^{2}=-id.$ Moreover, if $I$ is an integrable complex structure,
then\textit{\ }$\psi^{\ast}(I)$ is integrable too.

We will call a pair $($M, $\{\gamma_{1},...,\gamma_{b_{n}}\})$ a marked CY
manifold, if M is a Calabi-Yau manifold and $\{\gamma_{1},..,\gamma_{b_{n}}\}$
is a basis in $H_{n}(M,\mathbf{Z})/Tor.$ Over the Kuranishi space we have a
universal family of marked Calabi-Yau manifolds \ $\mathcal{X\rightarrow K}$
defined up to an action of a group that acts trivially on the middle homology
and preserves the polarizations class. And, moreover, as a C$^{\infty}$
manifold $\mathcal{X}$\ \ is diffeomorphic to \ $\mathcal{K}\times$M.

\begin{theorem}
\label{teich}The Teichm\"{u}ller space$\mathcal{\ T}$ $\left(  \text{M}%
\right)  $ of a Calabi Yau manifold M exists \ as a complex manifold of
dimension $h^{2,1}.$
\end{theorem}

\textbf{Proof: }For the proof of Theorem \ref{teich} see \cite{LTYZ}.
$\blacksquare$

\subsection{Construction of the Moduli Space}

\begin{definition}
We will define the mapping class group $\Gamma^{\prime}\left(  \text{M}%
\right)  $ \ as follows:
\[
\Gamma^{\prime}\left(  \text{M}\right)  :=\mathbf{Diff}^{+}\left(
\text{M}\right)  \left/  \mathbf{Diff}_{0}\left(  \text{M}\right)  \right.  ,
\]
\textit{where }$\mathbf{Diff}^{+}\left(  \text{M}\right)  $ \textit{is the
group of diffeomorphisms preserving the orientation of M and }$\mathbf{Diff}%
_{0}\left(  \text{M}\right)  $ \textit{is the group of diffeomorphisms
isotopic to identity.}
\end{definition}

D. Sullivan proved that the mapping class group of any C$^{\infty}$ manifold
of dimension greater or equal to 5 is an arithmetic group. (See \cite{sul}.)
It is easy to prove that the mapping class group $\Gamma^{\prime}\left(
\text{M}\right)  $ acts discretely on the Teichm\"{u}ller space $\mathcal{T}%
\left(  \text{M}\right)  $ of \ the CY manifold M.

We will consider from now on polarized CY manifolds, i.e. a pair
$(M,\omega(1,1)),$ where%
\[
\lbrack\omega(1,1)]\in H^{2}\left(  \text{M},\mathbb{Z}\right)  \cap
H^{1,1}\left(  \text{M},\mathbb{R}\right)
\]
is a fixed class of cohomology and it corresponds to the imaginary part of a
CY metric. We will define $\Gamma\left(  \text{M}\right)  $ as follows:%
\[
\Gamma_{\omega(1,1)}\left(  \text{M}\right)  :=\{\phi\in\Gamma^{\prime}\left(
\text{M}\right)  \left\vert \phi([\omega(1,1)])=[\omega(1,1)]\right.  \}.
\]

From now on we will work with this family.

\begin{theorem}
\label{mod}There exists a subgroup $\Gamma\left(  \text{M}\right)  $ in
$\Gamma_{\omega(1,1)}$ of finite index such that $\Gamma$ acts without fixed
points on the Teichm\"{u}ller space $\mathcal{T}\left(  \text{M}\right)  $.
The moduli space $\mathfrak{M}\left(  \text{M}\right)  =\mathcal{T}\left(
\text{M}\right)  /\Gamma\left(  \text{M}\right)  $ is a smooth
quasi-projective variety. There exists a family of polarized CY manifolds
$\mathcal{Y}\left(  \text{M}\right)  \rightarrow\mathcal{M}\left(
\text{M}\right)  =\mathcal{T}\left(  \text{M}\right)  \left/  \Gamma\left(
\text{M}\right)  \right.  $. The relative dualizing sheaf $\omega
_{\mathcal{Y}\left/  \mathcal{M}\left(  \text{M}\right)  \right.  }$ is a
trivial line bundle.
\end{theorem}

\textbf{Proof: }Viehweg proved in \cite{vw} that the moduli space
$\mathcal{M}\left(  \text{M}\right)  $\ is a quasi projective variety. In
\cite{LTYZ} it was proved that we can find a subgroup $\Gamma\left(
\text{M}\right)  $ in $\Gamma_{\omega(1,1)}\left(  \text{M}\right)  $ such
that the space $\mathcal{T}\left(  \text{M}\right)  /\Gamma\left(
\text{M}\right)  $ is a smooth complex manifold. We also proved that over
\ $\mathcal{T}\left(  \text{M}\right)  \left/  \Gamma\left(  \text{M}\right)
\right.  =\mathcal{M}\left(  \text{M}\right)  $ there exists a family of CY
manifolds $\mathcal{Y}\left(  \text{M}\right)  \rightarrow\mathcal{M}\left(
\text{M}\right)  .$ In \cite{BT} we proved the following Theorem:

\begin{theorem}
\label{BGS1}Let $\mathcal{M}\left(  \text{M}\right)  =\mathcal{T}\left(
\text{M}\right)  \left/  \Gamma\left(  \text{M}\right)  \right.  $. Then there
exists a global non vanishing holomorphic section $\eta_{\tau}$\ of the line
bundle $\omega_{\mathcal{Y}\left(  \text{M}\right)  \left/  \mathcal{M}\left(
\text{M}\right)  \right.  }$ whose $L^{2}$ norm $\left\Vert \eta_{\tau
}\right\Vert _{\mathbf{L}^{2}}^{2}$ is equal to $\left(  \det_{(0,1)}\right)
,$ where $\det_{(0,1)}$ is the regularized determinant of the Laplacian of a
CY acting on $\Omega_{\text{M}}^{0,1}$ of the CY metric with imaginary class
equal to the polarization class and $\omega_{\mathcal{Y}\left(  \text{M}%
\right)  \left/  \mathcal{M}\left(  \text{M}\right)  \right.  }$ is a trivial
holomorphic line bundle.
\end{theorem}

Theorem \ref{mod} follows from Theorem \ref{BGS1}. $\blacksquare$

\begin{corollary}
\label{bt1}$\eta_{\tau}$ defines a flat structure on $\omega_{\mathcal{Y}%
\left(  \text{M}\right)  \left/  \mathcal{M}\left(  \text{M}\right)  \right.
}.$
\end{corollary}

\subsection{Weil-Petersson Geometry}

In our paper \cite{to89} we define a metric on the Kuranishi space
$\mathcal{K}$\ and called this metric, the Weil-Petersson metric. We will
review the basic properties of the Weil-Petersson metric which were
established in \cite{to89}. In \cite{to89} we proved the following theorem:

\begin{theorem}
\label{hodge}Let M be a CY manifold of dimension n and let $\Omega_{\text{M}}$
be a non zero holomorphic n form on M such that%
\[
\left(  -1\right)  ^{\frac{n(n-1)}{2}}\left(  \sqrt{-1}\right)  ^{n}%
{\displaystyle\int\limits_{\text{M}}}
\Omega_{\text{M}}\wedge\overline{\Omega_{\text{M}}}=1.
\]
\textit{Let g be a Ricci flat (CY) metric on M. Then the map:}%
\[
\psi\in L^{2}(\text{M},\Omega_{\text{M}}^{0,k}\wedge^{m}T_{\text{M}}%
^{1,0})\rightarrow\psi\lrcorner\Omega_{\text{M}}\in L^{2}(\text{M}%
,\Omega_{\text{M}}^{n-m,k})
\]
\textit{gives an isomorphism between Hilbert spaces and this map preserves the
Hodge decomposition.}\cite{to89}.
\end{theorem}

\begin{corollary}
We can identify the tangent space $T_{\tau}=H^{1}($M$_{\tau},T_{\tau}^{1,0})$
\textit{at each point} $\ \tau\in\mathcal{T}\left(  \text{M}\right)
$\ \textit{with}$\ H^{1}($M$_{\tau},\Omega_{\tau}^{n-1})$, \textit{by using
the map }$\psi\rightarrow\psi\lrcorner\Omega_{\text{M}}.$
\end{corollary}

\begin{notation}
We will denote by
\begin{equation}
\left\langle \omega_{1},\omega_{2}\right\rangle :=%
{\displaystyle\int\limits_{\text{M}_{\tau}}}
\omega_{1}\wedge\overline{\omega_{2}}. \label{Sym}%
\end{equation}

\end{notation}

\begin{definition}
Let $\psi_{1},$ $\psi_{2}\in T_{\tau}=\mathbf{H}^{1}\left(  \text{M}_{\tau
},T_{\text{M}_{\tau}}^{1,0}\right)  $ $($the space of harmonic forms with
respect to the CY metric g.). We will define the Weil-Petersson metric as
follows:%
\[
\left\langle \psi_{1},\psi_{2}\right\rangle _{WP}:=\sqrt{-1}%
{\displaystyle\int\limits_{\text{M}_{\tau}}}
\left(  \psi_{1}\lrcorner\Omega_{\tau}\right)  \wedge\left(  \overline
{\psi_{2}\lrcorner\Omega_{\tau}}\right)  =\sqrt{-1}\left\langle \psi
_{1}\lrcorner\Omega_{\tau},\overline{\psi_{2}\lrcorner\Omega_{\tau}%
}\right\rangle
\]
and $\left\Vert \Omega_{\tau}\right\Vert ^{2}=1.$ Thus $\left\langle \psi
,\psi\right\rangle _{WP}>0$.
\end{definition}

The Weil-Petersson metric is a K\"{a}hler metric on the Teichm\"{u}ller space
\ $\mathcal{T}\left(  \text{M}\right)  $. It defines a natural connection,
namely the Levi-Civita connection$\;\not .  \;$ We will denote the covariant
derivatives in direction $\frac{\partial}{\partial\tau^{i}}$ at the tangent
space of a point $\tau\in\mathcal{T}\left(  \text{M}\right)  $ defined by
$\phi_{i}$ by $\nabla_{i}.$ In \cite{to89} we proved the following theorem:

\begin{theorem}
\label{WP} In the flat coordinate system introduced in Definition \ref{coor}
the following formulas hold for the curvature operator:%
\[
R_{i\overline{j},k\overline{l}}=\delta_{i\overline{j}}\delta_{k\overline{l}%
}+\delta_{i\overline{l}}\delta_{k\overline{j}}-\sqrt{-1}%
{\displaystyle\int\limits_{\text{M}}}
((\phi_{i}\wedge\phi_{k})\lrcorner\Omega_{\text{M}}))\wedge(\overline
{(\phi_{j}\wedge\phi_{l})\lrcorner\Omega_{\text{M}}})
\]%
\begin{equation}
=\delta_{i\overline{j}}\delta_{k\overline{l}}+\delta_{i\overline{l}}%
\delta_{k\overline{j}}-\sqrt{-1}\left\langle (\phi_{i}\wedge\phi_{k}%
)\lrcorner\Omega_{\text{M}})),(\phi_{j}\wedge\phi_{l})\lrcorner\Omega
_{\text{M}}\right\rangle . \label{wp0}%
\end{equation}

\end{theorem}

\section{Flat $\mathbb{S}$p(2h$^{2,1},\mathbb{R)}$ Structure on the Moduli
Space of CY Threefolds}

\subsection{A Flat Structure on the Line Bundle $\omega_{\mathcal{Y}\left(
\text{M}\right)  \left/  \mathcal{M}\left(  \text{M}\right)  \right.  }$}

The flat structure on the line bundle $\omega_{\mathcal{Y}\text{ }\left(
\text{M}\right)  \text{ /}\mathcal{M}\left(  \text{M}\right)  }$ is defined by
Corollary \ref{bt1}.

\subsection{Gauss-Manin Connection}

\begin{definition}
\label{manin}On the Teichm\"{u}ller space $\mathcal{T}\left(  \text{M}\right)
$ we have a trivial bundle namely%
\[
\mathcal{H}^{n}=H^{n}\left(  \text{M},\mathbb{C}\right)  \times\mathcal{T}%
\left(  \text{M}\right)  \rightarrow\mathcal{T}\left(  \text{M}\right)  .
\]
\textit{Theorem} \ref{mod}\textit{\ implies that we constructed the moduli
space }$\mathcal{M}\left(  \text{M}\right)  $ as $\mathcal{T}\left(
\text{M}\right)  $/$\Gamma\left(  \text{M}\right)  $. Thus \textit{we obtain a
natural representation of the group }$\Gamma\left(  \text{M}\right)  $ into
$H^{n}\left(  \text{M},\mathbb{C}\right)  $ and a \textit{flat connection on
the flat bundle}%
\[
\mathcal{H}^{n}\left/  \Gamma\left(  \text{M}\right)  \right.  \rightarrow
\mathcal{T}\left(  \text{M}\right)  \left/  \Gamma\left(  \text{M}\right)
\right.  =\mathcal{M}\left(  \text{M}\right)  .
\]
\textit{This connection is called the Gauss-Manin connection. The covariant
derivative in direction }$\phi_{i}$\textit{\ of the tangent space \ T}%
$_{\tau,\mathcal{M}\left(  \text{M}\right)  }$\textit{with respect to the
Gauss-Manin connection will be denoted by} \ $\mathcal{D}_{i}$.
\end{definition}

The Gauss-Manin connection $\mathcal{D}$ is defined in a much more general
situation and it is defined on the moduli space of CY manifolds of dimension
$n\geq3$. We will state explicit formulas for the covariant differentiation
$\mathcal{D}_{i}$ defined by the Gauss-Manin connection. We will fix a
holomorphic three form $\Omega_{0}$ such that%
\[
-\sqrt{-1}\left\langle \Omega_{0},\Omega_{0}\right\rangle =-\sqrt{-1}%
{\displaystyle\int\limits_{\text{M}}}
\Omega_{0}\wedge\overline{\Omega_{0}}=\left\Vert \Omega_{0}\right\Vert
^{2}=1.
\]
Using the form $\Omega_{0}$ and theorem \ref{hodge}, we can identify the
cohomology groups $H^{1}\left(  \text{M},T_{\text{M}}^{1,0}\right)  $ and
$H^{1}\left(  \text{M},\Omega_{\text{M}}^{2}\right)  $ on M:

\begin{proposition}
\label{ident}The map%
\begin{equation}
\mathfrak{i}:\psi\rightarrow\psi\lrcorner\Omega_{\text{M}} \label{Id}%
\end{equation}
is an isomorphism between the groups $H^{1}\left(  \text{M},T_{\text{M}}%
^{1,0}\right)  $ and $H^{1}\left(  \text{M},\Omega_{\text{M}}^{2}\right)  $.
\end{proposition}

\textbf{Proof: }Our proposition follows directly from Theorem \ref{hodge}.
$\blacksquare$

\begin{remark}
From now on in the map $\left(  \ref{Id}\right)  $ we will use for
$\Omega_{\text{M}}$ the restriction of the holomorphic form $\eta_{\tau}$ on M
defined by Theorem \ref{BGS1}.
\end{remark}

\begin{remark}
\label{ident1}Suppose that M is a three dimensional CY manifold. Then the
Poincare map identifies $H^{1}\left(  \text{M},\Omega_{\text{M}}^{2}\right)  $
with $H^{2}\left(  \text{M},\Omega_{\text{M}}^{1}\right)  .$ This
identification will be denoted by $\Pi$, i.e. $\Pi:H^{2}\left(  \text{M}%
,\Omega_{\text{M}}^{1}\right)  \rightarrow H^{1}\left(  \text{M}%
,\Omega_{\text{M}}^{2}\right)  $ and it is defined by identifying some basis
$\left\{  \Omega_{i}\right\}  $ of $\ H^{2}\left(  \text{M},\Omega_{\text{M}%
}^{1}\right)  $ with the basis $\left\{  \overline{\Omega_{i}}\right\}  $ of
$H^{1}\left(  \text{M},\Omega_{\text{M}}^{2}\right)  .$ So:%
\begin{equation}
\Pi(\Omega_{i}):=\overline{\Omega_{i}}. \label{Poin}%
\end{equation}

\end{remark}

\begin{notation}
\label{ident2} Using Proposition \ref{ident} and Remark \ref{ident1} one can
identify the spaces $H^{1}\left(  \text{M},T_{\text{M}}^{1,0}\right)  $ and
$H^{2}\left(  \text{M},\Omega_{\text{M}}^{1}\right)  $ for three dimensional
CY by using the map $F$, where%
\begin{equation}
F(\phi):=\Pi(\iota(\phi))\text{ } \label{Poin1}%
\end{equation}
for $\phi\in H^{1}($M$,T_{\text{M}}^{1,0}).$
\end{notation}

\begin{lemma}
\label{manin1}Let $\mathfrak{i}^{-1}:H^{1}\left(  \text{M},\Omega_{\text{M}%
}^{2}\right)  \overset{\llcorner\Omega_{\text{M}}^{\ast}}{\backsimeq}%
H^{1}\left(  \text{M},T_{\text{M}}^{1,0}\right)  $ be the inverse
identification defined by $\left(  \ref{Id}\right)  $. Then%
\begin{equation}
\mathcal{D}_{i}(\mathfrak{i(}\phi))=\mathfrak{\iota}\left(  \phi_{i}\right)
\lrcorner\phi\in H^{2}\left(  \text{M}_{0},\Omega_{\text{M}_{0}}^{1}\right)  .
\label{T3}%
\end{equation}

\end{lemma}

\textbf{Proof: }The proof of the lemma follows directly from formula$\left(
\ref{T2}\right)  $. Indeed since $\{\phi_{i}\}$ is a basis of $H^{1}%
($M$,T_{\text{M}}^{1,0})$, then $\{\phi_{i}\lrcorner\Omega_{\text{M}}\}$ will
be a basis of $H^{1}\left(  \text{M},\Omega_{\text{M}}^{2}\right)  $. In the
flat coordinates $(\tau^{1},...,\tau^{N})$ introduced in Definition \ref{coor}
and according to $\left(  \ref{T2}\right)  $ we have
\[
\left[  \Omega_{\tau}\right]  =[\Omega_{\text{M}}]-\sum_{i=1}^{N}\left[
\left(  \phi_{i}\lrcorner\Omega_{\text{M}}\right)  \right]  \tau^{i}%
+\sum_{i,j=1}^{N}\left[  \left(  \phi_{i}\wedge\phi_{j}\right)  \lrcorner
\Omega_{\text{M}}\right]  \tau^{i}\tau^{j}+O(\tau^{3}).\
\]
From Definition \ref{manin} and the expression of $[\Omega_{\tau}]$ given by
the above formula, we deduce formula $\left(  \ref{T3}\right)  $. Lemma
\ref{manin}\ is proved. $\blacksquare$

\subsection{Higgs Fields and the Tangent Space of $\mathcal{M}\left(
\text{M}\right)  $}

Let $\omega_{\mathcal{Y}\left(  \text{M}\right)  \left/  \mathcal{M}\left(
\text{M}\right)  \right.  }$ be the relative dualizing sheaf on $\mathcal{M}%
\left(  \text{M}\right)  .$

\begin{definition}
\label{Higgs0}We define the Higgs field of a holomorphic bundle $\mathcal{E}$
over a complex manifold M as a globally defined holomorphic map $\Phi
:\mathcal{E\rightarrow E\otimes}\Omega_{\text{M}}^{1}$ such that
\begin{equation}
\Phi\circ\Phi=0. \label{Poin2}%
\end{equation}

\end{definition}

\begin{lemma}
\label{tang}The tangent bundle $T_{\mathcal{T}\left(  \text{M}\right)  }$ of
the Teichm\"{u}ller space $\mathcal{T}\left(  \text{M}\right)  $ of any CY
manifold M is canonically isomorphic to the bundles%
\[
T_{\mathcal{T}\left(  \text{M}\right)  }\thickapprox Hom\left(  \pi_{\ast
}\omega_{\mathcal{Y}\left(  \text{M}\right)  \left/  \mathcal{M}\left(
\text{M}\right)  \right.  },R^{1}\pi_{\ast}\Omega_{\mathcal{Y}\left(
\text{M}\right)  \left/  \mathcal{M}\left(  \text{M}\right)  \right.  }%
^{2}\right)  \thickapprox
\]%
\[
\left(  \pi_{\ast}\omega_{\mathcal{Y}\left(  \text{M}\right)  \left/
\mathcal{M}\left(  \text{M}\right)  \right.  }\right)  ^{\ast}\otimes R^{1}%
\pi_{\ast}\Omega_{\mathcal{Y}\left(  \text{M}\right)  \left/  \mathcal{M}%
\left(  \text{M}\right)  \right.  }^{2}.
\]

\end{lemma}

\textbf{Proof: }The proof of this lemma is standard and follows from local
Torelli Theorem. $\blacksquare$

\subsection{Construction of a Higgs Field on $R^{1}\pi_{\ast}\Omega
_{\mathcal{Y}\left(  \text{M}\right)  \left/  \mathcal{M}\left(
\text{M}\right)  \right.  }^{2}$}

From now on we will consider only three dimensional CY manifolds.

\begin{definition}
\label{gm}We define a Higgs field $\tilde{\nabla}$ on $R^{1}\pi_{\ast}%
\Omega_{\mathcal{Y}\left(  \text{M}\right)  \left/  \mathcal{M}\left(
\text{M}\right)  \right.  }^{2}$ by using the Gauss-Manin connection \ and
Poincare duality $\Pi$ in the following manner:%
\begin{equation}
\tilde{\nabla}_{\phi_{i}}\Omega=\tilde{\nabla}_{i}\Omega:=\Pi(\mathcal{D}%
_{i}(\Omega))=\Pi(\Omega\lrcorner\phi_{i}), \label{Poin3}%
\end{equation}
where $\Omega\in H^{1}\left(  \text{M},\Omega^{2}\right)  $ $\ $and$\{\phi
_{i}\}\in H^{1}\left(  \text{M},T_{\text{M}}^{1,0}\right)  $ is an orthonormal
basis with respect to the Weil-Petersson metric on $T_{0}\thickapprox
\mathbb{H}^{1}\left(  \text{M},T_{\text{M}}^{1,0}\right)  $.
\end{definition}

\begin{lemma}
\label{Higgs-1}$\tilde{\nabla}$ as defined in Definition \ref{gm} is a Higgs field.
\end{lemma}

\textbf{Proof: }We need to check that $\tilde{\nabla}$ satisfies the
conditions in the definition \ref{Higgs0} of a Higgs field. The definition
$\left(  \ref{Poin3}\right)  $ of $\tilde{\nabla}$ implies that%

\begin{equation}
\tilde{\nabla}\in Hom\left(  R^{1}\pi_{\ast}\Omega_{\mathcal{Y}\left(
\text{M}\right)  \left/  \mathcal{M}\left(  \text{M}\right)  \right.  }%
^{2}\otimes T(\mathcal{M}\left(  \text{M}\right)  \right)  ,R^{1}\pi_{\ast
}\Omega_{\mathcal{Y}\left(  \text{M}\right)  \left/  \mathcal{M}\left(
\text{M}\right)  \right.  }^{2}). \label{Poin4}%
\end{equation}
On the other hand standard facts from commutative algebra imply%
\[
Hom\left(  R^{1}\pi_{\ast}\Omega_{\mathcal{Y}\left(  \text{M}\right)  \left/
\mathcal{M}\left(  \text{M}\right)  \right.  }^{2}\otimes T(\mathcal{M}\left(
\text{M}\right)  ,R^{1}\pi_{\ast}\Omega_{\mathcal{Y}\left(  \text{M}\right)
\left/  \mathcal{M}\left(  \text{M}\right)  \right.  }^{2}\right)  \backsimeq
\]%
\[
\left(  R^{1}\pi_{\ast}\Omega_{\mathcal{Y}\left(  \text{M}\right)  \left/
\mathcal{M}\left(  \text{M}\right)  \right.  }^{2}\right)  ^{\ast}%
\otimes\left(  T(\mathcal{M}\left(  \text{M}\right)  )\right)  ^{\ast}\otimes
R^{1}\pi_{\ast}\Omega_{\mathcal{Y}\left(  \text{M}\right)  \left/
\mathcal{M}\left(  \text{M}\right)  \right.  }^{2}\backsimeq
\]%
\[
\backsimeq\left(  R^{1}\pi_{\ast}\Omega_{\mathcal{Y}\left(  \text{M}\right)
\left/  \mathcal{M}\left(  \text{M}\right)  \right.  }^{2}\right)  ^{\ast
}\otimes\left(  \Omega_{\mathcal{M}\left(  \text{M}\right)  }^{1}\otimes
R^{1}\pi_{\ast}\Omega_{\mathcal{Y}\left(  \text{M}\right)  \left/
\mathcal{M}\left(  \text{M}\right)  \right.  }^{2}\right)  \backsimeq
\]

\begin{equation}
\backsimeq Hom\left(  R^{1}\pi_{\ast}\Omega_{\mathcal{Y}\left(  \text{M}%
\right)  \left/  \mathcal{M}\left(  \text{M}\right)  \right.  }^{2}%
,\Omega_{\mathcal{M}\left(  \text{M}\right)  }^{1}\otimes R^{1}\pi_{\ast
}\Omega_{\mathcal{Y}\left(  \text{M}\right)  \left/  \mathcal{M}\left(
\text{M}\right)  \right.  }^{2}\right)  . \label{Poin5}%
\end{equation}
So we obtain
\[
Hom\left(  R^{1}\pi_{\ast}\Omega_{\mathcal{Y}\left(  \text{M}\right)  \left/
\mathcal{M}\left(  \text{M}\right)  \right.  }^{2}\otimes T(\mathcal{M}\left(
\text{M}\right)  ),R^{1}\pi_{\ast}\Omega_{\mathcal{Y}\left(  \text{M}\right)
\left/  \mathcal{M}\left(  \text{M}\right)  \right.  }^{2}\right)  \backsimeq
\]%
\begin{equation}
\backsimeq Hom\left(  R^{1}\pi_{\ast}\Omega_{\mathcal{Y}\left(  \text{M}%
\right)  \left/  \mathcal{M}\left(  \text{M}\right)  \right.  }^{2}%
,\Omega_{\mathcal{M}\left(  \text{M}\right)  }^{1}\otimes R^{1}\pi_{\ast
}\Omega_{\mathcal{Y}\left(  \text{M}\right)  \left/  \mathcal{M}\left(
\text{M}\right)  \right.  }^{2}\right)  . \label{Poin6}%
\end{equation}
So the condition that
\[
\tilde{\nabla}\in Hom\left(  R^{1}\pi_{\ast}\Omega_{\mathcal{Y}\left(
\text{M}\right)  \left/  \mathcal{M}\left(  \text{M}\right)  \right.  }%
^{2},R^{1}\pi_{\ast}\Omega_{\mathcal{Y}\left(  \text{M}\right)  \left/
\mathcal{M}\left(  \text{M}\right)  \right.  }^{2}\otimes\Omega_{\mathcal{M}%
\left(  \text{M}\right)  }^{1}\right)
\]
follows directly from $\left(  \ref{Poin4}\right)  $ and $\left(
\ref{Poin6}\right)  $. Next we need to prove that $\tilde{\nabla}\circ
\tilde{\nabla}=0.$ It is a standard fact that the relation $\tilde{\nabla
}\circ\tilde{\nabla}=0$ is equivalent to the relations $[\tilde{\nabla}%
_{i},\tilde{\nabla}_{j}]=0.$ So in order to finish the proof of Lemma
\ref{Higgs-1} we need to prove the following Proposition:

\begin{proposition}
\label{Higgs1}The commutator of $\tilde{\nabla}_{i}$ and $\tilde{\nabla}_{j}$
is equal to zero, i.e.
\begin{equation}
\lbrack\tilde{\nabla}_{i},\tilde{\nabla}_{j}]=0. \label{Poin7}%
\end{equation}

\end{proposition}

\textbf{Proof: }Lemma \ref{Higgs1} follows directly from the definition of
$\ \tilde{\nabla}.$ Indeed since the Gauss-Manin connection is a flat, i.e.
$[\mathcal{D}_{i},\mathcal{D}_{j}]=0$ and the fact that the covariant
derivative of Poincare duality is zero we deduce that:
\[
\left[  \tilde{\nabla}_{i},\tilde{\nabla}_{j}\right]  =\left[  \Pi
(\mathcal{D}_{i}),\Pi\left(  \mathcal{D}_{j}\right)  \right]  =\Pi
\lbrack\mathcal{D}_{i},\mathcal{D}_{j}]=0.
\]
Proposition \ref{Higgs1} is proved. $\blacksquare$ Lemma \ref{Higgs-1} is
proved. $\blacksquare$

\subsection{Higgs Field on the Tangent Bundle of $\mathcal{M}\left(
\text{M}\right)  $}

Let $\omega_{\mathcal{Y}\left(  \text{M}\right)  \left/  \mathcal{M}\left(
\text{M}\right)  \right.  }$ be the relative dualizing sheaf on $\mathcal{M}%
\left(  \text{M}\right)  .$ We will define bellow a Higgs Field on the tangent
vector bundle \
\[
T_{\mathcal{M}\left(  \text{M}\right)  }\thickapprox\left(  \omega
_{\mathcal{Y}\left(  \text{M}\right)  \left/  \mathcal{M}\left(
\text{M}\right)  \right.  }\right)  ^{\ast}\otimes R^{1}\pi_{\ast}%
\Omega_{\mathcal{Y}\left(  \text{M}\right)  \left/  \mathcal{M}\left(
\text{M}\right)  \right.  }^{2}%
\]
of the moduli space \textit{$\mathcal{M}\left(  \text{M}\right)  $.} We denote
by $\theta=id\otimes\tilde{\nabla},$ where $\tilde{\nabla}$ is defined by
$\left(  \ref{Poin3}\right)  .$

\begin{lemma}
\label{gm1} Let $\{\phi_{i}\}$ be a basis of the tangent bundle $\mathcal{T}%
_{\mathcal{M}\left(  \text{M}\right)  }$ restricted on some open polydisk
$U\subset\mathcal{M}\left(  \text{M}\right)  .$ We can identify the fibre
$T_{\text{M}_{\tau}}^{1,0}$ of $\mathcal{T}_{\mathcal{M}\left(  \text{M}%
\right)  }$ with the harmonic forme $\mathbb{H}^{1}\left(  \text{M}_{\tau
},T_{\text{M}_{\tau}}^{1,0}\right)  $ with respect to the CY metric
corresponding to the [p;arization class $L.$ Let us define
\[
\vartheta:\mathcal{T}_{\mathcal{M}\left(  \text{M}\right)  }\rightarrow
\mathcal{T}_{\mathcal{M}\left(  \text{M}\right)  }\left(  \mathcal{M}\left(
\text{M}\right)  \right)  \otimes\Omega_{\mathcal{M}\left(  \text{M}\right)
}^{1}\left(  \mathcal{M}\left(  \text{M}\right)  \right)
\]
by
\[
\vartheta\left(  \eta_{\tau}^{-1}\otimes\phi_{j}\right)  :=\eta_{\tau}%
^{-1}\otimes\sum_{i=1}^{N}\left(  \digamma^{-1}\left(  \mathcal{D}_{i}%
(\phi_{j}\lrcorner\Omega_{\tau})\right)  \otimes(\phi_{i})^{\ast}\right)
\]
where $\mathcal{D}$ is the Gauss-Manin connection defined by $\left(
\ref{Poin3}\right)  $, $\digamma$ is defined by $\left(  \ref{Poin1}\right)  $
and $\eta_{\tau}$ is the holomorphic three form on M \ defined by Theorem
\ref{BGS1}. Then
\begin{equation}
\vartheta_{i}\left(  \eta_{\tau}^{-1}\otimes\phi_{j}\right)  =\eta_{\tau}%
^{-1}\otimes\left(  \left(  \left(  \phi_{j}\lrcorner\Omega_{\tau}\right)
\lrcorner\phi_{i}\right)  \lrcorner\left(  \Omega_{\tau}\right)  ^{\ast
}\right)  ,\label{GM}%
\end{equation}
and $\vartheta$ defines a Higgs field on the tangent bundle of $\mathcal{M}%
\left(  \text{M}\right)  .$
\end{lemma}

\textbf{Proof: }Lemma \ref{gm1} follows directly from Lemma \ref{Higgs-1}, the
definition of $\vartheta$ and the fact that
\[
T_{\mathcal{M}\left(  \text{M}\right)  }\thickapprox\left(  \omega
_{\mathcal{Y}\left(  \text{M}\right)  \left/  \mathcal{M}\left(
\text{M}\right)  \right.  }\right)  ^{\ast}\otimes R^{1}\pi_{\ast}%
\Omega_{\mathcal{Y}\left(  \text{M}\right)  \left/  \mathcal{M}\left(
\text{M}\right)  \right.  }^{2}.
\]
Lemma \ref{gm1} is proved. $\blacksquare$

We know that $h^{2,1}=\dim_{\mathbb{C}}H^{1}\left(  \text{M\thinspace,}%
\Omega_{\text{M}}^{2}\right)  =h^{1,2}=\dim_{\mathbb{C}}H^{2}\left(
\text{M\thinspace,}\Omega_{\text{M}}^{1}\right)  $ are constants. Therefore
$R^{1}\pi_{\ast}\Omega_{\mathcal{Y}\left(  \text{M}\right)  \left/
\mathcal{M}\left(  \text{M}\right)  \right.  }^{2}$ and $R^{2}\pi_{\ast}%
\Omega_{\mathcal{Y}\left(  \text{M}\right)  \left/  \mathcal{M}\left(
\text{M}\right)  \right.  }^{1}$ are holomorphic vector bundles over the
moduli space $\mathcal{M}\left(  \text{M}\right)  $ of polarized CY threefold.
We have a non degenerate pairing:%
\begin{equation}
R^{1}\pi_{\ast}\Omega_{\mathcal{Y}\left(  \text{M}\right)  \left/
\mathcal{M}\left(  \text{M}\right)  \right.  }^{2}\times R^{2}\pi_{\ast}%
\Omega_{\mathcal{Y}\left(  \text{M}\right)  \left/  \mathcal{M}\left(
\text{M}\right)  \right.  }^{1}\rightarrow R^{3}\pi_{\ast}\Omega
_{\mathcal{Y}\left(  \text{M}\right)  \left/  \mathcal{M}\left(
\text{M}\right)  \right.  }^{3}\label{hk0}%
\end{equation}
given by
\begin{equation}
\sqrt{-1}\left(  \left(  \phi_{j}(\tau)\lrcorner\Omega_{\tau}\right)  \right)
\wedge\omega_{j}(\tau)=h_{ij}(\tau)\Omega_{\tau}\wedge\overline{\Omega_{\tau}%
},\label{hk}%
\end{equation}
where $\left(  \phi_{j}(\tau)\lrcorner\Omega_{\tau}\right)  $ and $\omega
_{j}(\tau)$ are holomorphic sections of the holomorphic vector bundles
$R^{1}\pi_{\ast}\Omega_{\mathcal{Y}\left(  \text{M}\right)  \left/
\mathcal{M}\left(  \text{M}\right)  \right.  }^{2}$ and $R^{2}\pi_{\ast}%
\Omega_{\mathcal{Y}\left(  \text{M}\right)  \left/  \mathcal{M}\left(
\text{M}\right)  \right.  }^{1}.$ Clearly $h_{ij}(\tau)$ depends
holomorphically on $\tau\in\mathcal{M}\left(  \text{M}\right)  .$ By using the
non-degenerate pairing $\left(  \ref{hk0}\right)  $ we can identify $R^{1}%
\pi_{\ast}\Omega_{\mathcal{Y}\left(  \text{M}\right)  \left/  \mathcal{M}%
\left(  \text{M}\right)  \right.  }^{2}$ with $R^{2}\pi_{\ast}\Omega
_{\mathcal{Y}\left(  \text{M}\right)  \left/  \mathcal{M}\left(
\text{M}\right)  \right.  }^{1}$ as follows: to the basis of holomorphic
sections
\[
\left\{  \phi_{i}\lrcorner\Omega_{\tau_{0}}\right\}  \in R^{1}\pi_{\ast}%
\Omega_{\mathcal{Y}\left(  \text{M}\right)  \left/  \mathcal{M}\left(
\text{M}\right)  \right.  }^{2}%
\]
we will assign
\[
\left(  \phi_{i}\lrcorner\Omega_{\tau_{0}}\right)  ^{\ast}\in Hom\left(
R^{2}\pi_{\ast}\Omega_{\mathcal{Y}\left(  \text{M}\right)  \left/
\mathcal{M}\left(  \text{M}\right)  \right.  }^{1},R^{3}\pi_{\ast}%
\Omega_{\mathcal{Y}\left(  \text{M}\right)  \left/  \mathcal{M}\left(
\text{M}\right)  \right.  }^{3}\right)
\]
such that for the pairing defined by $\left(  \ref{hk}\right)  $ satisfies:%
\[
\phi_{i}\lrcorner\Omega_{\tau}\wedge\left(  \phi_{i}\lrcorner\Omega_{\tau_{0}%
}\right)  ^{\ast}=
\]%
\begin{equation}
\phi_{i}\lrcorner\Omega_{\tau}\wedge\left(  \phi_{i}\lrcorner\Omega_{\tau_{0}%
}\right)  ^{\ast}=\delta_{i,j}\Omega_{\tau}\wedge\overline{\Omega_{\tau}%
}.\label{hm2}%
\end{equation}
So $\left(  \phi_{i}\lrcorner\Omega_{\tau_{0}}\right)  $ and $\omega_{l}%
^{\ast}$ are given by the formula:%
\begin{equation}
\omega_{l}^{\ast}=%
{\displaystyle\sum\limits_{k}}
h^{lk}\left(  \phi_{k}\lrcorner\Omega_{\tau}\right)  \label{hm3}%
\end{equation}%
\begin{equation}
\left(  \phi_{l}\lrcorner\Omega_{\tau}\right)  ^{\ast}=%
{\displaystyle\sum\limits_{k}}
h^{kl}\omega_{k}.\label{hm3a}%
\end{equation}

\begin{lemma}
\label{flat} Let us a fix a point $\tau_{0}\in\mathcal{M}\left(
\text{M}\right)  .$ Suppose that
\[
\left\{  \phi_{i},\text{ }i=1,...,N\right\}  \text{ and }\left\{  \omega
_{i},\text{ }i=1,...,N\right\}
\]
are some bases of $T_{\tau,\mathcal{M}\left(  \text{M}\right)  }%
(U)=\Omega_{\tau}\otimes R^{1}\pi_{\ast}\Omega_{\mathcal{Y}\left(
\text{M}\right)  \left/  \mathcal{M}\left(  \text{M}\right)  \right.  }^{2}$
and $R^{2}\pi_{\ast}\Omega_{\mathcal{Y}\left(  \text{M}\right)  \left/
\mathcal{M}\left(  \text{M}\right)  \right.  }^{1}$ in some open polydisk of
$\tau_{0}.$ Let us define%
\begin{equation}
C_{ijl}=-\sqrt{-1}%
{\displaystyle\int\limits_{\text{M}}}
(\Omega_{\tau}\wedge\left(  (\phi_{i}\wedge\phi_{j}\wedge\phi_{l}%
)\lrcorner\Omega_{\tau}\right)  ). \label{H7}%
\end{equation}
Let $h_{ij}:=\left\langle \phi_{i}\lrcorner\Omega_{\tau},\omega_{j}%
\right\rangle $ be the pairing between the holomorphic vector bundles
$R^{1}\pi_{\ast}\Omega_{\mathcal{Y}\left(  \text{M}\right)  \left/
\mathcal{M}\left(  \text{M}\right)  \right.  }^{2}$ and $R^{2}\pi_{\ast}%
\Omega_{\mathcal{Y}\left(  \text{M}\right)  \left/  \mathcal{M}\left(
\text{M}\right)  \right.  }^{1}$ defined by $\left(  \ref{hk}\right)  .$ Then
\begin{equation}
\vartheta_{i}\phi_{j}=\sum_{k,l=1}^{N}C_{ij}^{k}\left(  \omega_{k}^{\ast
}\lrcorner\Omega_{\tau}^{\ast}\right)  , \label{H1}%
\end{equation}%
\[
\vartheta_{i}\phi_{j}=\sum_{k=1}^{N}C_{ijk}\phi_{k},
\]
where $C_{ij}^{k}$ and $C_{ijk}$ are holomorphic functions in $U$. The
relations between $C_{ij}^{k}$ and $C_{ijl}$ are given by\
\begin{equation}
C_{ij}^{k}=%
{\displaystyle\sum\limits_{m=1}^{n}}
C_{ijm}h^{mk}, \label{h1}%
\end{equation}
and,
\begin{equation}
C_{ijl}=C_{jil}=C_{ilj}=C_{jli}=C_{lij}=C_{lij}. \label{H2}%
\end{equation}

\end{lemma}

\textbf{Proof: }$\left(  \ref{GM}\right)  $ implies that
\[
\vartheta_{i}\phi_{j}=\left(  \left(  \phi_{j}\lrcorner\Omega_{\tau}\right)
\lrcorner\phi_{i}\right)  \lrcorner\left(  \Omega_{\tau}\right)  ^{\ast
}=\left(  \phi_{i}\wedge\phi_{j}\lrcorner\Omega_{\tau}\right)  \lrcorner
\Omega_{\tau},
\]
where $\left(  \phi_{i}\wedge\phi_{j}\lrcorner\Omega_{\tau}\right)  \in
R^{2}\pi_{\ast}\Omega_{\mathcal{Y}\left(  \text{M}\right)  \left/
\mathcal{M}\left(  \text{M}\right)  \right.  }^{1}\left\vert _{U}\right.  .$
Since $\left\{  \omega_{k}\right\}  $ is a basis of the holomorphic vector
bundle $R^{2}\pi_{\ast}\Omega_{\mathcal{Y}\left(  \text{M}\right)  \left/
\mathcal{M}\left(  \text{M}\right)  \right.  }^{1}\left\vert _{U}\right.  $ we
have
\begin{equation}
\left(  \phi_{i}\wedge\phi_{j}\lrcorner\Omega_{\tau}\right)  =%
{\displaystyle\sum\limits_{k}}
C_{ij}^{k}\omega_{k}.\label{H6}%
\end{equation}
Poincare duality and $\left(  \ref{H6}\right)  $ imply that we can identify
\[
\left(  \phi_{i}\wedge\phi_{j}\lrcorner\Omega_{\tau}\right)  \in R^{2}%
\pi_{\ast}\Omega_{\mathcal{Y}\left(  \text{M}\right)  \left/  \mathcal{M}%
\left(  \text{M}\right)  \right.  }^{1}%
\]
with the holomorphic section%
\[%
{\displaystyle\sum\limits_{k}}
C_{ij}^{k}\omega_{k}^{\ast}\in R^{1}\pi_{\ast}\Omega_{\mathcal{Y}\left(
\text{M}\right)  \left/  \mathcal{M}\left(  \text{M}\right)  \right.  }^{2},
\]
where $\omega_{k}^{\ast}\in R^{1}\pi_{\ast}\Omega_{\mathcal{Y}\left(
\text{M}\right)  \left/  \mathcal{M}\left(  \text{M}\right)  \right.  }%
^{2}\left\vert _{U}\right.  $ are defined by $\left(  \ref{hm3}\right)  $ and
are the Poincare dual of $\omega_{k}\in R^{2}\pi_{\ast}\Omega_{\mathcal{Y}%
\left(  \text{M}\right)  \left/  \mathcal{M}\left(  \text{M}\right)  \right.
}^{1}\left\vert _{U}\right.  $. Thus we have%
\[
\left\langle \phi_{l}\lrcorner\Omega_{\tau},\phi_{i}\wedge\phi_{j}%
\lrcorner\Omega_{\tau}\right\rangle =\left\langle \phi_{l}\lrcorner
\Omega_{\tau},%
{\displaystyle\sum\limits_{k}}
C_{ij}^{k}\omega_{k}^{\ast}\right\rangle =
\]
\begin{equation}
\left\langle \phi_{l}\lrcorner\Omega_{\tau},%
{\displaystyle\sum\limits_{k}}
C_{ij}^{k}\omega_{k}\right\rangle =%
{\displaystyle\sum\limits_{k}}
C_{ij}^{k}h_{lk}.\label{H8}%
\end{equation}
Combining $\left(  \ref{H7}\right)  $ and $\left(  \ref{H8}\right)  $ we get:
\[
\left\langle \phi_{l}\lrcorner\Omega_{\tau},%
{\displaystyle\sum\limits_{k}}
C_{ij}^{k}\omega_{k}\right\rangle =\left\langle \phi_{l}\lrcorner\Omega_{\tau
},\phi_{i}\wedge\phi_{j}\lrcorner\Omega_{\tau}\right\rangle =
\]%
\begin{equation}
-\sqrt{-1}%
{\displaystyle\int\limits_{\text{M}_{\tau_{0}}}}
(\Omega_{\tau}\wedge\left(  (\phi_{i}\wedge\phi_{j}\wedge\phi_{l}%
)\lrcorner\Omega_{\tau}\right)  =C_{ijl}.\label{hm1}%
\end{equation}
So we can conclude from $\left(  \ref{hk}\right)  $ and $\left(
\ref{hm3a}\right)  $ that
\begin{equation}
C_{ij}^{k}=%
{\displaystyle\sum\limits_{l=1}^{N}}
C_{ijl}h^{lk}.\label{hm1a}%
\end{equation}
Thus $\left(  \ref{hm1}\right)  $ and $\left(  \ref{hm1a}\right)  $ imply
$\left(  \ref{H1}\right)  $ and $\left(  \ref{h1}\right)  .$

Next we will prove $\left(  \ref{H2}\right)  .$ We can multiply the global
section $\eta_{\tau}$ of $\omega_{\mathcal{Y}\left(  \text{M}\right)
\text{/}\mathcal{M}\left(  \text{M}\right)  }$ defined in \cite{to1} by a
constant and so we can assume that at the point $\tau_{0}\in\mathcal{M}\left(
\text{M}\right)  $ we have
\begin{equation}
\eta_{\tau_{0}}=\Omega_{0}\text{ and }\eta_{\tau}=\lambda(\tau)\Omega_{\tau}.
\label{H7a}%
\end{equation}
From the definition of $\iota$ given by $\left(  \ref{Id}\right)  $ and
formula $\left(  \ref{H7}\right)  $, we conclude that $\left(  \ref{H2}%
\right)  $ holds. Thus Lemma \ref{flat} is proved. $\blacksquare$

\begin{lemma}
\label{flata}We have
\begin{equation}
\lbrack\vartheta_{i},\vartheta_{j}]=0. \label{H3}%
\end{equation}

\end{lemma}

\textbf{Proof: }The formula for\ $\vartheta_{i}$ given by $\left(
\ref{GM}\right)  $ implies%

\begin{equation}
\lbrack\vartheta_{i},\vartheta_{j}]\phi_{k}=\digamma^{-1}([\mathcal{D}%
_{i},\mathcal{D}_{j}](\phi_{k}\lrcorner\Omega_{\tau})). \label{H8a}%
\end{equation}
Since Gauss-Manin connection $\mathcal{D}$ is a flat connection then
\begin{equation}
\lbrack\mathcal{D}_{i},\mathcal{D}_{j}]=0. \label{H9}%
\end{equation}
Combining formula $\left(  \ref{H9}\right)  $ with $\left(  \ref{H8a}\right)
$ we get $[\vartheta_{i},\vartheta_{j}]\phi_{k}=0.$ Thus Lemma \ref{flata} is
proved. $\blacksquare$

\subsection{Relations with Frobenius Algebras}

One can use Lemma \ref{manin} to define an associative product on the tangent
bundle of the moduli space $\mathcal{M}\left(  \text{M}\right)  $ of three
dimensional CY manifolds as follows: Let $\{\phi_{i}\}$ be a basis of
T$_{\tau_{0},\mathcal{M}\left(  \text{M}\right)  }=H^{1}\left(  \text{M}%
,T_{\text{M}}^{1,0}\right)  ,$ then we define the product as:%
\begin{equation}
\phi_{i}\times\phi_{j}=\mathfrak{i}^{-1}\left(  \Pi\left(  \mathcal{D}%
_{i}(\mathfrak{i(}\phi_{j}))\right)  \right)  =F_{ijk}\phi_{k}. \label{Frob}%
\end{equation}

\begin{lemma}
\label{Fr}Let $C_{ij}^{k}$ be defined by $\left(  \ref{H7}\right)  ,$ then
\begin{equation}
F_{ijk}=\sqrt{-1}C_{ijk}. \label{Fr1}%
\end{equation}

\end{lemma}

\textbf{Proof:} Lemma \ref{Fr} follows directly from the formulas for
$F_{ijk}$ and $C_{ijk}.$ $\blacksquare$

\begin{corollary}
\label{FrI}The relations $\left(  \ref{H2}\right)  $ and $\left(
\ref{H3}\right)  $ shows that $F_{ijk}$ define a structure of a commutative
algebra on the tangent bundle of the moduli space $\mathcal{M}\left(
\text{M}\right)  $ of three dimensional CY manifolds.
\end{corollary}

\subsection{The Analogue of Cecotti-Vafa tt* Equations on $\mathcal{M}\left(
\text{M}\right)  $}

\begin{definition}
\label{csv}Let $\vartheta$ be the Higgs field defined by $\left(
\ref{GM}\right)  .$ Theorem \ref{mod} implies the existence of a global
holomorphic non vanishing section $\eta_{\tau}$ of the line bundle $\pi_{\ast
}\omega_{\mathcal{Y}\left(  \text{M}\right)  \left/  \mathcal{M}\left(
\text{M}\right)  \right.  }.$ Then $\eta_{\tau}$ defines a metric on the flat
line bundle $\pi_{\ast}\omega_{\mathcal{Y}\left(  \text{M}\right)  \left/
\mathcal{M}\left(  \text{M}\right)  \right.  }$ with curvature zero$.$ Let us
define the Weil-Petersson metric on%
\[
\mathcal{T}_{\mathcal{M}\left(  \text{M}\right)  }\thickapprox Hom\left(
\pi_{\ast}\omega_{\mathcal{Y}\left(  \text{M}\right)  \left/  \mathcal{M}%
\left(  \text{M}\right)  \right.  },R^{1}\pi_{\ast}\Omega_{\mathcal{Y}\left(
\text{M}\right)  \left/  \mathcal{M}\left(  \text{M}\right)  \right.  }%
^{2}\right)  \thickapprox
\]%
\[
\left(  \pi_{\ast}\omega_{\mathcal{Y}\left(  \text{M}\right)  \left/
\mathcal{M}\left(  \text{M}\right)  \right.  }\right)  ^{\ast}\otimes R^{1}%
\pi_{\ast}\Omega_{\mathcal{Y}\left(  \text{M}\right)  \left/  \mathcal{M}%
\left(  \text{M}\right)  \right.  }^{2}%
\]
as follows: Let
\[
\phi_{\tau}=\left(  \eta_{\tau}\right)  ^{\ast}\otimes\omega_{\tau}%
(1,1)\in\mathcal{T}_{\mathcal{M}\left(  \text{M}\right)  }\thickapprox\left(
\pi_{\ast}\omega_{\mathcal{Y}\left(  \text{M}\right)  \left/  \mathcal{M}%
\left(  \text{M}\right)  \right.  }\right)  ^{\ast}\otimes R^{1}\pi_{\ast
}\Omega_{\mathcal{Y}\left(  \text{M}\right)  \left/  \mathcal{M}\left(
\text{M}\right)  \right.  }^{2}.
\]
Let us define the function
\begin{equation}
\lambda(\tau):=\frac{\eta_{\tau}}{\Omega_{\tau}}. \label{eta}%
\end{equation}
Then%
\[
\left\Vert \phi_{\tau}\right\Vert ^{2}:=\left\vert \lambda(\tau)\right\vert
^{-2}\otimes\sqrt{-1}%
{\displaystyle\int\limits_{\text{M}}}
\omega_{\tau}(1,2)\wedge\overline{\omega_{\tau}(1,2)}=
\]%
\begin{equation}
\frac{\sqrt{-1}\left\langle \omega_{\tau}(1,2),\overline{\omega_{\tau}%
(1,2)}\right\rangle }{\left\vert \lambda(\tau)\right\vert ^{2}}. \label{fwp}%
\end{equation}
Let $\nabla$ be the standard connection of the metric defined by $\left(
\ref{fwp}\right)  .$ We will define the Cecotti-Hitchin-Vafa-Simpson (CHVS)
connection $D=D+\overline{D}$ on the tangent bundle of the moduli space
$\mathcal{M}\left(  \text{M}\right)  $ of three dimensional CY manifolds as
follows:%
\begin{equation}
D_{i}:\nabla_{i}+t\vartheta_{i}\text{ and }D_{\overline{j}}=\nabla
_{\overline{j}}+t^{-1}\overline{\vartheta}_{j}, \label{CSV1}%
\end{equation}
\textit{where } $t\in\mathbb{C}^{\ast}.$
\end{definition}

\begin{theorem}
\label{WP1}The curvature of the metric defined by $\left(  \ref{fwp}\right)  $
in the flat coordinates defined by Definition \ref{coor} is given by:%
\[
R_{i\overline{j},k\overline{l}}=-\sqrt{-1}%
{\displaystyle\int\limits_{\text{M}}}
\left(  (\phi_{i}\wedge\phi_{k})\lrcorner\Omega_{\text{M}}\right)
\wedge\left(  \overline{(\phi_{j}\wedge\phi_{l})\lrcorner\Omega_{\text{M}}%
}\right)  =
\]%
\begin{equation}
=-\sqrt{-1}\left\langle (\phi_{i}\wedge\phi_{k})\lrcorner\Omega_{\text{M}%
},(\phi_{j}\wedge\phi_{l})\lrcorner\Omega_{\text{M}}\right\rangle .
\label{cwp}%
\end{equation}

\end{theorem}

\textbf{Proof: }Since the metric on $(\pi_{\ast}\omega_{\mathcal{Y}\left(
\text{M}\right)  \left/  \mathcal{M}\left(  \text{M}\right)  \right.  }%
)^{\ast}$ is flat, then it has a zero curvature. The connection of the metric
defined by $\left(  \ref{fwp}\right)  $ will be
\[
\nabla:=\widetilde{\nabla}\otimes id\oplus id\otimes\nabla_{1},
\]
where $\widetilde{\nabla}$ is the flat connection on $\left(  \pi_{\ast}%
\omega_{\mathcal{Y}\left(  \text{M}\right)  \left/  \mathcal{M}\left(
\text{M}\right)  \right.  }\right)  ^{\ast}$ and $\nabla_{1}$ is the
connection on $R^{1}\pi_{\ast}\Omega_{\mathcal{Y}\left(  \text{M}\right)
\left/  \mathcal{M}\left(  \text{M}\right)  \right.  }^{2}.$ Thus we get that
for the curvature $\left[  \nabla,\nabla\right]  $ we have that $\left[
\nabla,\nabla\right]  =\left[  \nabla_{1},\nabla_{1}\right]  .$ So the
curvature of the metric defined by $\left(  \ref{fwp}\right)  $ is equal equal
to the curvature on $R^{1}\pi_{\ast}\Omega_{\mathcal{Y}\left(  \text{M}%
\right)  \left/  \mathcal{M}\left(  \text{M}\right)  \right.  }^{2}.$ Formula
$\left(  \ref{T2}\right)  $ implies%
\[
\sqrt{-1}\partial_{i}\overline{\partial}_{\overline{j}}\left(  \left\langle
\Omega_{\tau},\Omega_{\tau}\right\rangle \right)  =
\]%
\begin{equation}
\sqrt{-1}\left\langle \Omega_{0}\lrcorner\phi_{i},\Omega_{0}\lrcorner\phi
_{j}\right\rangle -\frac{\sqrt{-1}}{2}%
{\displaystyle\sum\limits_{k,l}}
\left\langle (\phi_{i}\wedge\phi_{k})\lrcorner\Omega_{0})),(\phi_{j}\wedge
\phi_{l})\lrcorner\Omega_{0}\right\rangle \tau^{k}\overline{\tau}%
^{l}+O(\left\vert \tau\right\vert ^{3}). \label{wp1}%
\end{equation}
Thus $\left(  \ref{wp1}\right)  $ implies $\left(  \ref{cwp}\right)  .$
Theorem \ref{WP1} is proved. $\blacksquare$ \ 

We will show that the Cecotti-Hitchin-Vafa-Simpson connection is flat.

\begin{theorem}
\label{csvflat}The connection $D$ defined by $\left(  \ref{CSV1}\right)  $ is
a flat one, i.e.:%
\begin{equation}
\lbrack D_{i},D_{j}]=[D_{i},D_{\overline{j}}]=[D_{\overline{i}},D_{\overline
{j}}]=0 \label{CSV2}%
\end{equation}
for all $0\leqq i,j\leqq N.$
\end{theorem}

\textbf{Proof: }First we will prove that $[D_{i},D_{j}]=0.$ The definition of
$D_{i}=\nabla_{i}+t\vartheta_{i}$ implies that
\begin{equation}
\lbrack D_{i},D_{j}]=[\nabla_{i},\nabla_{j}]+t[\nabla_{i},\vartheta_{j}%
]+t^{2}[\vartheta_{i},\vartheta_{j}]. \label{CSV3}%
\end{equation}
We know that $\triangledown_{i}$ is a Hermitian connection of the
Weil-Petersson metric defined by $\left(  \ref{fwp}\right)  $, which is a
K\"{a}hler metric and thus the (2,0) part of its curvature is zero. This
implies that $[\nabla_{i},\nabla_{j}]=0.$ From Lemma \ref{flat} we know that
$[\vartheta_{i},\vartheta_{j}]=0.$

\begin{lemma}
\label{Gl}We have
\begin{equation}
\lbrack\nabla_{i},\vartheta_{j}]=0. \label{CSVa}%
\end{equation}

\end{lemma}

\textbf{Proof: }In the flat coordinates $(\tau^{1},...,\tau^{N})$ at a fixed
point $\tau_{0}=0\in\mathcal{M}\left(  \text{M}\right)  $ of the moduli space
we have that $\nabla_{i}=\partial_{i}$. The definition of $\vartheta_{i}$
given by the formula $\left(  \ref{GM}\right)  $ and applied to
\begin{equation}
\phi_{k}(\tau):=\left(  \frac{\partial}{\partial\tau^{k}}\Omega_{\tau}\right)
\lrcorner\left(  \Omega_{\tau}\right)  ^{\ast},\text{ }k=1,..,N, \label{orth}%
\end{equation}
where $\Omega_{\tau}$ is defined by $\left(  \ref{T1}\right)  $ gives that we
have at the point $\tau_{0}=0:$%
\[
\nabla_{i}\left(  \vartheta_{j}\left(  \phi_{k}(\tau)\right)  \right)
\left\vert _{\tau=\tau_{0}}\right.  =
\]%
\begin{equation}
\vartheta_{i}\left(  \nabla_{j}\left(  \phi_{k}(\tau)\right)  \right)
=\left(  \left(  \frac{\partial^{2}}{\partial\tau^{i}\partial\tau^{j}}\left(
\frac{\partial}{\partial\tau^{k}}\Omega_{\tau}\right)  \right)  \lrcorner
\left(  \Omega_{\tau}\right)  ^{\ast}\right)  . \label{GSVb}%
\end{equation}
So $\left(  \ref{GSVb}\right)  $ implies $\left(  \ref{CSVa}\right)  .$ Lemma
\ref{Gl} is proved. $\blacksquare$

\begin{corollary}
\label{Gl1}$[D_{i},D_{j}]=[D_{\overline{i}},D_{\overline{j}}]=0.$
\end{corollary}

\begin{lemma}
\label{Gl2b}We have $[\nabla_{i}+t\vartheta_{i},\nabla_{\overline{j}}%
+t^{-1}\overline{\vartheta_{j}}]=0.$
\end{lemma}

\textbf{Proof: }We will identify $T_{0,\mathcal{M}\left(  \text{M}\right)  }$
with $H^{1}\left(  \text{M}_{0},\Omega_{\text{M}_{0}}^{2}\right)  $ as in
Proposition \ref{ident}. We assumed that $\{\phi_{i}\}$ given by $\left(
\ref{orth}\right)  $ is an orthonormal basis of the tangent space
$T_{0,\mathcal{K}}=H^{1}\left(  \text{M},T_{\text{M}}^{1,0}\right)  $ at at
point $0\in\mathcal{M}\left(  \text{M}\right)  $. We have in the flat
coordinates $(\tau^{1},...,\tau^{N})$ that $\phi_{k}(\tau)=\left(  \nabla
_{k}\Omega_{\tau})\lrcorner\eta_{\tau}^{-1}\right)  .$ We will need the
following Propositions

\begin{proposition}
\label{Gl2c}We have
\begin{equation}
\left\langle \nabla_{i}\phi_{k},\nabla_{j}\phi_{l}\right\rangle =\frac
{1}{\left\vert \lambda\left(  \tau\right)  \right\vert ^{2}}%
{\displaystyle\int\limits_{\text{M}}}
\left(  \left(  \phi_{i}\wedge\phi_{k}\right)  \lrcorner\Omega_{\tau}\right)
\wedge\overline{((\phi_{j}\wedge\phi_{l})\lrcorner\Omega_{\tau})},
\label{con2}%
\end{equation}
where $\left(  \phi_{i}\wedge\phi_{k}\right)  \lrcorner\Omega_{\tau}$ and
$\left(  \phi_{j}\wedge\phi_{l}\right)  \lrcorner\Omega_{\tau}$ are forms of
type $(1,2)$ on M$_{\tau}.$
\end{proposition}

\textbf{Proof: }Since $\phi_{i}(\tau)=\left(  \left(  \nabla_{i}\Omega_{\tau
}\right)  \lrcorner\eta_{\tau}^{-1}\right)  $ we get that%
\[
\left\langle \nabla_{i}\phi_{k},\nabla_{j}\phi_{l}\right\rangle \left\vert
_{\tau=0}\right.  =\left\langle \left(  \nabla_{i}\nabla_{k}\Omega_{\tau
}\right)  \lrcorner\eta_{\tau}^{-1},(\nabla_{j}\nabla_{l}\Omega_{\tau
})\lrcorner\eta_{\tau}^{-1}\right\rangle \left\vert _{\tau=0}\right.  .
\]
From $\left(  \ref{T3}\right)  $ we derive that
\[
\nabla_{i}\nabla_{k}\Omega_{\tau}\left\vert _{\tau=0}\right.  =\left(
\phi_{i}\wedge\phi_{k}\right)  \lrcorner\Omega_{0}\text{ and }\nabla_{j}%
\nabla_{l}\Omega_{\tau}\left\vert _{\tau=0}\right.  =\left(  \phi_{j}%
\wedge\phi_{l}\right)  \lrcorner\Omega_{0}%
\]
are form of type $(1,2).$ Thus we get
\[
\left\langle \nabla_{i}\phi_{k},\nabla_{j}\phi_{l}\right\rangle \left\vert
_{\tau=0}\right.  =\left(  \frac{1}{\left\vert \lambda\left(  \tau\right)
\right\vert ^{2}}%
{\displaystyle\int\limits_{\text{M}}}
\left(  \left(  \phi_{i}\wedge\phi_{k}\right)  \lrcorner\Omega_{\tau}\right)
\wedge\overline{((\phi_{j}\wedge\phi_{l})\lrcorner\Omega_{\tau})}\right)
\left\vert _{\tau=0}\right.  ,
\]
where $\left(  \phi_{i}\wedge\phi_{k}\right)  \lrcorner\Omega_{0}$ and
$\left(  \phi_{j}\wedge\phi_{l}\right)  \lrcorner\Omega_{0}$ are forms of type
$(1,2).$ Proposition \ref{Gl2c} is proved. $\blacksquare$

\begin{proposition}
\label{Gl2d}We have
\[
\vartheta_{i}(\phi_{k}(\tau))\left\vert _{\tau=0}\right.  =\left(  \left(
\Pi\left(  (\phi_{i}(\tau)\wedge\phi_{k}(\tau))\lrcorner\Omega_{\tau}\right)
\right)  \lrcorner\left(  \Omega_{\tau}\right)  ^{-1}\right)  \left\vert
_{\tau=0}\right.
\]
and%
\[
\left\langle \vartheta_{i}\phi_{k},\vartheta_{j}\phi_{l}\right\rangle
\left\vert _{\tau=0}\right.  =
\]%
\begin{equation}
\left(  \frac{1}{\left\vert \lambda\left(  \tau\right)  \right\vert ^{2}}%
{\displaystyle\int\limits_{\text{M}}}
\left(  \Pi\left(  (\phi_{i}(\tau)\wedge\phi_{k}(\tau))\lrcorner\Omega_{\tau
}\right)  \right)  \wedge\overline{(\Pi\left(  (\phi_{j}(\tau)\wedge\phi
_{l}(\tau))\lrcorner\Omega_{\tau})\right)  }\right)  \left\vert _{\tau
=0}\right.  , \label{con3}%
\end{equation}
where $\Pi\left(  \left(  \phi_{i}\wedge\phi_{k}\right)  \lrcorner\Omega
_{0}\right)  $ and $\Pi\left(  \left(  \phi_{j}\wedge\phi_{l}\right)
\lrcorner\Omega_{0}\right)  $ are forms of type $(2,1).$
\end{proposition}

\textbf{Proof: }It follows from the definition of $\vartheta$
\begin{equation}
\vartheta_{i}(\left(  \nabla_{k}\Omega_{\tau})\lrcorner\Omega_{\tau}%
^{-1}\right)  \left\vert _{\tau=0}\right.  =\left(  \left(  \Pi\left(
\phi_{i}\wedge\phi_{k}\lrcorner\Omega_{0}\right)  \right)  \lrcorner\Omega
_{0}^{-1}\right)  , \label{con3a}%
\end{equation}
where $\Pi\left(  \phi_{i}\wedge\phi_{k}\lrcorner\Omega_{\tau}\right)  $ is
the Poincare dual of the form $\left(  \phi_{i}\wedge\phi_{k}\lrcorner
\Omega_{0}\right)  $ of type $(1,2).$ Thus $\Pi\left(  \phi_{i}\wedge\phi
_{k}\lrcorner\Omega_{0}\right)  $ is a form of type $(2,1)$ and $\left(
\ref{con3a}\right)  $ implies $\left(  \ref{con3}\right)  .$ Proposition
\ref{Gl2d} is proved. $\blacksquare$

\begin{proposition}
\label{Gl2e}We have the following formula:%
\[
\left\langle ((\phi_{i}\wedge\phi_{k})\lrcorner\Omega_{0}),\overline
{((\phi_{j}\wedge\phi_{l})\lrcorner\Omega_{0})}\right\rangle =%
{\displaystyle\int\limits_{\text{M}}}
((\phi_{i}\wedge\phi_{k})\lrcorner\Omega_{0})\wedge\overline{((\phi_{j}%
\wedge\phi_{l})\lrcorner\Omega_{0})}=
\]%
\[
-\left\langle \Pi((\phi_{i}\wedge\phi_{k})\lrcorner\Omega_{0}),\overline
{\Pi((\phi_{j}\wedge\phi_{l})\lrcorner\Omega_{0})}\right\rangle =
\]%
\begin{equation}
=-%
{\displaystyle\int\limits_{\text{M}}}
\Pi((\phi_{i}\wedge\phi_{k})\lrcorner\Omega_{0})\wedge\overline{\Pi((\phi
_{j}\wedge\phi_{l})\lrcorner\Omega_{0})}. \label{com}%
\end{equation}

\end{proposition}

\textbf{Proof:} Since $\{\phi_{i}\}$ is an orthonormal basis in
$T_{0,\mathcal{M}\left(  \text{M}_{0}\right)  }$ with respect to the W.-P.
metric we get that $\{\phi_{i}\lrcorner\Omega_{0}\}$ is orthonormal basis in
$H^{1}\left(  \text{M}\right)  $ and $\{\overline{\phi_{i}\lrcorner\Omega_{0}%
}\}$ is an orthonormal basis in $H^{2}\left(  \text{M}_{0},\Omega
_{\text{M}_{0}}^{1}\right)  .~$\ Thus
\begin{equation}
\sqrt{-1}\left\langle \phi_{i}\lrcorner\Omega_{0},(\phi_{j}\lrcorner\Omega
_{0})\right\rangle =\sqrt{-1}%
{\displaystyle\int\limits_{\text{M}}}
(\phi_{i}\lrcorner\Omega_{0})\wedge\overline{(\phi_{j}\lrcorner\Omega_{0}%
)}=\delta_{i\overline{j}} \label{I}%
\end{equation}
and
\begin{equation}
\sqrt{-1}\left\langle \overline{\phi_{i}\lrcorner\Omega_{0}},\overline
{\phi_{j}\lrcorner\Omega_{0}}\right\rangle =\sqrt{-1}%
{\displaystyle\int\limits_{\text{M}}}
\overline{\phi_{i}\lrcorner\Omega_{0}}\wedge\phi_{j}\lrcorner\Omega
_{0}=-\delta_{i\overline{j}}. \label{II}%
\end{equation}
Let
\begin{equation}
(\phi_{i}\wedge\phi_{k})\lrcorner\Omega_{0}=%
{\displaystyle\sum\limits_{\nu=1}^{N}}
\alpha_{\nu}\overline{\phi_{\nu}\lrcorner\Omega_{0}},\text{ }(\phi_{j}%
\wedge\phi_{l})\lrcorner\Omega_{0}=%
{\displaystyle\sum\limits_{i=1}^{N}}
\beta_{\mu}\left(  \overline{\phi_{\mu}\lrcorner\Omega_{0}}\right)  ,
\label{III}%
\end{equation}
then
\begin{equation}
\Pi\left(  (\phi_{i}\wedge\phi_{k})\lrcorner\Omega_{0}\right)  =%
{\displaystyle\sum\limits_{\nu=1}^{N}}
\alpha_{\nu}\left(  \phi_{\nu}\lrcorner\Omega_{0}\right)  \text{ and }%
\Pi\left(  (\phi_{j}\wedge\phi_{l})\lrcorner\Omega_{0}\right)  =%
{\displaystyle\sum\limits_{i=1}^{N}}
\beta_{\mu}\left(  \phi_{\mu}\lrcorner\Omega_{0}\right)  . \label{IV}%
\end{equation}
Combining $\left(  \ref{I}\right)  ,$ $\left(  \ref{II}\right)  ,$ $\left(
\ref{III}\right)  $ and $\left(  \ref{IV}\right)  $ we get that
\begin{equation}
\sqrt{-1}\left\langle ((\phi_{i}\wedge\phi_{k})\lrcorner\Omega_{0}),((\phi
_{j}\wedge\phi_{l})\lrcorner\Omega_{0})\right\rangle =%
{\displaystyle\sum\limits_{\nu=1}^{N}}
\alpha_{\nu}\beta_{\nu} \label{V}%
\end{equation}
and
\begin{equation}
\sqrt{-1}\left\langle \Pi((\phi_{i}\wedge\phi_{k})\lrcorner\Omega_{0}%
),(\Pi(\phi_{j}\wedge\phi_{l})\lrcorner\Omega_{0})\right\rangle =-%
{\displaystyle\sum\limits_{\nu=1}^{N}}
\alpha_{\nu}\beta_{\nu}. \label{VI}%
\end{equation}
Thus $\left(  \ref{V}\right)  $ and $\left(  \ref{VI}\right)  $ imply
Proposition \ref{Gl2e}. $\blacksquare$

We have
\begin{equation}
\left\langle \left[  \nabla_{i}+t\vartheta_{i},\nabla_{\overline{j}}%
+t^{-1}\overline{\vartheta_{j}}\right]  \phi_{k},\phi_{l}\right\rangle
=\left\langle \left[  \nabla_{i},\nabla_{\overline{j}}\right]  \phi_{k}%
,\phi_{l}\right\rangle +\left\langle \left[  \vartheta_{i},\overline
{\vartheta_{j}}\right]  \phi_{k},\phi_{l}\right\rangle . \label{cwp1}%
\end{equation}
From $\left(  \ref{cwp}\right)  $ we derive that
\begin{equation}
R_{i\overline{j}.k\overline{l}}=\left\langle [\nabla_{i},\nabla_{\overline{j}%
}]\phi_{k},\phi_{l}\right\rangle =-\frac{\sqrt{-1}}{\left\vert \lambda
(0)\right\vert ^{2}}\left\langle (\phi_{i}\wedge\phi_{k})\lrcorner\Omega
_{0})),(\phi_{j}\wedge\phi_{l})\lrcorner\Omega_{0}\right\rangle . \label{cwp2}%
\end{equation}
$\left(  \ref{con3}\right)  $ implies
\begin{equation}
\left\langle \left[  \vartheta_{i},\overline{\vartheta_{j}}\right]  \phi
_{k},\phi_{l}\right\rangle =\frac{1}{\left\vert \lambda(0)\right\vert ^{2}%
}\left\langle \left(  \Pi\left(  (\phi_{i}\wedge\phi_{k})\lrcorner\Omega
_{0}\right)  \right)  ,\overline{(\Pi\left(  (\phi_{j}\wedge\phi_{l}%
)\lrcorner\Omega_{0})\right)  }\right\rangle . \label{cwp3}%
\end{equation}
Combining $\left(  \ref{cwp1}\right)  ,$ $\left(  \ref{cwp2}\right)  ,$
$\left(  \ref{cwp3}\right)  $ with $\left(  \ref{com}\right)  $ we conclude
that $[\nabla_{i}+t\vartheta_{i},\nabla_{\overline{j}}+t^{-1}\overline
{\vartheta_{j}}]=0.$ Theorem \ref{csvflat} is proved. $\blacksquare$

\begin{corollary}
\label{CHSV1}The connection constructed in Theorem \ref{csvflat} is a flat
$\mathbb{S}p(2h^{2,1},\mathbb{R)}$ connection on the tangent space of the
moduli space $\mathcal{M}\left(  \text{M}\right)  $ of three dimensional CY
manifolds. The imaginary form of the Weil-Petersson metric is a parallel form
with respect to the CHSV connection.
\end{corollary}

\textbf{Proof: }It is an well known fact that the imaginary part of the
Weil-Petersson metric $\omega_{\tau}(1,1)$ is parallel with respect to the
CHSV connection since it is the imaginary part of a K\"{a}hler metric. On the
other hand $\omega_{\tau}(1,1)$ is just the restriction of the intersection
form on $H^{3}\left(  \text{M}\right)  $ and so it is parallel with respect to
the Gauss-Manin connection and so with respect to the connection $\theta.$
From here Corollary \ref{CHSV1} follows directly. $\blacksquare$

\begin{remark}
\label{chvs} It is easy to see that Cecotti-Vafa tt* equations are exactly the
Hitchin-Simpson self duality equations studied in \ \cite{H87}, \cite{sim1}
and \cite{sim2}.
\end{remark}

We will call the connection that we constructed a Cecotti-Hitchin-Simpson-Vafa
connection and will denote it as a CHSV connection.

\section{Review of the Geometric Quantization}

\subsection{Basic Notions of ADW\ Geometric Quantization}

In this paragraph we are going to review the method of the geometric
quantization described in \cite{ADW} and \cite{W93}. We will consider a linear
space $\mathbf{W}\approxeq\mathbb{R}^{2n}$ with a constant symplectic structure%

\begin{equation}
\omega=\frac{1}{2}\omega_{ij}dt^{i}\wedge dt^{j}, \label{SY}%
\end{equation}
where $\omega_{ij}$ is a constant invertible matrix and the $x^{i}$ linear
coordinates on $\mathbb{R}^{2n}=\mathbf{W}.$ We will denote by $\omega^{-1}$
the matrix inverse to $\omega$ and obeying $\omega_{ij}\left(  \omega
^{-1}\right)  ^{jk}=\delta_{i}^{k}.$

\begin{definition}
\label{qlb}The \textquotedblright prequantum line bundle\textquotedblright\ is
a unitary line bundle $\mathcal{L}$ \ over $\mathbf{W}$ with a connection
whose curvature is $\sqrt{-1}\omega$. Up to an isomorphism, there is only one
such choice of $\mathcal{L}$. One can take $\mathcal{L}$ to be the trivial
unitary line bundle, with a connection given by the covariant derivatives%
\begin{equation}
\frac{D}{Dt^{i}}=\frac{\partial}{\partial t^{i}}+\frac{\sqrt{-1}}{2}%
\omega_{ij}t^{j}. \label{Wit1}%
\end{equation}

\end{definition}

\begin{definition}
\label{preqhs}We define the $\mathbf{L}^{2}$ norm of the sections of
$\mathcal{L}$ as follows; Let $h$ be a positive function on $\mathbf{W}$ which
define a metric on $\mathcal{L}$ and $dd^{c}\log h=\omega.$ Then we will say
\ that the $\mathbf{L}^{2}$ norm of a section $\phi$ of $\mathcal{L}$ is
defined as
\begin{equation}
\left\Vert \phi\right\Vert _{\mathbf{L}^{2}}^{2}=\left(  -1\right)
^{n(n-1)\left/  2\right.  }\left(  -\frac{\sqrt{-1}}{2}\right)  ^{n}%
{\displaystyle\int\limits_{\mathbf{W}}}
\exp\left(  -h\right)  \left\vert \phi\right\vert ^{2}dz^{1}\wedge...\wedge
dz^{n}\wedge\overline{dz^{1}\wedge...\wedge dz^{n}}. \label{Wit1a}%
\end{equation}
Then we define the\textquotedblright prequantum Hilbert
space\textquotedblright\ $\mathcal{H}_{0}$ as the Hilbert space that consists
of sections of $\mathcal{L}$ with a finite $\mathbf{L}^{2}$ norm.
\end{definition}

In order to define the quantum Hilbert space, we will introduce the notion of polarization.

\begin{definition}
\label{pol}We will define the polarization as a choice of a complex structure
$J$ on $\mathbf{W}$ with the following properties: \textbf{a.} $J$ is a
translation invariant, so it is defined by a constant matrix $J_{j}^{i}$ with
$J^{2}=-1$. \textbf{b.} The two-form $\omega$ is of type (1,1) with respect to
the complex structure $J.$ \textbf{c.} $J$ is positive in the sense that the
bilinear form g defined by g$(u,v)=\omega(u,Jv)$ is strictly positive.
\end{definition}

\begin{definition}
\label{qhs}Given such a complex structure $J$, we define the quantum Hilbert
space $H_{J}$ to be \ space of all holomorphic functions $\phi_{J}%
(z^{1},...,z^{n})$ on $\mathbf{W}$ with respect to the complex structure $J$
with a finite $\mathbf{L}^{2}$ norm.
\end{definition}

It is well known that the Heisenberg group of $\mathbf{W}$ has an irreducible
projective representation in $\mathcal{H}_{J}.$ Thus each such representation
of the Heisenberg group of $\mathbf{W}$ depends on the choice of the complex
structure $J$ of $\mathbf{W}.$ We want to construct a projectively flat
connection on the infinite dimensional dimensional vector space over the
parameter space of the complex structures of $\mathbf{W}.$ Construction of
such a connection enables one to identify all the irreducible projective
representation of the Heisenberg group in $\mathcal{H}_{J}.$

\subsection{Siegel Space}

We will introduce some notations following \cite{W93}. First of all, one has
the projection operators $\frac{1}{2}(1\mp\sqrt{-1}J)$ on \textbf{$W$}$^{1,0}$
and on $\overline{\mathbf{W}^{1,0}}=$\textbf{$W$}$^{0,1},$ where%
\[
\mathbf{W}^{1,0}:=\{u\in\mathbf{W}\otimes\mathbb{C}\left\vert Ju=\sqrt
{-1}u\right.  \}\text{ .}%
\]
We know that the vector space \textbf{$W$} with a complex structure $J$ can be
identified as a complex vector space canonically with the spaces \textbf{$W$%
}$^{1,0}$\ and \textbf{$W$}$^{0,1}$ by the maps%
\[
u\rightarrow\frac{1}{2}(1\mp\sqrt{-1}J)u.
\]
We need to write down explicitly in fixed coordinates these two
identifications. We will follow the notations from \cite{W93} in the above
identifications. For any vector $v=(...,v^{i},...)$, we denote by%
\[
v^{\underline{i}}=\frac{1}{2}(1-\sqrt{-1}J)_{j}^{i}v^{j}\text{ and
}v^{\overline{i}}=\frac{1}{2}(1+\sqrt{-1}J)_{j}^{i}v^{j}.
\]
For one forms $w=(...,w_{i},...)$ we have%
\[
w_{\underline{j}}=\frac{1}{2}(1-\sqrt{-1}J)_{j}^{i}w_{i}\text{ and
}w_{\overline{j}}=\frac{1}{2}(1+\sqrt{-1}J)_{j}^{i}w_{j}.
\]
Thus $J_{\underline{j}}^{\underline{i}}=\sqrt{-1}\delta_{\underline{j}%
}^{\underline{i}}$ and $J_{\overline{j}}^{\overline{i}}=\sqrt{-1}%
\delta_{\overline{j}}^{\overline{i}}.$ This means that the projections of
$J_{j}^{i}$ and $\delta_{j}^{i}$ on \textbf{$W$}$^{1,0}$ and \textbf{$W$%
}$^{0,1}$ are proportional.

Let $\mathfrak{Z}$ be the space of all $J$ obeying the conditions in
Definition \ref{pol}. Then $\mathfrak{Z}$ is a symmetric space, i.e.
$\mathfrak{Z}=\mathbb{S}p(2n,\mathbb{R)}/\mathbb{U}(n).$ $\mathfrak{Z}$ is
called Siegel space. It is a well known fact that we have the following
realization of $\mathfrak{Z}$ $=\mathbb{S}p(2n,\mathbb{R)}/\mathbb{U}(n)$ as a
tube domain;%
\[
\mathfrak{Z=}\{Z|Z\text{ }is\text{ }a\text{ }(n\times n)\text{ complex matrix
such that }Z^{t}=Z\text{ }and\text{ }\operatorname{Im}Z>0\}.
\]
$\mathfrak{Z}$ \ has a natural complex structure, defined as follows. The
condition $J^{2}=-1$ implies that for a first order variation $\delta J$ of
$J$, one must have%
\[
J\circ\delta J+\delta J\circ J=0
\]
This means that the non-zero projections of $\delta J$ on (1,0) vectors and
(0,1) vectors are $\delta J_{\overline{j}}^{\underline{i}}$ and $\delta
J_{\underline{j}}^{\overline{i}}$. We define the complex structure on
$\mathfrak{Z}$ \ by declaring $\delta J_{\overline{j}}^{\underline{i}}$ to be
of type $(1,0)$ and $\delta J_{\underline{j}}^{\overline{i}}$ to be of type
$(0,1)$. Notice that the projection of $\delta J$ on $(0,1)$ vectors is a map
from $(1,0)$ vectors to $(0,1)$ vectors. So the $(1,0)$ part of $\delta J$ is
the Beltrami differential as it was defined in \textbf{Section 2 }of this article.

\subsection{Construction of Witten Projective Connection}

Over $\mathfrak{Z}$ we introduce two Hilbert space bundles. One of them, say
$\mathcal{H}^{0}$, is the trivial bundle $\mathfrak{Z\times}\mathcal{H}_{0}$,
where $\mathcal{H}_{0}$ is the Hilbert space of all function $\psi(t^{i};J) $
with finite $\mathbf{L}^{2}$ norm. The definition of $\mathcal{H}_{0}$ is
independent of $J$. The second is the bundle $\mathcal{H}^{Q}$, whose fibre
over a point $J\in\mathfrak{Z}$ are functions $\psi(t^{i};J)$ of $t^{i}$,
holomorphic in the complex structure defined by $J,$ i.e. the following
equation holds:
\[
\frac{D}{Dt^{\overline{i}}}\psi(t^{i};J)=0.
\]
This equation has a dependence on $J$ coming from the projection operators
used in defining $t^{\overline{i}}.$

A connection on the bundle $\mathcal{H}^{0}$ restricts to a connection on
$\mathcal{H}^{Q}$ if and only if its commutator with $D_{\overline{i}}$ is a
linear combination of the $D_{\overline{j}}$. Since $\mathcal{H}^{0}$ is
defined as a product bundle $\mathfrak{Z\times}\mathcal{H}_{0}$, there is a
trivial connection $\delta$ on this bundle%

\begin{equation}
\delta:=\sum_{i,j}dJ_{j}^{i}\frac{\partial}{\partial J_{j}^{i}}. \label{Wit2}%
\end{equation}
We can expand $\delta$ in $(1,0)$ and $(0,1)$ pieces, $\delta=\delta
^{(1,0)}+\delta^{(0,1)}$, with%

\begin{equation}
\delta^{(1,0)}=\sum_{i,j}dJ_{\overline{j}}^{\underline{i}}\frac{\partial
}{\partial J_{\overline{j}}^{\underline{i}}}. \label{Wit3}%
\end{equation}
Unfortunately, as it was shown in \cite{W93} the commutator of $\delta
^{(1,0)}$ with $D_{\overline{i}}$ is not a linear combination of the
$D_{\overline{j}}$. So one needs to modify the trivial connection so that its
commutator with $D_{\overline{i}}$ will be a linear combination of the
$D_{\overline{j}}$

\begin{definition}
\label{pc}Witten defined in \cite{W93} the following connection $\nabla$ on
the bundle $\mathcal{H}^{0}=\mathfrak{Z\times}\mathcal{H}_{0}\rightarrow
\mathfrak{Z:}$%
\begin{equation}
\mathfrak{\nabla}^{(1,0)}:=\delta^{(1,0)}-\frac{1}{4}(dJ\omega^{-1}%
)^{\underline{i}\underline{j}}\frac{D}{Dt^{\underline{i}}}\frac{D}%
{Dt^{\underline{j}}}\text{ and }\mathfrak{\nabla}^{(0,1)}:=\delta^{(0,1)}.
\label{Wit4}%
\end{equation}

\end{definition}

Witten proved the following theorem in \cite{W93}:

\begin{theorem}
\textbf{A. }The connection $\nabla$ descends to a connection on $\mathcal{H}%
^{Q}$. \textbf{B.} The curvature of the connection $\nabla$ on $\mathcal{H}%
^{Q}$ is of type $(1,1)$ and it is equal to%
\begin{equation}
\lbrack\nabla^{(0,1)},\nabla^{(1,0)}]=\left[  \delta^{(0,1)},-\frac{1}%
{4}(dJ\omega^{-1})^{\underline{i}\underline{j}}\frac{D}{Dt^{\underline{i}}%
}\frac{D}{Dt^{\underline{j}}}\right]  =-\frac{1}{8}dJ_{\overline{k}%
}^{\underline{i}}dJ_{\underline{j}}^{\overline{k}}. \label{Wit5}%
\end{equation}

\end{theorem}

This theorem shows that the curvature of $\nabla$ is not zero, even when it is
restricted to $\mathcal{H}^{Q}$. The curvature is a c-number, that is, it
depends only on $J$, and not on the variables $t^{i}$ that are being
quantized. The fact that the curvature of $\ \nabla$ is a c-number means that
parallel transport by $\nabla$ is unique up to a scalar factor which,
moreover, it is of modulus 1 since the curvature is real or more fundamentally
since $\nabla$ is unitary. So up to this factor one can identify the various
$\mathcal{H}_{J}$ 's, and regard them as a different realization of the
quantum Hilbert space $\mathcal{H}$.

\section{Geometric Quantization of Moduli Space of CY Manifolds and Quantum
Background Independence}

\subsection{Symplectic Structures and the Flat Coordinates}

Our goal is to quantize geometrically the cotangent bundle of the moduli space
$\mathcal{M}\left(  \text{M}\right)  $. This means to define a flat
$\mathbb{S}p(2h^{2,1},\mathbb{R})$ connection on the cotangent bundle of
$\mathcal{M}\left(  \text{M}\right)  $. We have done this in the
\textbf{Section 3. }Then will construct the prequantum line bundle. Once the
prequatum line bundle is constructed we will compute explicitly the projective
connection defined in \cite{ADW} and \cite{W93}. After the geometric
quantization is done we will derive BCOV holomorphic anomaly equations
established in \cite{BCOV} as the projective flat connections on the Hilbert
bundle $\mathcal{H}^{Q}$. To perform these computations, it is important to
fix first the flat symplectic structure on the tangent bundle $\mathcal{T}%
_{\mathcal{M}\left(  \text{M}\right)  }$, then the local coordinates on
$\mathcal{M}\left(  \text{M}\right)  $ that describe the change of the complex
structures on the CY manifold M$_{\tau}$ and the linear coordinates in
\textbf{$W$}$_{\tau}^{\ast}=\Omega_{\tau,\mathcal{M}\left(  \text{M}\right)
}^{1}.$

\begin{remark}
\label{0}We will define the symplectic form $\omega_{1}(\tau)$ on the vector
bundle
\[
R^{1}\pi_{\ast}\Omega_{\mathcal{Y}\left(  \text{M}\right)  \left/
\mathcal{M}\left(  \text{M}\right)  \right.  }^{2}\rightarrow\mathcal{M}%
\text{$\left(  \text{M}\right)  $}%
\]
as the imaginary part of the metric%
\begin{equation}
g_{i,\overline{j}}:=\left\langle \frac{\partial}{\partial\tau^{i}}\Omega
_{\tau},\frac{\partial}{\partial\tau^{j}}\Omega_{\tau}\right\rangle =\sqrt{-1}%
{\displaystyle\int\limits_{\text{M}}}
\frac{\partial}{\partial\tau^{i}}\Omega_{\tau}\wedge\overline{\frac{\partial
}{\partial\tau^{j}}\Omega_{\tau}}\text{ .} \label{Wit6}%
\end{equation}
\textit{The symplectic structure on the tangent space }$\mathcal{T}%
_{\mathcal{M}\left(  \text{M}\right)  }$\textit{\ is defined by the imaginary
part of the Weil-Petersson metric, i.e.:}
\begin{equation}
\omega(\tau)=\operatorname{Im}G\left\vert _{T_{\tau,\mathcal{K}}}\right.  ,
\label{Wit14}%
\end{equation}
\textit{where }%
\begin{equation}
G_{i,\overline{j}}:=\left\langle \left(  \frac{\partial}{\partial\tau^{i}%
}\Omega_{\tau}\right)  \lrcorner\eta_{\tau}^{-1},\left(  \frac{\partial
}{\partial\tau^{j}}\Omega_{\tau}\right)  \lrcorner\eta_{\tau}^{-1}%
\right\rangle _{WP}=\frac{\sqrt{-1}}{\left\Vert \eta_{\tau}\right\Vert ^{2}}%
{\displaystyle\int\limits_{\text{M}}}
\frac{\partial}{\partial\tau^{i}}\Omega_{\tau}\wedge\overline{\frac{\partial
}{\partial\tau^{j}}\Omega_{\tau}}\text{ .} \label{Wit14a}%
\end{equation}
\textit{Thus} $\left(  \ref{Wit14a}\right)  $\textit{\ implies the following
relations between the two symplectic structures: }
\begin{equation}
g_{i,\overline{j}}=e^{-K}G_{i,\overline{j}}, \label{Wit12a}%
\end{equation}
\textit{where}%
\begin{equation}
\left\Vert \eta_{\tau}\right\Vert ^{-2}=-\sqrt{-1}%
{\displaystyle\int\limits_{\text{M}}}
\eta_{\tau}\wedge\overline{\eta_{\tau}}=\exp(-K). \label{Wit12b}%
\end{equation}
\textit{We know from Remark \ref{chvs} that the forms} $\eta_{\tau}$
\textit{are parallel with respect to the Cecotti-Hitchin-Simpson-Vafa
connection. Thus we can identified the symplectic structures of the tangent
bundle by the flat} $\mathbb{S}p(2h^{2},\mathbb{R)}$
\textit{Cecotti-Hitchin-Simpson-Vafa connection.}
\end{remark}

\begin{remark}
\label{1} We know that if we fix a basis of the orthogonal vectors\textit{\ }%
$\{\phi_{i,\tau}\}$ in the holomorphic tangent space\textit{\ }%
\[
T_{\tau,\mathcal{K}}=H^{1}\left(  \text{M}_{\tau},T_{\text{M}_{gt}}%
^{1,0}\right)
\]
then we obtain a linear coordinate system $(\tau^{1},...,\tau^{N})$ in the
dual $\Omega_{\tau,\mathcal{K}}^{1}$ of $T_{\tau,\mathcal{K}}.$ $\Omega
_{\tau,\mathcal{K}}^{1}$ can be canonically identified with $T_{\tau
,\mathcal{K}}$ by using the parallel symplectic form. According to the results
obtained in \cite{to89} then linear coordinate system $(\tau^{1},...,\tau
^{N})$ defined by the choice of the orthogonal vectors\textit{\ }%
$\{\phi_{i,\tau_{0}}\}$ in the holomorphic tangent space\textit{\ }%
\[
T_{\tau_{\tau},\mathcal{K}}=H^{1}\left(  \text{M}_{\tau_{0}},T_{\text{M}%
_{\tau_{0}}}^{1,0}\right)
\]
defines the same local coordinate system $(\tau^{1},...,\tau^{N})$ in
$\mathcal{K}$ since
\[
\mathcal{K}\subset H^{1}\left(  \text{M,}T_{\text{M}_{\tau_{0}}}^{1,0}\right)
.
\]
See Theorem \ref{tod1}.
\end{remark}

\begin{remark}
\label{1a}In \cite{to89} we construct in a canonical way a family of
holomorphic form $\Omega_{\tau}$ given by $\left(  \ref{T1}\right)  $. Then
the family of holomorphic form $\Omega_{\tau}$ defines a choice of a basis
\begin{equation}
\{\phi_{i,\tau}\lrcorner\Omega_{\tau}\}\text{ and }\{\overline{(\phi_{i,\tau
}\lrcorner\Omega_{\tau})}\} \label{De}%
\end{equation}
of the complexified tangent space: $H^{1}\left(  \text{M,}\Omega_{\text{M}%
}^{2}\right)  \oplus H^{2}\left(  \text{M,}\Omega_{\text{M}}^{1}\right)  .$
Thus the flat local coordinate system $(\tau^{1},...,\tau^{N})$ is defined by
a choice of the basis in the holomorphic tangent space\textit{\ }%
$T_{\tau,\mathcal{K}}$ defines in a canonical way a coordinate system on
$H^{1}\left(  \text{M}\right)  $ which is the same as $(\tau^{1},...,\tau
^{N})$ in case we choose the basis $\left(  \ref{De}\right)  .$ We identified
the dual of $H^{1}\left(  \text{M,}\right)  $ with $H^{1}\left(
\text{M}\right)  $ by using the parallel symplectic form.
\end{remark}

\begin{remark}
\label{1b}The complex structure $J_{\tau}$ on $T_{\tau,\mathcal{M}\left(
\text{M}\right)  }$ is defined by the complex structure on \ $H^{1}\left(
\text{M,}T_{\text{M}}^{1,0}\right)  .$ Indeed by local deformation theory we
have
\[
T_{\tau,\mathcal{M}\left(  \text{M}\right)  }^{1,0}=H^{1}\left(
\text{M,}T_{\text{M}}^{1,0}\right)  .
\]
On CY manifold $H^{1}($M$_{\tau},T_{\tau}^{1,0})$ can be identified with
$H^{1}\left(  \text{M,}\Omega_{\text{M}}^{2}\right)  $ by contraction with the
non zero holomorphic form $\Omega_{\tau}.$ Thus the complex structure operator
$J_{\tau}$ on the CY manifolds acts on
\[
T_{\tau,\mathcal{M}\left(  \text{M}\right)  }^{1,0}=H^{1}\left(
\text{M,}T_{\text{M}}^{1,0}\right)
\]
in a natural way as follows
\[
J\left(  (dz^{i})\wedge(J\left(  dz^{j}\right)  \wedge(J\left(  \overline
{dz}^{k}\right)  \right)  =\sqrt{-1}(dz^{i}\wedge dz^{j}\wedge\overline
{dz}^{k})
\]
since $J(dz^{i})=\sqrt{-1}dz^{i}.$
\end{remark}

\begin{remark}
\label{2}We are using two different identifications of $H^{1}($M$_{\tau
},T_{\text{M}_{\tau}}^{1,0})$ with $H^{1}\left(  \text{M,}\Omega_{\text{M}%
}^{2}\right)  $ by using using the contraction with the families of the two
holomorphic $3-$forms $\eta_{\tau}$ and $\Omega_{\tau}.$ Recall that
$\eta_{\tau}$ was defined globally by Theorem \ref{BGS1}. We obtained two
coordinate systems $(\tau^{1},...,\tau^{N})$ and $(t^{1},...,t^{N})$ on
$H^{1}\left(  \text{M,}\Omega_{\text{M}}^{2}\right)  .$ The relation between
them is given by the relations between $\Omega_{\tau}$ and $\eta_{\tau},$ i.e.
the two coordinate systems are proportional with the coefficients of
proportionality $\lambda(\tau),$ where $\Omega_{\tau}=\lambda(\tau)\eta_{\tau
}.$
\end{remark}

Next we will define the prequantum line bundles over $R^{1}\pi_{\ast}%
\Omega_{\mathcal{Y}\left(  \text{M}\right)  \text{/}\mathcal{M}\left(
\text{M}\right)  }^{2}$ and $\mathcal{T}_{\mathcal{M}\left(  \text{M}\right)
}$.

\begin{theorem}
\label{W!!!!}The prequantum line bundles on the tangent bundle $\mathcal{T}%
_{\mathcal{M}\left(  \text{M}\right)  }$ is $p^{\ast}\left(  \pi_{\ast}\left(
\omega_{\mathcal{Y}(\text{M)/}\mathcal{M}(\text{M)}}^{\ast}\right)  \right)
,$ where $p:\mathcal{T}_{\mathcal{M}\left(  \text{M}\right)  }\rightarrow
\mathcal{M}\left(  \text{M}\right)  ,$ $\omega_{\mathcal{Y}(\text{M)/}%
\mathcal{M}(\text{M)}}$ is the relative dualizing line bundle and the metric
on its fibre is defined by the $L^{2}$ norm of the holomorphic three form. The
Chern class of $\pi^{\ast}\left(  \omega_{\mathcal{Y}(\text{M) /}%
\mathcal{M}(\text{M)}}^{\ast}\right)  $ is given by the restriction of
\[
\omega(\tau):=\{\text{Imaginary Part of the }W-P\text{ metric}\}
\]
on $T_{\tau,\mathcal{M}\left(  \text{M}\right)  }.$
\end{theorem}

\textbf{Proof: }We proved in \cite{to89} that the natural metric
\[
\Vert\Omega_{\tau}\Vert^{2}=-\sqrt{-1}%
{\displaystyle\int\limits_{\text{M}}}
\Omega_{\tau}\wedge\overline{\Omega_{\tau}}%
\]
on $\omega_{\mathcal{Y}\text{ }\left(  \text{M}\right)  \text{ /}%
\mathcal{M}\left(  \text{M}\right)  }$ is such that its Chern form is
\[
\sqrt{-1}\frac{\partial^{2}}{\partial\tau^{i}\partial\overline{\tau^{j}}}%
\log(\Vert\Omega_{\tau}\Vert^{2})=-\omega(\tau).
\]
This shows that $\mathcal{L}=\pi^{\ast}(\omega_{\mathcal{Y}(\text{M)
/}\mathcal{M(}\text{M)}}^{\ast}\mathcal{)}$ is the prequantum line bundle on
the tangent bundle $\mathcal{T}_{\mathcal{M}\left(  \text{M}\right)  }$. This
proves Theorem \ref{W!!!!}. $\blacksquare$

\begin{corollary}
\label{W!a}The prequantum line bundles on the bundle $R^{1}\pi_{\ast}%
\Omega_{\mathcal{Y}\left(  \text{M}\right)  \left/  \mathcal{M}\left(
\text{M}\right)  \right.  }^{2}$ is $p^{\ast}\left(  \pi_{\ast}\left(
\omega_{\mathcal{Y}\left(  \text{M}\right)  \text{ /}\mathcal{M}\left(
\text{M}\right)  }^{\ast}\right)  \right)  ,$ where%
\[
p:R^{1}\pi_{\ast}\Omega_{\mathcal{Y}\left(  \text{M}\right)  \left/
\mathcal{M}\left(  \text{M}\right)  \right.  }^{2}\rightarrow\mathcal{M}%
\left(  \text{M}\right)  ,\omega_{\mathcal{Y}\left(  \text{M}\right)  \left/
\mathcal{M}\left(  \text{M}\right)  \right.  }%
\]
is the relative dualizing line bundle and the metric on its fibre is defined
by the $L^{2}$ norm of the holomorphic three form. The Chern class of
$\pi^{\ast}\left(  \omega_{\mathcal{Y}\left(  \text{M}\right)  \left/
\mathcal{M}\left(  \text{M}\right)  \right.  }^{\ast}\right)  $ is given by
the restriction of $\omega_{1}(\tau)$ on \textbf{$W$}$_{\tau}=H^{1}\left(
\text{M,}T_{\text{M}}^{1.0}\right)  .$
\end{corollary}

\subsection{The Geometric Quantization of $R^{1}\pi_{\ast}\Omega
_{\mathcal{Y}\left(  \text{M}\right)  \left/  \mathcal{M}\left(
\text{M}\right)  \right.  }^{2}$}

We will compute the projective flat connection on the Hilbert space bundle
$\mathcal{H}^{Q}$ over the vector bundle $R^{1}\pi_{\ast}\Omega_{\mathcal{Y}%
\left(  \text{M}\right)  \text{/}\mathcal{M}\left(  \text{M}\right)  }^{2}.$
The prequantum line bundle is the line bundle $p^{\ast}\left(  \omega
_{\mathcal{Y}\left(  \text{M}\right)  /\mathcal{M}\left(  \text{M}\right)
}^{\ast}\right)  $ over $R^{1}\pi_{\ast}\Omega_{\mathcal{Y}\left(
\text{M}\right)  \mathcal{M}\left(  \text{M}\right)  }^{2}.$

\begin{theorem}
\label{expl}\textit{The matrix of the operator }$(dJ)$ \textit{in the basis}
\[
\left\{  \frac{\partial}{\partial\tau^{i}}\Omega_{\tau},\text{ }%
i=1,...,N\right\}
\]
\textit{is given on each fibre }of $R^{1}\pi_{\ast}\Omega_{\mathcal{Y}\left(
\text{M}\right)  \left/  \mathcal{M}\left(  \text{M}\right)  \right.  }^{2}%
$\textit{\ by},%
\begin{equation}
(dJ)_{\underline{a}}^{\overline{b}}=\sum_{c,d}C_{\underline{a}\underline
{c}\underline{d}}g^{d,\overline{b}}. \label{Wit12}%
\end{equation}

\end{theorem}

\textbf{Proof:} The proof of Theorem \ref{expl} follows directly from the
following Lemma:

\begin{lemma}
\label{W!}We have
\begin{equation}
\left(  \frac{\partial}{\partial\tau^{i}}J\right)  \left(  \frac{\partial
}{\partial\tau^{j}}\Omega_{\tau}\right)  =\sum_{k,l=1}^{N}C_{ijk}%
g^{k,\overline{l}}\overline{\frac{\partial}{\partial\tau^{l}}\Omega_{\tau}%
},\text{ } \label{Wit7a}%
\end{equation}
where%
\[
C_{ijk}=-\sqrt{-1}%
{\displaystyle\int\limits_{\text{M}}}
\left(  \frac{\partial^{2}}{\partial\tau^{i}\partial\tau^{j}}\Omega_{\tau
}\right)  \wedge\left(  \frac{\partial}{\partial\tau^{k}}\Omega_{\tau}\right)
\left\vert _{\tau=0}\right.  =
\]%
\begin{equation}
\sqrt{-1}%
{\displaystyle\int\limits_{\text{M}}}
\left(  \Omega_{\tau}\right)  \wedge\left(  \frac{\partial^{3}}{\partial
\tau^{i}\partial\tau^{j}\partial\tau^{k}}\Omega_{\tau}\right)  \left\vert
_{\tau=0}\right.  . \label{Wit7}%
\end{equation}

\end{lemma}

The idea of \ the proof of formula $\left(  \ref{Wit7a}\right)  $ is the
following one; we know that in the basis
\[
\left\{  \frac{\partial}{\partial\tau^{i}}\Omega_{\tau},\text{ }%
i=1,...,N\right\}
\]
of $H^{1}($M$_{\tau},\Omega_{\tau}^{2})$ the complex structure operator is
given by the matrix:%
\begin{equation}
\left(
\begin{array}
[c]{ll}%
\sqrt{-1}I_{h^{1,2}} & 0\\
0 & -\sqrt{-1}I_{h^{1,2}}%
\end{array}
\right)  . \label{Wit7b}%
\end{equation}
Thus $\left(  \ref{Wit7b}\right)  $ implies
\[
\left(  \frac{\partial}{\partial\tau^{i}}J\right)  \left(  \frac{\partial
}{\partial\tau^{j}}\Omega_{\tau}\right)  \left\vert _{\tau=0}\right.
=J\left(  \frac{\partial^{2}\Omega_{\tau}}{\partial\tau^{i}\partial\tau^{j}%
}\right)  \left\vert _{\tau=0}\right.  .
\]
So we need to compute the expression of the vectors $\left\{  \frac
{\partial^{2}}{\partial\tau^{i}\partial\tau^{j}}\Omega_{\tau},\text{
}i,j=1,...,N\right\}  $ as a linear combinations of%
\begin{equation}
\left\{  \overline{\frac{\partial}{\partial\tau^{i}}\Omega_{\tau}},\text{
}i=1,...,N\right\}  . \label{Wit7c}%
\end{equation}
We know from\ \cite{to89} that$\left\{  \frac{\partial^{2}\Omega_{\tau}%
}{\partial\tau^{i}\partial\tau^{j}}\right\}  \in H^{2}\left(  \text{M}\right)
$ where $\tau=(\tau^{1},...,\tau^{N})$ are the flat coordinates$.$ Therefore
if we express the vectors
\[
\left\{  \frac{\partial^{2}\Omega_{\tau}}{\partial\tau^{i}\partial\tau^{j}%
}\right\}  \in H^{2}\left(  \text{M}_{\tau},\Omega_{\text{M}_{\tau}}%
^{1}\right)
\]
as linear combination of the basis $\left(  \ref{Wit7c}\right)  ,$ we will get
explicitly the matrices of the operators $\frac{\partial}{\partial\tau^{i}%
}J_{\tau=0},$ $i=0,...,N.$ From the explicit formulas of the operators
$\frac{\partial}{\partial\tau^{i}}J\left\vert _{\tau=0}\right.  $ we will get
the formula $\left(  \ref{Wit7a}\right)  $.

\textbf{Proof: }We know that%

\begin{equation}
-\sqrt{-1}%
{\displaystyle\int\limits_{\text{M}}}
\left(  \frac{\partial}{\partial\tau^{i}}\Omega_{\tau}\right)  \wedge\left(
\overline{\frac{\partial}{\partial\tau^{j}}\Omega_{\tau}}\right)  \left\vert
_{\tau=0}\right.  =\delta_{i,\overline{j}}. \label{Wit8}%
\end{equation}
By using the expression $\left(  \ref{T1}\right)  $ for $\Omega_{\tau}$ and
the natural $L^{2}$ metric on $H^{2}\left(  \text{M}\right)  $ we get the
following expression of the vector $\left(  \frac{\partial^{2}}{\partial
\tau^{i}\partial\tau^{j}}\Omega_{\tau}\right)  \left\vert _{\tau=0}\right.  $
in the orthogonal basis
\[
\left\{  \overline{\frac{\partial}{\partial\tau^{k}}\Omega_{\tau}}\left\vert
_{\tau=0}\right.  ,\text{ }k=1,...,N\right\}
\]
of $H^{2,1}\left(  \text{M}\right)  :$%

\[
\left(  \frac{\partial^{2}}{\partial\tau^{i}\partial\tau^{j}}\Omega_{\tau
}\right)  |_{\tau=0}=-\sqrt{-1}\sum_{k=1}^{N}\left(  \left\langle
\frac{\partial^{2}}{\partial\tau^{i}\partial\tau^{j}}\Omega_{\tau}%
,\overline{\frac{\partial}{\partial\tau^{k}}\Omega_{\tau}}\right\rangle
\right)  \left(  \overline{\frac{\partial}{\partial\tau^{k}}\Omega_{\tau}%
}\right)  \left\vert _{\tau=0}\right.  =
\]%
\[
-\sqrt{-1}\sum_{k=1}^{N}\left(
{\displaystyle\int\limits_{\text{M}}}
\left(  \frac{\partial^{2}}{\partial\tau^{i}\partial\tau^{j}}\Omega_{\tau
}\right)  \wedge\left(  \frac{\partial}{\partial\tau^{k}}\Omega_{\tau}\right)
\right)  \left(  \overline{\frac{\partial}{\partial\tau^{k}}\Omega_{\tau}%
}\right)  \left\vert _{\tau=0}\right.  =
\]

\begin{equation}
-\sqrt{-1}\sum_{k=1}^{N}C_{ijk}g^{k,\overline{l}}\left(  \overline
{\frac{\partial}{\partial\tau^{l}}\Omega_{\tau}}\right)  \left\vert _{\tau
=0}\right.  \label{Wit9}%
\end{equation}
Formula $\left(  \ref{Wit9}\right)  $ implies that for any $\tau\in
\mathcal{K}$ we have%

\begin{equation}
\frac{\partial^{2}}{\partial\tau^{i}\partial\tau^{j}}\Omega_{\tau}=-\sqrt
{-1}\sum_{k=1}^{N}C_{ijk}g^{k,\overline{l}}\left(  \overline{\frac{\partial
}{\partial\tau^{l}}\Omega_{\tau}}\right)  . \label{Wi9}%
\end{equation}
We know that%

\begin{equation}
J\left(  \frac{\partial}{\partial\tau^{i}}\Omega_{\tau}\right)  =\sqrt
{-1}\left(  \frac{\partial}{\partial\tau^{i}}\Omega_{\tau}\right)  .
\label{Wit10}%
\end{equation}
Combining $\left(  \ref{Wi9}\right)  $ and $\left(  \ref{Wit10}\right)  $ we
get that%
\[
\frac{\partial}{\partial\tau^{j}}\left(  J\left(  \frac{\partial}{\partial
\tau^{i}}\Omega_{\tau}\right)  \right)  =\frac{\partial}{\partial\tau^{j}%
}\left(  \sqrt{-1}\left(  \frac{\partial}{\partial z^{i}}\Omega_{\tau}\right)
\right)  =
\]

\begin{equation}
\sqrt{-1}\left(  \frac{\partial^{2}}{\partial\tau^{j}\partial z^{i}}%
\Omega_{\tau}\right)  =\sum_{k=1}^{N}C_{ijk}g^{k,\overline{l}}\left(
\overline{\frac{\partial}{\partial\tau^{l}}\Omega_{\tau}}\right)  .
\label{Wit11}%
\end{equation}
Lemma \ref{W!} is proved. $\blacksquare$

Lemma \ref{W!} implies directly Theorem \ref{expl}. $\blacksquare$

\begin{corollary}
\label{expla}The projective flat connection on the Hilbert space bundle
$\mathcal{H}^{Q}$ over the vector bundle $R^{1}\pi_{\ast}\Omega_{\mathcal{X}%
\text{/}\mathcal{K}}^{2}$ is given by%
\begin{equation}
\frac{\partial}{\partial\tau^{a}}+\frac{\sqrt{-1}}{4}\sum\overline{C_{abc}%
}g^{\overline{b},i}g^{\overline{c},j}\frac{D}{Dt^{i}}\frac{D}{Dt^{j}}.
\label{Wit17b}%
\end{equation}

\end{corollary}

\textbf{Proof: }According to Theorem \ref{expl} on $R^{1}\pi_{\ast}%
\Omega_{\mathcal{X}\text{/}\mathcal{K}}^{2}$ we have
\begin{equation}
\left(  dJ\omega_{1}^{-1}\right)  ^{ij}D_{i}D_{j}=-\sqrt{-1}\sum_{a=1}%
^{N}\overline{C_{abc}}g^{\overline{b},i}g^{\overline{c},j}D_{i}D_{j}.
\label{Wit17}%
\end{equation}
Thus the projective connection on $R^{1}\pi_{\ast}\Omega_{\mathcal{X}%
\text{/}\mathcal{K}}^{2}$ is given by $\left(  \ref{Wit17b}\right)  $.
$\blacksquare$

\subsection{Computations on the Tangent Bundle of the Moduli Space}

To quantize geometrically the tangent bundle $\mathcal{T}_{\mathcal{M}\left(
\text{M}\right)  }$ on the moduli space $\mathcal{M}\left(  \text{M}\right)  $
means to compute explicitly the prequantum line bundle and then projective
connection on the Hilbert space bundle $\mathcal{H}^{Q}$ related to the
prequantum line bundle on the tangent bundle $\mathcal{T}_{\mathcal{M}\left(
\text{M}\right)  }$ of the moduli space $\mathcal{M}\left(  \text{M}\right)
.$

\begin{theorem}
\label{W!!}We have on $\mathcal{T}_{\mathcal{M}\left(  \text{M}\right)  }$%
\begin{equation}
(dJ)_{i}^{\overline{l}}=e^{K}\sum_{j,k,l=1}^{N}C_{ijk}G^{k,\overline{l}}.
\label{Wit15}%
\end{equation}

\end{theorem}

\textbf{Proof:} Formula $\left(  \ref{Wit15}\right)  $ follows directly from
formula $\left(  \ref{Wit17}\right)  $, the globally defined isomorphism
$\iota_{\tau}:H^{1}\left(  \text{M,}T_{\text{M}}^{1,0}\right)  \approxeq
H^{1}\left(  \text{M,}\Omega_{\text{M}}^{2}\right)  ,$ given by $\phi
\rightarrow\phi\lrcorner(\eta_{\tau}),$ and the relation $g_{i,\overline{j}%
}=e^{-K}G_{i,\overline{j}},$ where $e^{-K}$ is given by $\left(
\ref{Wit12a}\right)  .$ Theorem \ref{W!!} is proved. $\blacksquare$

\subsection{BCOV\ Anomaly Equations}

\begin{definition}
\label{wcoor}\textit{The following data}, \textbf{a. }\textit{The moduli space
}$\mathcal{M}\left(  \text{M}\right)  $\textit{\ of CY threefolds, }\textbf{b.
}\textit{The CHSV} $\mathbb{S}p(2h^{2,1},\mathbb{R)}$ \textit{flat connection
on the tangent bundle constructed in }, \textbf{c. }\textit{The "prequantized
line bundle" }$\pi^{\ast}\mathcal{(}\omega_{\mathcal{Y}\left(  \text{M}%
\right)  \text{/}\mathcal{M}\left(  \text{M}\right)  }^{\ast})$\textit{\ on
the tangent bundle }$\mathcal{T}_{\mathcal{M}\left(  \text{M}\right)  }$,
\textbf{d. }The imaginary part of the Weil-Petersson metric $\omega$ and
\textbf{e. }\textit{the bundle of the Hilbert spaces }$\mathcal{H}^{Q}$
\textit{over the tangent bundle }$\mathcal{T}_{\mathcal{M}\left(
\text{M}\right)  }$ with a projective flat connection on it, \textit{will be
called a CY quantum system.}
\end{definition}

\begin{theorem}
\label{W!!!}The expression of Witten projective flat connection as defined in
Definition \ref{pc} in the flat coordinates for the CY\ quantum system defined
in Definition \ref{wcoor} coincides with the BCOV anomaly equations $\left(
\ref{Z0}\right)  $ in \cite{BCOV}.
\end{theorem}

\textbf{Proof: }In order to prove Theorem \ref{W!!!} we need to compute
explicitly the Witten projective connection constructed in \cite{ADW} on the
Hilbert vector bundle $\mathcal{H}^{Q}$ over the tangent bundle of the moduli
space $\mathcal{M}\left(  \text{M}\right)  $. The explicit formula $\left(
\ref{Wit4}\right)  $ and since BCOV anomaly equations were established in the
flat coordinate system $(\tau^{1},...,\tau^{N})$ imply that we need to compute
the expression of
\[
(dJ\omega^{-1})^{ij}D_{i}D_{j}%
\]
on $\mathcal{T}_{\mathcal{M}\left(  \text{M}\right)  }$ in the same coordinate
system. As we pointed out the flat coordinate system $(\tau^{1},..,\tau
^{N})\in\mathcal{K}$ introduced on the basis of Theorem \ref{tod1} is same
coordinate system used in \cite{BCOV}.

We already established in Theorem \ref{expl} the explicit expression of
$(dJ\omega_{1}^{-1})^{ij}D_{i}D_{j}$ on $R^{1}\pi_{\ast}\Omega_{\mathcal{Y}%
\left(  \text{M}\right)  \left/  \mathcal{M}\left(  \text{M}\right)  \right.
}^{2}$ in the coordinates $(\tau^{1},...,\tau^{N}).$ We will establish first
the local expression of $(dJ\omega^{-1})^{ij}D_{i}D_{j}$ on $\mathcal{T}%
_{\mathcal{M}\left(  \text{M}\right)  }$ in the coordinate system
$(t^{1},...,t^{N})$ defined by the identification $\mathcal{T}_{\mathcal{M}%
\left(  \text{M}\right)  }$ with $\omega_{\mathcal{Y}\left(  \text{M}\right)
\left/  \mathcal{M}\left(  \text{M}\right)  \right.  }^{\ast}\otimes R^{1}%
\pi_{\ast}\Omega_{\mathcal{Y}\left(  \text{M}\right)  \left/  \mathcal{M}%
\left(  \text{M}\right)  \right.  }^{2}$ that uses the parallel section
$\eta_{\tau}^{-1}$ of $\omega_{\mathcal{Y}\left(  \text{M}\right)  \left/
\mathcal{M}\left(  \text{M}\right)  \right.  }^{\ast}.$ Then we will compute
$(dJ\omega^{-1})^{ij}D_{i}D_{j}$ on $\mathcal{T}_{\mathcal{M}\left(
\text{M}\right)  }$ in the flat coordinate $(\tau^{1},...,\tau^{N})$ defined
by the identification of $\mathcal{T}_{\mathcal{M}\left(  \text{M}\right)  }$
with $\omega_{\mathcal{Y}\left(  \text{M}\right)  \left/  \mathcal{M}\left(
\text{M}\right)  \right.  }^{\ast}\otimes R^{1}\pi_{\ast}\Omega_{\mathcal{Y}%
\left(  \text{M}\right)  \left/  \mathcal{M}\left(  \text{M}\right)  \right.
}^{2}$ by the holomorphic tensor $\Omega_{\tau}^{-1}.$ Thus we will establish
BCOV anomaly equations $\left(  \ref{Z0}\right)  $ in \cite{BCOV} as the
projective Witten connection. Theorem \ref{W!!!} will be proved.

In order to compute $\left(  dJ\omega_{1,\tau}^{-1}D_{i}D_{j}\right)  $ we
need to established the relations between $g^{\overline{i},j}$ and
$G^{\overline{i},j}.$ $\left(  \ref{Wit12a}\right)  $ implies that these
relations given by
\begin{equation}
g^{\overline{i},j}=e^{K}G^{\overline{i},j}, \label{Wit17a}%
\end{equation}
where $e^{-K}=\left\Vert \eta_{\tau}\right\Vert ^{2}.$\textbf{ }

Let $(t^{1},...,t^{N})$ be the complex linear coordinates on $T_{\tau
,\mathcal{M}\left(  \text{M}\right)  }$ defined by the identification of
$\mathcal{T}_{\mathcal{M}\left(  \text{M}\right)  }$ with $\omega
_{\mathcal{Y}\left(  \text{M}\right)  \left/  \mathcal{M}\left(
\text{M}\right)  \right.  }^{\ast}\otimes R^{1}\pi_{\ast}\Omega_{\mathcal{Y}%
\left(  \text{M}\right)  \left/  \mathcal{M}\left(  \text{M}\right)  \right.
}^{2}$ that uses that parallel section $\eta$ of $\omega_{\mathcal{Y}\left(
\text{M}\right)  \left/  \mathcal{M}\left(  \text{M}\right)  \right.  }^{\ast
}.$ Remark \ref{2} and the relation $\Omega_{\tau}=\lambda(\tau)\eta_{\tau}$
imply that we have
\begin{equation}
(\tau^{1},..,\tau^{N})=\lambda(t^{1},..,t^{N}). \label{W0}%
\end{equation}
We know that the expression of $\left(  dJ\omega_{1}^{-1}\right)  ^{ij}%
D_{i}D_{j}$ on the fibre of $R^{1}\pi_{\ast}\Omega_{\mathcal{X}\text{/}%
\mathcal{K}}^{2}$ is given by $\left(  \ref{Wit17b}\right)  $ in the
coordinates $(\tau^{1},..,\tau^{N}).$ Thus combining $\left(  \ref{Wit17a}%
\right)  $ with $\left(  \ref{Wit17b}\right)  $ we get that on $T_{\tau
,\mathcal{M}\left(  \text{M}\right)  }$ in the coordinates $(t^{1},..,t^{N})$
we have:%
\begin{equation}
\left(  dJ\omega^{-1}\right)  ^{ij}D_{i}D_{j}=-\sqrt{-1}e^{2K}\sum_{a=1}%
^{N}\overline{C_{abc}}G^{\overline{b},i}G^{\overline{c},j}D_{i}D_{j}.
\label{Wit18}%
\end{equation}
From $\left(  \ref{Wit18}\right)  $ and $\left(  \ref{Wit4}\right)  $ we get
that the projective connection on the tangent bundle of $\mathcal{M}\left(
\text{M}\right)  $ is given by%
\begin{equation}
\frac{\partial}{\partial\tau^{a}}+\frac{\sqrt{-1}}{4}e^{2K}\sum\overline
{C_{abc}}G^{\overline{b},i}G^{\overline{c},j}\frac{D}{Dt^{i}}\frac{D}{Dt^{j}}
\label{Wit18a}%
\end{equation}
By using the symplectic identifications of the fibres of $\mathcal{T}%
_{\mathcal{M}\left(  \text{M}\right)  }$ by using the flat $\mathbb{S}%
p(2h^{2,1},\mathbb{R})$ and then by using the projective flat connection on
the Hilbert space fibration $\mathcal{H}^{Q}$ we obtain that the quantum state
represented by a parallel vector $\Psi(\tau,t)$ is independent of $t$. This
means that $\left(  \ref{Wit18a}\right)  $ implies that if $\Psi(\tau,t)$ is
independent of $t$ then $\Psi(\tau,t)$ satisfies the following equations:%

\begin{equation}
\left(  \frac{\partial}{\partial\tau^{a}}+\frac{\sqrt{-1}}{4}e^{2K}%
\sum\overline{C_{abc}}G^{\overline{b},i}G^{\overline{c},j}\frac{D}{Dt^{i}%
}\frac{D}{Dt^{j}}\right)  \Psi(\tau,t)=0 \label{Wit16}%
\end{equation}
and
\begin{equation}
\frac{\partial}{\partial\overline{t^{k}}}\Psi(\tau,t)=0. \label{Wit16a}%
\end{equation}
Based on $\left(  \ref{W0}\right)  $ we have%
\[
\frac{D}{Dt^{i}}=\lambda\frac{D}{D\tau^{i}}.
\]
Thus formulas $\left(  \ref{Wit16}\right)  $ and $\left(  \ref{Wit16a}\right)
$ can be written as follows:%

\begin{equation}
\left(  \frac{\partial}{\partial\tau^{a}}+\frac{\sqrt{-1}}{4}\lambda^{2}%
e^{2K}\sum\overline{C_{abc}}G^{\overline{b},i}G^{\overline{c},j}\frac{D}%
{D\tau^{i}}\frac{D}{D\tau^{j}}\right)  \Psi(\tau)=0 \label{Wit20}%
\end{equation}
and
\begin{equation}
\frac{\partial}{\partial\overline{\tau^{k}}}\Psi(\tau)=0. \label{Wit20a}%
\end{equation}
As we pointed out the flat coordinate system $(\tau^{1},..,\tau^{N}%
)\in\mathcal{K}$ introduced on the basis of Theorem \ref{tod1} is same
coordinates are use in \cite{BCOV}. Thus the BCOV anomaly equations $\left(
\ref{Z0}\right)  $ in \cite{BCOV} are the same as the equations $\left(
\ref{Wit20}\right)  $ and $\left(  \ref{Wit20a}\right)  .$

Direct computations show that if $F_{g}$ satisfy the equations $\left(
\ref{Z0}\right)  $ then
\[
\Psi(\tau)=\exp\left(
{\displaystyle\sum\limits_{g=1}^{\infty}}
\lambda^{2g-2}F_{g}\right)
\]
satisfy $\left(  \ref{Wit20}\right)  $. Theorem \ref{W!!!} is proved.
$\blacksquare$

\subsection{Comments}

\begin{itemize}
\item Bershadsky, Cecotti, Ooguri and Vafa used for the free energy
$Z(\lambda;t,\overline{t})$ the following expression:%
\begin{equation}
\mathfrak{F}(\lambda;t,\overline{t}):=%
{\displaystyle\sum\limits_{g=1}^{\infty}}
\lambda^{2g-2}F_{g}\text{ and }Z(\lambda;t,\overline{t})=\exp\left(
\mathfrak{F}(\lambda;t,\overline{t})\right)  . \label{Zb}%
\end{equation}
Compare this with the expression $\left(  \ref{Z}\right)  $ for $Z,$ i.e.
$Z=\exp\left(  \frac{1}{2}\mathfrak{F}(\lambda;t,\overline{t})\right)  $ used
by Witten in \cite{W93}.

\item It is proved in \cite{BCOV} by using physical arguments that $g$-genus
partition function $F_{g}$ satisfy the equation $\left(  \ref{Z0}\right)  $.
The holomorphic anomaly equation $\mathbf{(3.8)}$ derived in \cite{BCOV} is%
\begin{equation}
\left(  \overline{\frac{\partial}{\partial t^{i}}}-\overline{\frac{\partial
}{\partial t^{i}}}F_{1}\right)  \exp\mathfrak{F}(\lambda;t,\overline{t}%
)=\frac{\lambda^{2}}{2}\overline{C}_{ijk}e^{2K}G^{i,\overline{j}%
}G^{k,\overline{k}}\hat{D}_{j}\hat{D}_{k}\exp\mathfrak{F}(\lambda
;t,\overline{t}), \label{Zc}%
\end{equation}
where
\[
\hat{D}_{j}\mathfrak{F}(\lambda;t,\overline{t})=%
{\displaystyle\sum\limits_{g}}
\lambda^{2g-2}D_{j}F_{g}=
\]%
\[%
{\displaystyle\sum\limits_{g}}
\lambda^{2g-2}\left(  \partial_{j}-\left(  2g-2\right)  \partial_{j}K\right)
F_{g}=
\]%
\[
\left(  \partial_{j}-\partial_{j}K\lambda\partial_{\lambda}\right)
\mathfrak{F}(\lambda;t,\overline{t}).
\]
Thus the equation $\left(  \ref{Zc}\right)  $ is different by term involving
$F_{1}$ from the equation $\left(  \ref{Wit20}\right)  $.

\item In \cite{BCOV} the authors used the normalized holomorphic form, namely
they normalized $\Omega_{\tau}$ in such a way that $%
{\displaystyle\int\limits_{\gamma}}
\Omega_{\tau}=1,$ where $\gamma$ is the invariant vanishing cycle. This
normalized form is the same as the form defined in \ref{tod4}. They showed by
using string theory that the functions $F_{g}$ \textquotedblright
count\textquotedblright\ curves of genus g. It seems to me that it is a very
deep mathematical fact.
\end{itemize}

\section{The Extended Period Space of CY Manifolds.}

\subsection{Definition of the Extended Period Space and Basic Properties}

In this paragraph we will study the extended period space $\mathfrak{h}%
_{2,2h^{2,1}}\subset\mathbb{P(}H^{3}($M$,\mathbb{C)},$ which parametrizes all
possible filtrations of the type:%
\[
F^{0}=H^{3,0}\subset F^{1}=H^{3,0}+H^{2,1}+H^{1,2}\subset F^{2}=H^{3}%
(\text{M},\mathbf{C}),
\]
where $\dim F^{0}=1$ plus some extra properties which are motivated from the
above filtration and Variations of Hodge Structures on K3 surfaces.

We will use the following notation for the cup product for $\Omega_{1}$ and
$\Omega_{2}\in$ $H^{3}($M,$\mathbb{C)}$, i.e. $\left\langle \Omega_{1}%
,\Omega_{2}\right\rangle =%
{\displaystyle\int\limits_{\text{M}}}
\Omega_{1}\wedge\Omega_{2}.$

\begin{definition}
\label{Hodge}$\mathfrak{h}_{2,2h^{2,1}}$ by definition is the set of lines in
$H^{3}($M,$\mathbb{C)}$ spanned by the cohomology class $[\Omega]$ of the
holomorphic three form in $H^{3}($M,$\mathbb{C)}$ such that%
\begin{equation}
-\sqrt{-1}%
{\displaystyle\int\limits_{\text{M}}}
\Omega\wedge\overline{\Omega}=-\sqrt{-1}\left\langle \Omega,\overline{\Omega
}\right\rangle >0. \label{Ho1}%
\end{equation}

\end{definition}

\begin{theorem}
\label{char}There is a one to one map between the points $\tau$ of
$\mathfrak{h}_{2,2h^{1,2}}$ and the two dimensional oriented planes $E_{\tau}$
in $H^{3}$(M,$\mathbb{R)}$, where $E_{\tau}$ is generated by $\gamma_{1}$ and
$\mu_{1}$ such that $\left\langle \gamma_{1},\mu_{1}\right\rangle =1.$
\end{theorem}

\textbf{Proof: }Let $[\Omega_{\tau}]\in H^{3}$(M,$\mathbb{C}$) be a non-zero
vector corresponding to a point $\tau\in\mathfrak{h}_{2,2h^{1,2}}.$ Since the
class of cohomology $[\Omega_{\tau}]$ satisfies $\left(  \ref{Ho1}\right)  $
we may chose $\Omega_{\tau}$ such that
\begin{equation}
-\sqrt{-1}\left\langle \Omega_{\tau},\overline{\Omega_{\tau}}\right\rangle =2.
\label{Ho2}%
\end{equation}
Then $\left(  \ref{Ho1}\right)  $ and $\left(  \ref{Ho2}\right)  $ imply$%
{\displaystyle\int\limits_{\text{M}}}
\operatorname{Im}\Omega_{\tau}\wedge\operatorname{Re}\Omega_{\tau
}=\left\langle \operatorname{Re}\Omega_{\tau},\operatorname{Im}\Omega_{\tau
}\right\rangle =1.$ We will define $E_{\tau}$ to be the the two dimensional
oriented subspace in $H^{3}$(M,$\mathbb{R}$) spanned by $\operatorname{Re}%
\Omega_{\tau}$ and $\operatorname{Im}\Omega_{\tau}.$ So to each point $\tau
\in\mathfrak{h}_{2,2h^{1,2}}$ we have assigned an oriented two plane $E_{\tau
}$ in $H^{3}$(M,$\mathbb{R)}.$

Suppose that $E$ is a two dimensional oriented plane in $H^{3}$(M,$\mathbb{R}%
$) spanned by $\gamma$ and $\mu\in H^{3}$(M,$\mathbb{R}$) be such that
$\left\langle \gamma,\mu\right\rangle =1.$ Let $\Omega_{E}:\mu+\sqrt{-1}%
\gamma.$ Then we have%
\[
-\sqrt{-1}%
{\displaystyle\int\limits_{\text{M}}}
\Omega_{E}\wedge\overline{\Omega_{E}}=-\sqrt{-1}\left\langle \Omega
_{E},\overline{\Omega_{E}}\right\rangle =2\left\langle \gamma,\mu\right\rangle
=2.
\]
So to the plane $E$ we assign the line in $H^{3}$(M,$\mathbb{C}$) spanned by
$\Omega_{E}.$ This proves Theorem \ref{char}. $\blacksquare$

\begin{corollary}
$\mathfrak{h}_{2,2h^{1,2}}$ is an open set in $\mathbf{Grass}(2,2h^{1,2}+2).$
So it has a complex dimension $2h^{1,2}.$
\end{corollary}

\begin{remark}
\label{weight2}It is easy to see that each point $\tau\in\mathfrak{h}%
_{2,2h^{1,2}}$ defines a natural filtration of length two in $H^{3}%
$(M,$\mathbb{C}$). Indeed,\ let $H_{\tau}^{3,0}$ be the subspace in $H^{3}%
$(M,$\mathbb{C}$) spanned by a non-zero element $\Omega_{\tau}\in F^{0}.$ Let
$\gamma_{0}=\operatorname{Re}\Omega_{\tau}$ and $\operatorname{Im}\Omega
_{\tau}=\mu_{0}.$ From Theorem \ref{char} we know that $\left\langle
\gamma_{0},\mu_{0}\right\rangle $ is a positive number. Let%
\[
\{\gamma_{0},\mu_{0},..,\gamma_{h^{1,2}},\mu_{h^{1,2}}\}
\]
\textit{be a symplectic basis of }$H^{3}$(M,$\mathbb{R}$)\textit{\ such that
}$\left\langle \gamma_{i},\mu_{j}\right\rangle =-\delta_{ij}.$ Let $\Omega
_{i}:=\mu_{i}+\sqrt{-1}\gamma_{i}.$ \textit{We will define }$H_{\tau}^{2,1}$
\textit{to be the subspace in }$H^{3}($M,$\mathbb{C)}$\textit{\ spanned by the
vectors \ } $\Omega_{i}$ \textit{for }$i=1,..,h^{1,2}.$ \textit{Then we define
}$H_{\tau}^{1,2}:=\overline{H_{\tau}^{2,1}}.$ \textit{It is easy to see that}%
\[
\mathit{\ }\left(  H_{\tau}^{2,1}+H_{\tau}^{1,2}\right)  ^{\perp}=H_{\tau
}^{3,0}+H_{\tau}^{0,3},
\]
\textit{where }$\mathit{\ }H_{\tau}^{0,3}:=\overline{H_{\tau}^{3,0}}$ and the
orthogonality is with respect to%
\[
\left\langle \omega_{1},\overline{\omega_{2}}\right\rangle =%
{\displaystyle\int\limits_{\text{M}}}
\omega_{1}\wedge\overline{\omega_{2}}.
\]
\textit{The natural filtration in }$H^{3}$(M,$\mathbb{R}$) \textit{that
corresponds} to $\tau\in\mathfrak{h}_{2,2h^{1,2}}$ \textit{is defined as
follows:}%
\begin{equation}
F_{\tau}^{0}=H_{\tau}^{3,0}\subset F_{\tau}^{1}=H_{\tau}^{3,0}+H_{\tau}%
^{2,1}+H_{\tau}^{1,2}\subset H^{3}(M,\mathbb{C)}. \label{Ho3a}%
\end{equation}

\end{remark}

Next we will introduce a metric G on $H^{3}($M,$\mathbf{C}).$ We will use the
metric G to show that the filtration defined by \ref{Ho3a} is a Hodge
filtration of weight two.

\begin{definition}
\label{metric}Let M be a fixed CY manifold and let%
\[
\Omega=\Omega^{3,0}+\Omega^{2,1}+\Omega^{1,2}+\Omega^{0,3}\in H^{3}%
(\text{M},\mathbf{C)}%
\]
\textit{be the Hodge decomposition of }\ $\Omega$, \textit{then}
$G(\Omega,\overline{\Omega})$ \textit{is defined as follows}:%
\[
G(\Omega,\overline{\Omega}):=-\sqrt{-1}\left(
{\displaystyle\int\limits_{\text{M}}}
\Omega^{3,0}\wedge\overline{\Omega^{3,0}}\ +%
{\displaystyle\int\limits_{\text{M}}}
\Omega^{2,1}\wedge\overline{\Omega^{2,1}}\right)  +
\]%
\begin{equation}
-\sqrt{-1}\left(
{\displaystyle\int\limits_{\text{M}}}
\overline{\Omega^{1,2}}\ \wedge\Omega^{1,2}+%
{\displaystyle\int\limits_{\text{M}}}
\overline{\Omega^{0,3}}\wedge\Omega^{0,3}\right)  . \label{Ho4}%
\end{equation}

\end{definition}

From the definition of the metric, it follows that it has a signature
$(2,2h^{2,1})$ on H$^{3}$(M,$\mathbb{R}$\textbf{)}. We will denote the
quadratic form of this metric by Q.

\begin{lemma}
\label{metrica}The metric defined by $\left(  \ref{Ho4}\right)  $ does not
depend on the choice of the complex structure on M.
\end{lemma}

\textbf{Proof: }Let M$_{0}$ and M$_{\tau}$ be two different complex structures
on M. Let $\Omega_{0}$ and $\Omega_{\tau}$ be two non zero holomorphic three
forms on M$_{0}$ and M$_{\tau}$ respectively. Let $\{\Omega_{0,i}\}$ and
$\{\Omega_{\tau,i}\}$ be two bases of $H^{2,1}($M$_{0})$ and $H^{2,1}%
($M$_{\tau})$ respectively, where $i=1,..,h^{2,1}$ such that $\left\langle
\Omega_{0,i},\overline{\Omega_{0,j}}\right\rangle =\left\langle \Omega
_{\tau,i},\overline{\Omega_{\tau,j}}\right\rangle =\delta_{ij}.$ Then one see
immediately that%
\[
\left\{  \operatorname{Re}\Omega_{0},\operatorname{Im}\Omega_{0}%
,...,\operatorname{Re}\Omega_{0,i},\operatorname{Im}\Omega_{0,i},..\right\}
\ and\text{ }\left\{  \operatorname{Re}\Omega_{\tau},\operatorname{Im}%
\Omega_{\tau},..,\operatorname{Re}\Omega_{\tau,i},\operatorname{Im}%
\Omega_{\tau,i},..\right\}
\]
are two different symplectic bases of $H^{3}($M,$\mathbb{R}$). So there exists
an element g$\in\mathbb{S}p(2h^{2,1}+2)$ such that
\begin{equation}
\text{g}(\Omega_{0})=\Omega_{\tau}\text{ and g}(\Omega_{0,i})=\Omega_{\tau,i}.
\label{Ho5}%
\end{equation}
So $\left(  \ref{Ho5}\right)  $ implies that%

\begin{equation}
\text{g}(H^{2,1}(\text{M}_{0}))=H^{2,1}(\text{M}_{\tau}) \label{Ho5a}%
\end{equation}
and
\begin{equation}
\text{g}(H^{3,0}(\text{M}_{0}))=H^{3,0}(\text{M}_{\tau}) \label{Ho5b}%
\end{equation}
Lemma \ref{metrica} follows directly from $\left(  \ref{Ho5a}\right)  $,
$\left(  \ref{Ho5b}\right)  $ and the definition of the metric G.
$\blacksquare$

\begin{theorem}
\label{Hodge1}Let $\mathfrak{h}_{2,2h^{1,2}}$ be the space that parametrizes
all filtrations in $H^{3}($M$,\mathbb{C}\mathbf{)}$\textit{\ defined in
Remark} \ref{weight2}. Then these filtration are Hodge filtrations of weight
two and their moduli space \textit{is isomorphic to the symmetric space}%
\[
\mathit{\ }\mathfrak{h}_{2,2h^{1,2}}:=\mathbb{SO}_{0}(2,2h^{1,2}%
)/\mathbb{SO}(2)\times\mathbb{SO}(2h^{1,2}).
\]

\end{theorem}

\textbf{Proof: }Since the signature of the metric G is $(2,2h^{2,1})$ the
proof of Theorem \ref{Hodge1} is standard and follows directly from definition
of the variations of Hodge structures of weight two. See for example \cite{G}.
$\blacksquare$

We will use the fact that the space $\mathfrak{h}_{2,2h^{1,2}}=\mathbb{SO}%
_{0}(2,2h^{1,2})/\mathbb{SO}(2)\times\mathbb{SO}(2h^{1,2})$ is as an open set
in the Grassmannian \textbf{$Grass$}$(2,b_{3})$ identified with all two
dimensional oriented subspaces in $H^{3}\left(  \text{M},\mathbb{R}\right)  $
such that the restriction of the quadratic form Q on them is positive, i.e.
\[
\mathfrak{h}_{2,h^{1,2}}:=\{E\subset H^{3}\left(  \text{M},\mathbb{R}\right)
\dim E=2,Q\left\vert _{E}\right.  >0+orientation\}.
\]

\begin{definition}
\label{charb}\textit{We will define a canonical map from }$\mathfrak{h}%
_{2,2h^{1,2}}\subset Gr(2,2h^{1.2}+2)$ \textit{to} $\mathbb{P}\mathbf{(}%
H^{3}\left(  \text{M},\mathbb{R}\right)  \mathbf{\otimes\mathbb{C})}$
\textit{as follows}; \textit{let} $E_{\tau}\in\mathfrak{h}_{2,2h^{1,2}},$ i.e.
$E_{\tau}$ is an oriented two dimensional subspace in $H^{3}\left(
\text{M},\mathbb{R}\right)  $ on which the restriction of Q is
positive.\textit{\ Let} $e_{1}$ and $e_{2}$ \textit{be an orthonormal basis of
}$E_{\tau}$. \textit{Let }$\Omega_{\tau}:=\mathit{\ }e_{1}+\sqrt{-1}e_{2}.$
Then $\Omega_{\tau}$ \textit{defines a point }$\tau\in\mathbb{P}%
\mathbf{(}H^{3}($M$,\mathbb{Z)\otimes C)})$ that corresponds to the line in
$H^{3}\left(  \text{M},\mathbb{R}\right)  \mathbb{\otimes C)}$ spanned by
$\Omega_{\tau}.$ It is a standard fact that \textit{the points }$\tau$ in
$\mathbb{P}\mathbf{(}H^{3}\left(  \text{M},\mathbb{R}\right)  \mathbb{\otimes
C)})$ is such\textit{\ that} Q($\tau,\tau$)$=0$ \textit{and }Q($\tau
,\overline{\tau})>0$ are in one to one corresponds with the points in
$\mathfrak{h}_{2,h^{1,2}}.$ See \cite{to1}. This follows from the arguments
used in the proof of Theorem \ref{char} or see \cite{to1}.
\end{definition}

It is a well known fact that $\mathfrak{h}_{2,h^{1,2}}$ is isomorphic to one
of the irreducible component of the open set of the quadric in $\mathbb{P}%
\mathbf{(}H^{3}\left(  \text{M},\mathbb{R}\right)  \mathbf{\otimes\mathbb{C}%
)}$ defined as follows:%
\[
\mathfrak{h}_{2,2h^{1,2}}\approx\{\tau\in\mathbb{P}(H^{3}(\text{M}%
,\mathbb{Z}\mathbf{)\otimes\mathbb{C})}\left\vert Q(\tau,\tau)=0\text{
}and\text{ }Q(\tau,\overline{\tau})>0\right.  \}.
\]
(See \cite{G}.)

We will consider the family $\mathcal{X\times}\overline{\mathcal{X}%
}\rightarrow\mathcal{K\times}\overline{\mathcal{K}},$ where the
family\textit{\ }$\ \overline{\mathcal{X}}\rightarrow\overline{\mathcal{K}}%
$\textit{\ \ }is the family that corresponds to the conjugate complex
structures, i.e. the point $(\tau_{1},\overline{\tau_{2}})\in\mathcal{K\times
K}$ $\tau_{1}$ corresponds to the complex structure M$_{\tau_{1}}$ and
$\overline{\tau_{2}}$ corresponds to the $\overline{\text{M}_{\tau_{2}}},$
where $\overline{\text{M}_{\tau_{2}}}$ is the conjugate complex structure on
M$_{\tau_{2}}.$

We will define the period map%

\begin{equation}
p:\mathcal{K\times}\overline{\mathcal{K}}\mathcal{\rightarrow}\mathbb{P}%
\mathbf{(}H^{3}(\text{M},\mathbb{Z}\mathbf{)\otimes\mathbb{C})} \label{per}%
\end{equation}
as follows; to each point $(\tau,\upsilon)\in\mathcal{K\times}\overline
{\mathcal{K}}$ we will assign\textit{\ }the complex line $p(\tau,\upsilon)$ in
$H^{3}\left(  \text{M},\mathbb{R}\right)  \mathbf{\otimes}\mathbb{C}$ defined
by the oriented two plane $E_{\tau,\upsilon}\subset H^{3}($M$,\mathbb{Z}%
\mathbf{)\otimes}\mathbb{R}$ spanned by $\operatorname{Re}(\Omega_{\tau
}+\overline{\Omega_{\upsilon}})$ \textit{\ }and\textit{\ \ }$\operatorname{Im}%
(\Omega_{\tau}-\Omega_{\upsilon})$, where $\Omega_{\tau}$ and $\Omega_{\nu}$
are defined as in Theorem \ref{tod2}. $\mathcal{K}$ is the Kuranishi space
defined in Definition \ref{coor}. We will show that the analogue of local
Torelli Theorem holds, i.e. we will show that the period map $p$ is a local
embedding $\ p:\mathcal{K\times}\overline{\mathcal{K}}\mathcal{\subset
}\mathfrak{h}_{2,2h^{1,2}}.$

\begin{remark}
\label{emb}We will define an embedding of the Kuranishi family $\mathcal{K}$
into $\mathcal{K\times}\overline{\mathcal{K}}$ as follows; to each $\tau
\in\mathcal{K}$ we will associate the complex structure ($I_{\tau},-I_{\tau})$
on M$\times$M.
\end{remark}

\begin{theorem}
\label{torelli}The period map $p$ defined by $\left(  \ref{per}\right)  $ is a
local isomorphism. Moreover the image $p(\mathcal{K\times}\overline
{\mathcal{K}}\mathcal{)}$ is contained in $\mathfrak{h}_{2,2h^{1,2}}%
\subset\mathbb{P}\left(  H^{3}(\text{M,}\mathbb{Z)}\otimes\mathbb{C}\right)  $
for small enough $\varepsilon.$
\end{theorem}

\textbf{Proof: }The fact that the period map $p$ is a local isomorphism
follows directly from the local Torelli theorem for CY manifolds proved in
\cite{gr}. Let $(\tau,\overline{\tau})\in\bigtriangleup\subset\mathcal{K\times
}\overline{\mathcal{K}},$ then clearly the point $p(\tau,\overline{\tau}%
)\in\mathfrak{h}_{2,2h^{1,2}}\subset\mathbb{P(}H^{3}\left(  \text{M}%
,\mathbb{R}\right)  \mathbb{\otimes C)},$ i.e. Q%
$\vert$%
$_{E_{\tau,\overline{\tau}}}>0.$ This follows directly from the definition of
Q and the fact that $E_{\tau,\overline{\tau}}$ is the subspace in $H^{3}%
($M$_{\tau},\mathbb{R})$ spanned by $\operatorname{Re}\Omega_{\tau}$ and
$\operatorname{Im}\Omega_{\tau}.$

Let $(\tau_{1},\overline{\tau_{2}})\in\mathcal{K\times}\overline{\mathcal{K}}$
be a point "close" to the diagonal $\Delta$ in $\mathcal{K\times}%
\overline{\mathcal{K}}$ then Q$\left\vert _{E_{\tau,\upsilon}}\right.  >0.$
Indeed this follows from the fact that the condition Q$\left\vert
_{E_{\tau,\upsilon}}\right.  >0$ is an open one on $Gr(2,2h^{1,2}+2)$. Then
the two dimensional oriented space $E_{(\tau_{1},\overline{\tau_{2}})}\subset
H^{3}\left(  \text{M},\mathbb{R}\right)  $ spanned by $\left\{
\operatorname{Re}\Omega_{\tau_{1}}+\operatorname{Re}\Omega_{\tau_{2}%
},\operatorname{Im}\Omega_{\tau_{1}}+\operatorname{Im}\Omega_{\tau_{2}%
}\right\}  ,$ where $\Omega_{\tau_{1}}$ and $\Omega_{\tau_{2}}$ are defined by
formula $\left(  \ref{T1}\right)  $ in Theorem \ref{tod2} will be such that
Q$\left\vert _{E_{\tau,\upsilon}}\right.  >0$. From here we deduce that
$p(\tau_{1},\overline{\tau_{2}})\in\mathfrak{h}_{2,2h^{1,2}}.$ Theorem
\ref{torelli} is proved. $\blacksquare$

We defined the Kuranishi space $\mathcal{K}$ to be a open polydisk $|\tau
^{i}|<\varepsilon$ for $i=1,...,N$ in $H^{1}\left(  \text{M}\right)  .$
$\,$Since $p$ is a local isomorphism we may assume that $\mathcal{K\times
}\overline{\mathcal{K}}$ is contained in $\mathfrak{h}_{2,2h^{1,2}}$ for small
enouph $\varepsilon.$

\subsection{\textbf{Construction of a }$\mathbb{Z}$\textbf{\ Structure on the
Tangent Space of }$\mathcal{M}\left(  \text{M}\right)  $}

\begin{definition}
\label{zstr}To define a $\mathbb{Z}$ structure on a complex vector space $V$
means the construction of a free abelian group $A\subset V$ such that the rank
of $A$ is equal to the dimension of V, i.e. $A\otimes\mathbb{C}=V.$
\end{definition}

\begin{definition}
We define $\mathfrak{h}_{2,2h^{1,2}}(\mathbb{Q)}$ as follows; A point $\tau
\in\mathfrak{h}_{2,2h^{1,2}}(\mathbb{Q)}$ if the two dimensional oriented
subspace $E_{\tau}=H_{\tau}^{3,0}+H_{\tau}^{0,3}$ that corresponds to $\tau$
constructed in Theorem \ref{char} is such that $E_{\tau}\subset H^{3}%
($M,$\mathbb{Z)}\otimes\mathbb{Q}.$
\end{definition}

\begin{theorem}
\label{dense}$\mathfrak{h}_{2,2h^{1,2}}(\mathbb{Q)}$ is an everywhere dense
subset in $\mathfrak{h}_{2,2h^{1,2}}.$
\end{theorem}

\textbf{Proof: }Our claim follows directly from two facts. The first one is
that the set of the points in $Gr(2,2+2h^{2,1})$ that corresponds to two
dimensional subspaces in $H^{3}($M,$\mathbb{Z)}\otimes\mathbb{Q}$ form an
everywhere dense subset in $Gr(2,2+2h^{2,1})$ and the second one is that
$\mathfrak{h}_{2,2h^{2,1}}$ is an open set in $Gr(2,2+2h^{2,1}).$ Theorem
\ref{dense} is proved. $\blacksquare$

\begin{theorem}
\label{zstr1}For each $\tau\in\mathfrak{h}_{2,2h^{2,1}}(\mathbb{Q)}$ a natural
$\mathbb{Z}$ structure is defined on the tangent space T$_{\tau,\mathfrak{h}%
_{2,2h^{2,1}}}$ at the point $\tau\in\mathfrak{h}_{2,2h^{2,1}}.$ This means
that there exists a subspace $\mathbb{Z}^{2h^{2,1}}\subset H^{3}%
($M$,\mathbb{Z)}$ such that $H_{\tau}^{2,1}+H_{\tau}^{1,2}\approxeq
T_{\tau,\mathfrak{h}_{2,2h^{2,1}}}=\mathbb{Z}^{2h^{2,1}}\mathbb{\otimes R}.$
\end{theorem}

\textbf{Proof: }From the theory of Grassmannians we know that $T_{\tau
,\mathfrak{h}_{2,2h^{2,1}}}$ can be identified with $H_{\tau}^{2,1}+H_{\tau
}^{1,2}.$ Our corollary follows directly from the construction of $H_{\tau
}^{2,1}+H_{\tau}^{1,2}$ described in Remark \ref{weight2}. $\ $Indeed the
point $\tau\in\mathfrak{h}_{2,2h^{2,1}}(\mathbb{Q)}$ defines two vectors
$\gamma_{0}$ and $\mu_{0}\in H^{3}($M$,\mathbb{Z})/Tor$ that form a basis of
$H^{3,0}\left(  \text{M}\right)  \oplus H^{0,3}\left(  \text{M}\right)  $ such
that $\left\langle \mu_{0},\gamma_{0}\right\rangle \in\mathbb{Z}$ and
$\left\langle \mu_{0},\gamma_{0}\right\rangle >0.$ We choose the vectors
\[
\{\gamma_{0},\mu_{0},\gamma_{1},\mu_{1},..,\gamma_{h^{2,1}},\mu_{h^{2,1}}\}
\]
to be in $H^{3}($M,$\mathbb{Z}$)/Tor and we require that $\left\langle \mu
_{i},\gamma_{j}\right\rangle =\delta_{ij}.$ Then, from the way we defined
$H_{\tau}^{2,1}$ and $H_{\tau}^{1,2},$ it follows that%
\[
H_{\tau}^{2,1}+H_{\tau}^{1,2}=(\mathbb{Z\gamma}_{1}\oplus\mathbb{Z\mu}%
_{1}\oplus...\mathbb{Z\gamma}_{h^{2,1}}\oplus\mathbb{Z\mu}_{h^{2,1}}%
)\otimes\mathbb{R}.
\]
Theorem \ref{zstr1} is proved$.$ $\blacksquare$

We will consider the embedding of $\mathcal{K}$ in $\mathcal{K\times}%
\overline{\mathcal{K}}$ defined in Remark \ref{emb}. Next we choose a point
$\kappa\in\mathcal{K\times}\overline{\mathcal{K}}\mathcal{\subset}%
\mathfrak{h}_{2,2h^{2,1}}\subset\mathbb{P}(H^{3}($M$,\mathbb{Z)\otimes C})$
such that $\kappa\in(\mathcal{K\times}\overline{\mathcal{K}}\mathcal{)}%
\cap\mathfrak{h}_{2,2h^{2,1}}(\mathbb{Q)}.$ We know that $\kappa$ corresponds
to a two dimensional space $E_{\kappa}\subset H^{3}($M,$\mathbb{Q)}$, with the
additional condition, that there exists vectors $\gamma_{0}$ and $\mu_{0}\in
H^{3}($M,$\mathbb{Z}$)/Tor \ that span $E_{\kappa}$ and \ $\left\langle
\mu_{0},\gamma_{0}\right\rangle =1$. The existence of such points follows from
the fact that the set of all two dimensional space, $E_{\kappa}\subset H^{3}%
($M$,\mathbb{Q)}$ such that there exists vectors $\gamma_{0}$ and $\mu_{0}\in
H^{3}($M$,\mathbb{Z})/Tor$ and\ $\left\langle \mu_{0},\gamma_{0}\right\rangle
=1$ is an everywhere dense subset in $\mathfrak{h}_{2,2h^{2,1}}$. Let%
\[
\left\{  \gamma_{0},\mu_{0},\gamma_{1},\mu_{1},..,\gamma_{h^{2,1}}%
,\mu_{h^{2,1}}\right\}  \in H^{3}(\text{M},\mathbb{Z})/Tor
\]
be such that $\gamma_{0}$ and $\mu_{0}$ span $E_{\kappa}$ and $\left\langle
\mu_{i},\gamma_{j}\right\rangle =\delta_{ij}.$ It follows from the
construction in Remark \ref{weight2} that the vectors%
\[
\gamma_{1},\mu_{1},..,\gamma_{h^{2,1}},\mu_{h^{2,1}}\in H^{3}(\text{M}%
,\mathbb{Z})/Tor
\]
span $H_{\kappa}^{2,1}+H_{\kappa}^{1,2},$ i.e. they span the tangent space
$T_{\kappa,\mathfrak{h}_{2,2h^{2,1}}}=H_{\kappa}^{2,1}+H_{\kappa}^{1,2}.$ We
know from Corollary \ref{CHSV1} that there exists an $\mathbb{S}%
p(2h^{2,1},\mathbb{R)}$ flat connection on $\ \mathcal{K}$ and so we define a
flat connection on the product $\mathcal{K\times}\overline{\mathcal{K}}$ as
the direct sum of the two connections. Using this $\mathbb{S}p(4h^{2,1}%
,\mathbb{R)}$ flat connection we can perform a parallel transport of the
vectors $\gamma_{1,}\mu_{1},..,\gamma_{h^{2,1}},\mu_{h^{2,1}}\in T_{\kappa}$
to a basis%
\[
\gamma_{1,\tau},\mu_{1,\tau},..,\gamma_{h^{2,1},\tau},\mu_{h^{2,1},\tau}%
\]
in the tangent space$\ H_{\tau}^{2,1}+H_{\tau}^{1,2}\approxeq T_{\tau
,\mathfrak{h}_{2,2h^{2,1}}}$ to each point $\tau\in\mathcal{K\subset K\times
}\overline{\mathcal{K}}.$ Thus we can conclude that $\left\langle \mu_{i,\tau
},\gamma_{j,\tau}\right\rangle =\delta_{ij}$ and\ the free abelian group%
\[
A_{\tau}:=\mathbb{Z\gamma}_{1,\tau}\oplus\mathbb{Z\mu}_{1,\tau}\oplus
...\mathbb{Z\gamma}_{h^{2,1},\tau}\oplus\mathbb{Z\mu}_{h^{2,1},\tau}%
\]
in $H_{\tau}^{2,1}+H_{\tau}^{1,2}$ is such that%
\[
H_{\tau}^{2,1}+H_{\tau}^{1,2}=(\mathbb{Z\gamma}_{1,\tau}\oplus\mathbb{Z\mu
}_{1,\tau}\oplus...\mathbb{Z\gamma}_{h^{2,1},\tau}\oplus\mathbb{Z\mu}%
_{h^{2,1},\tau})\otimes\mathbb{R}.
\]
So we defined for each $\tau\in\mathcal{K}$ an abelian subgroup $A_{\tau
}\subset$ T$_{\tau,\mathcal{K}}=\left(  H_{\tau}^{2,1}+H_{\tau}^{1,2}\right)
$ such that
\[
A_{\tau}\otimes\mathbb{C=}T_{\tau,\mathcal{K}}=(H_{\tau}^{2,1}+H_{\tau}%
^{1,2})
\]
and $\left\langle \gamma,\mu\right\rangle \in\mathbb{Z}$ for $\gamma$ and
$\mu\in A_{\tau}$.

\begin{definition}
\label{abel} The image of the projection of \ the abelian subgroup $A_{\tau}$
of T$_{\tau,\mathcal{K}}=(H_{\tau}^{2,1}+H_{\tau}^{1,2})$ to $H_{\tau}^{2,1}$
will be denoted by $\Lambda_{\tau}$ for each $\tau\in\mathcal{M}\left(
\text{M}\right)  $.
\end{definition}

\begin{theorem}
\label{abel1}There exists a holomorphic map $\phi$ from the moduli space
$\mathcal{M}\left(  \text{M}\right)  $ of CY manifolds to the moduli space of
principally polarized abelian varieties $\mathbb{S}p(2h^{2,1},\mathbb{Z)}%
\backslash\mathfrak{Z}_{h^{2,1}}$, where $\mathfrak{Z}_{h^{2,1}}%
:=\mathbb{S}p(2h^{2,1},\mathbb{R)}/U(h^{2,1}).$
\end{theorem}

\textbf{Proof: }From \ref{abel} we know that there exists a lattice
$\Lambda_{\tau}\subset$ T$_{\tau,\mathcal{K}}$ such that the restriction of
the imaginary part of the Weil-Petersson metric
\[
\operatorname{Im}(g)(u,v)=\left\langle u,v\right\rangle
\]
on $\Lambda_{\tau}$ is such that $\left\langle u,v\right\rangle \in\mathbb{Z}
$ and $|\det\left\langle \gamma_{i},\gamma_{j}\right\rangle |=1$ for any
symplectic basis of $\Lambda_{\tau}.$ Thus over $\mathcal{K}$ we can construct
a family of principally polarized abelian varieties
\begin{equation}
\mathcal{A}_{K}\rightarrow\mathcal{K}. \label{AB}%
\end{equation}
In fact we constructed a family of principally polarized abelian varieties%

\begin{equation}
\mathcal{A\rightarrow M}\left(  \text{M}\right)  \label{ABG}%
\end{equation}
over the moduli space $\mathcal{M}\left(  \text{M}\right)  $ since CHSV
connection is a flat connection globally defined over $\mathcal{M}\left(
\text{M}\right)  $. This means that we defined the holomorphic map $\phi$
between the quasi-projective varieties%
\[
\phi:\mathcal{M}\left(  \text{M}\right)  \rightarrow\mathbb{S}p(2h^{2,1}%
,\mathbb{Z)}\backslash\mathfrak{Z}_{h^{2,1}}.
\]
The existence of $\phi$ follows from the fact that there exists a versal
family of principally polarized abelian varieties%
\[
\mathfrak{A\rightarrow}\mathbb{S}p(2h^{2,1},\mathbb{Z)}\backslash
\mathfrak{Z}_{h^{2,1}}.
\]
Theorem \ref{abel1} is proved. $\blacksquare$

Notice that the family of principally polarized varieties
$\mathcal{A\rightarrow M}\left(  \text{M}\right)  $ is constructed by using
the vector bundle $R^{1}\pi_{\ast}\Omega_{\mathcal{Y}\left(  \text{M}\right)
\left/  \mathcal{M}\left(  \text{M}\right)  \right.  }^{2}.$ Using the
identification between $\mathcal{T}_{\mathcal{M}\left(  \text{M}\right)  }$
and $R^{1}\pi_{\ast}\Omega_{\mathcal{Y}\left(  \text{M}\right)  \left/
\mathcal{M}\left(  \text{M}\right)  \right.  }^{2}$ given by $\phi
\rightarrow\phi\lrcorner\eta_{\tau},$ we define a family \ of principally
polarized varieties isomorphic to the family $\mathcal{A\rightarrow M}\left(
\text{M}\right)  .$

\subsection{Holomorphic Symplectic Structure on the Extended Period Domain}

\begin{definition}
Let us fix a symplectic basis
\[
\{\gamma_{0},\gamma_{1},..,\gamma_{h^{1,2}};\upsilon_{0},...,\upsilon
_{h^{1,2}}\}
\]
\textit{in} $H^{3}($M,$\mathbb{Z}$)/Tor, \textit{i.e. }$\left\langle
\gamma_{i},\upsilon_{j}\right\rangle =\delta_{ij}.$ \textit{This basis defines
a coordinate system in }$\mathbb{P(}H^{3}($M,$\mathbb{Z)\otimes C)}%
$\textit{\ which we will denote by \ }$(z_{0}:..:z_{2h^{1,2}}).\mathit{\ }%
$\textit{On the open set: }%
\[
U_{0}:=\{(z_{0}:..:z_{2h^{1,2}+1})\left\vert z_{0}\neq0\right.  \},
\]
\textit{we define a holomorphic one forms:}%
\begin{equation}
\alpha_{0}:=dt_{h^{1,2}+1}+t_{1}dt_{h^{1,2}+2}+..+t_{h^{1,2}}dt_{2h^{1,2}+1},
\label{s0}%
\end{equation}
\textit{where }$t_{i}=\frac{z_{i}}{z_{0}},$ for $i=1,..,h^{1,2}.$ \textit{Let
us restrict }$\alpha_{0}$ on $\mathfrak{h}_{2,2h^{1,2}}\cap U_{0}$ \textit{and
denote this restriction by } $\alpha_{0}.$ \textit{In the same way we can
define the forms} $\alpha_{i}$ \textit{on the open set }$U_{i}:=\{(z_{0}%
:..:z_{2h^{1,2}})|$ $z_{i}\neq0\}.$
\end{definition}

\begin{theorem}
\label{sympl}There exists a closed holomorphic non degenerate two form $\psi$
on $\mathfrak{h}_{2,2h^{1,2}}$ such that
\begin{equation}
\psi\left\vert _{U_{i}\cap\mathfrak{h}_{2,2h^{1,2}}}\right.  =d\alpha
_{i}\left\vert _{U_{i}\cap\mathfrak{h}_{2,2h^{1,2}}}\right.  . \label{s1}%
\end{equation}

\end{theorem}

\textbf{Proof: }The proof of Theorem \ref{sympl} is based on the following Proposition:

\begin{proposition}
\label{nondeg}We have%
\begin{equation}
d\alpha_{i}=d\alpha_{j} \label{s2}%
\end{equation}
on $\mathfrak{h}_{2,2h^{1,2}}\cap(U_{i}\cap U_{j}).$ Thus there exists a
holomorphic non degenerate form $\psi$ such that $\psi\left\vert _{U_{i}%
}\right.  =d\alpha_{i}.$
\end{proposition}

\textbf{Proof: }It is easy to see that since the extended period domain
$\mathfrak{h}_{2,2h^{1,2}}$ is an open set of a quadric in $\mathbb{P}(H^{3}%
($M,$\mathbb{Z)\otimes C)}$ then that the tangent space $T_{\tau
,\mathfrak{h}_{2,2h^{1,2}}}$ to any point $\tau\in\mathfrak{h}_{2,2h^{1,2}}$
can be identified with the orthonormal complement $\left(  H_{\tau}%
^{3,0}+\overline{H^{3,0}}\right)  ^{\perp}\subset H^{3}($M$,\mathbb{C)}$ with
respect to the metric G defined in Definition \ref{metric}. Let $\upsilon$ and
$\mu\in T_{\tau,\mathfrak{h}_{2,2h^{1,2}}}\subset H^{3}($M,$\mathbb{C)}$. From
the definition of the form $d\alpha_{i},$ it follows that
\begin{equation}
d\alpha_{i}(\upsilon,\mu)=\left\langle \upsilon,\mu\right\rangle , \label{s4}%
\end{equation}
where $\left\langle \upsilon,\mu\right\rangle $ is the symplectic form defined
by the intersection form on $H^{3}($M,$\mathbb{Z})$. From here Proposition
\ref{nondeg} follows directly. $\blacksquare$

Theorem \ref{sympl} follows directly from Proposition \ref{nondeg}.
$\blacksquare$

\begin{corollary}
\label{paral}The holomorphic two form $\psi$ is a parallel form when
restricted to $\mathcal{K\times K\subset}\mathfrak{h}_{2,2h^{1,2}}$
\textit{with respect to the CHSV connection. }(See Definition \ref{csv}.)
\end{corollary}

\textbf{Proof: }The corollary follows directly from Definition \ref{csv} and
the fact that the imaginary part of the Weil-Petersson metric when restricted
to the tangent space $T_{\tau,\mathcal{K}}=H^{2,1}\subset H^{3}($%
M$,\mathbb{C)}$ is just the restriction of intersection form $\left\langle
\upsilon,\mu\right\rangle $ on $H^{3}($M,$\mathbb{C)}$. $\blacksquare$

\section{Algebraic Integrable System on the Moduli Space of CY Manifolds.}

\begin{definition}
\label{dm}Let N be an algebraic variety. An algebraic integrable system is a
holomorphic map $\pi:X\rightarrow$N where \textbf{a.} X is a complex
symplectic manifold with holomorphic symplectic form $\psi\in\Omega^{2,0}(X)$;
\textbf{b}. The fibres of \ $\pi$ are compact Lagrangian submanifolds, hence
affine tori; \textbf{c. }There is a family of smoothly varying cohomology
classes $[\rho_{n}]\in H^{1,1}(X_{n})\cap H^{2}(X_{n},\mathbb{Z)},n\in N,$
\textit{such that }[$\rho_{n}$] \textit{is a positive polarization of the
fibre }$X_{n}$. \textit{Hence} $X_{n}$ \textit{is an abelian tersor, i.e. on
}$X_{n}$ we do not have a point which represents zero to define a structure of
a group on $X_{n}$. See \cite{DM}.
\end{definition}

This\ notion is the complex analogue of completely integrable (finite
dimensional) systems in classical mechanics was introduced by R. Donagi and E.
Markman in \cite{DM}. We will show that the family $\mathcal{A\rightarrow
M}\left(  \text{M}\right)  $ as defined in Definition \ref{abel} is an
algebraic integrable system in the sense of Donagi-Markrman.

\begin{theorem}
\label{fr}The holomorphic family $\mathcal{A\rightarrow M}\left(
\text{M}\right)  $ defines an algebraic integrable system on the moduli space
of three dimensional CY manifolds $\mathcal{M}\left(  \text{M}\right)  $ in
the sense of R. Donagi and Markman.
\end{theorem}

\textbf{Proof: }We must check properties\ \textbf{a}, \textbf{b }\ and
\textbf{c} stated in \textbf{Definition }\ref{dm}. In order to check property
\textbf{a }\ and \textbf{b, }we need to construct a non-degenerate closed
holomorphic two form $\Omega_{1}$ on the cotangent space $T^{\ast}%
\mathcal{K}\left(  \text{M}\right)  $. The cotangent space $T_{\tau}^{\ast
}\left(  \text{M}\right)  $ at a point $\tau\in\mathcal{K}\left(
\text{M}\right)  $ can be identified with $H^{1,2}($M$_{\tau})$ by contraction
with $\overline{\Omega_{\tau}},$ where $\Omega_{\tau}$ is a holomorphic three
form such that
\[
-\sqrt{-1}%
{\displaystyle\int\limits_{\text{M}}}
\Omega_{\tau}\wedge\overline{\Omega_{\tau}}=1.
\]
Then the local Torelli theorem for CY manifolds shows that the restriction of
the symplectic form $\left\langle \upsilon,\mu\right\rangle $ defined by the
intersection form on $H^{3}($M,$\mathbb{Z)}$/Tor by formulas $\left(
\ref{s1}\right)  $ and $\left(  \ref{s4}\right)  $\ will give a globally
defined holomorphic two form $\Omega_{1}$ on the cotangent bundle of
$\mathcal{K}\left(  \text{M}\right)  $. Then the properties \textbf{a }and
\textbf{b} as stated in \cite{DM} are obvious.

Next we will construct the smoothly varying cohomology classes $[\rho_{\tau}]$
which fulfill property \textbf{c. }Let $\rho(1,1)$ be the imaginary form of
the Weil-Petersson metric on $\mathcal{M}\left(  \text{M}\right)  $. Since we
proved in \cite{to89} that the potential of the Weil-Petersson metric is
defined from a metric on the relative dualizing line bundle of the family
$\mathcal{X\rightarrow M}\left(  \text{M}\right)  $, we deduce that
$\rho(1,1)$ is a smoothly varying cohomology class of type $(1,1)$. From here
we deduce that for each $\tau\in\mathcal{M}\left(  \text{M}\right)  $,
$[\rho_{\tau}]\in H^{1,1}(\mathcal{A}_{\tau},\mathbb{R)}$ and $[\rho_{\tau}]$
varies smoothly. We need to show that $[\rho_{\tau}]\in H^{2}(\mathcal{A}%
_{\tau},\mathbb{Z)}.$ This statement is equivalent to saying that if $\nu$ and
$\mu$ are any two vectors in the lattice%
\[
\Lambda_{\tau}\subset H^{2,1}(\text{M}_{\tau})\oplus H^{1,2}(\text{M}_{\tau
}),
\]
then $\rho_{\tau}(\nu,\mu)\in\mathbb{Z}$ and if $\gamma_{1},..,\gamma
_{2h^{1,2}}$ is a $\mathbb{Z}$-basis of the lattice $\Lambda_{\tau}$, then
$\det(\rho_{\tau}(\gamma_{i},\gamma_{j}))=1.$ We proved that the imaginary
part of the Weil-Petersson metric is a parallel with respect to the
Cecotti-Hitchin-Simpson-Vafa connection. (See Remark \ref{chvs}.) We used the
Cecotti-Hitchin-Simpson-Vafa parallel transport to define the $\mathbb{Z}$
structure on%
\[
T_{\tau,\mathcal{K\subset K\times K}}=H^{2,1}\left(  \text{M}\right)
+H^{1,2}\left(  \text{M}\right)  .
\]
From here it follows that the number $[\rho_{\tau}](\nu,\mu)$ is equal to the
cup product of the parallel transport of the vectors $\nu$ and $\mu$ at a
point $(\tau,\upsilon)\in\mathfrak{h}_{2,2h^{1,2}}(\mathbb{Q)},$ which is an
integer. Exactly the same arguments show that
\[
\ \det(\rho_{\tau}(\gamma_{i},\gamma_{j}))=1.
\]
So the family $\mathcal{A\rightarrow M}\left(  \text{M}\right)  $ fulfills
properties \textbf{b} and \textbf{c.} Our Theorem is proved. $\blacksquare$

\begin{corollary}
\label{fr1}On the tangent bundle of the moduli space of polarized CY
threefolds there exists a canonical Hyper-K\"{a}hler metric.
\end{corollary}

\textbf{Proof:} Cor. \ref{fr1} follows directly from \cite{f}. $\blacksquare$

\end{document}